\documentclass[nofootinbib,aip,amsmath,amssymb,reprint,twocolumn,10pt]{revtex4-2}
\usepackage{graphicx}
\usepackage{dcolumn}
\usepackage{bm}
\usepackage[utf8]{inputenc}
\usepackage[T1]{fontenc}
\usepackage{mathptmx}
\usepackage{verbatim}
\usepackage{fdsymbol}
\usepackage{graphicx}
\usepackage{sidecap}
\usepackage{placeins}
\usepackage{xcolor}
\usepackage{amsfonts,amsthm,amsmath,amssymb,amscd,graphicx}
\usepackage{mathtools}
\usepackage{hyperref}
\usepackage{amsmath}
\usepackage{amssymb}
\usepackage{bbold}
\usepackage{float}
\usepackage{hyperref}
\usepackage{tabularx}
\usepackage[ruled,lined]{algorithm2e}

\newcommand{\llbracket}{\mathopen{[\![}}
\newcommand{\rrbracket}{\mathclose{]\!]}}

\begin{document}

\title{Widespread neuronal chaos induced by slow oscillating currents}

\author{James Scully}
\email{jamesjscully@gmail.com}
\affiliation{Neuroscience Institute, Georgia State University,
100 Piedmont Ave., Atlanta, GA 30303, USA.}
\author{Carter Hinsley}
\email{chinsley1@student.gsu.edu}
\affiliation{Department of Mathematics \& Statistics, Georgia State University, 
25 Park Pl., Atlanta, GA 30303, USA.}
\author{David Bloom}
\email{bloomdt@gmail.com}
\affiliation{Neuroscience Institute, Georgia State University, 
100 Piedmont Ave., Atlanta, GA 30303, USA.}
\author{Hil G.E. Meijer}
\email{h.g.e.meijer@utwente.nl}
\affiliation{Department of Applied Mathematics, TechMed Centre, University of Twente, The Netherlands}
\author{Andrey L. Shilnikov}
\email{ashilnikov@gsu.edu}
\affiliation{Neuroscience Institute and Department of Mathematics \& Statistics, Georgia State University, \\
100 Piedmont Ave., Atlanta, GA 30303, USA.}

\date{\today}
\begin{abstract}

This paper investigates the origin and onset of chaos in a mathematical model of an individual neuron, arising from the intricate interaction between 3D fast and 2D slow dynamics governing its intrinsic currents. Central to the chaotic dynamics are multiple homoclinic connections and bifurcations of saddle equilibria and periodic orbits. This neural model reveals a rich array of codimension-2 bifurcations, including Shilnikov-Hopf, Belyakov, Bautin, and Bogdanov-Takens points, which play a pivotal role in organizing the complex bifurcation structure of the parameter space. We explore various routes to chaos occurring at the intersections of quiescent, tonic-spiking, and bursting activity regimes within this space, and provide a thorough bifurcation analysis. Despite a high dimensionality of the model, its fast-slow dynamics allow a reduction to a one-dimensional return map, accurately capturing and explaining the complex dynamics of the neural model. Our approach integrates parameter continuation analysis, newly developed symbolic techniques, and Lyapunov exponents, collectively unveiling the intricate dynamical and bifurcation structures present in the system.

\end{abstract}

\maketitle

{\bf This study investigates the intrinsic mechanisms underlying widespread chaos in a biologically plausible model of a neuron, specifically of the Hodgkin-Huxley type, driven by two slowly oscillating currents. We focus on the role of homoclinic bifurcations in triggering and sustaining chaotic behavior that extends well beyond typical expectations in neurodynamics. By developing and implementing analytical methods and simulation tools, we offer new approaches for studying high-dimensional deterministic systems with structurally unstable chaotic dynamics due to complex global bifurcations of saddle equilibria and periodic orbits. Our findings deepen the understanding of the origins and universal nature of deterministic chaos, with potential applications across a wide range of systems, including the life sciences. We aim to make these insights and tools accessible to a broad, interdisciplinary audience, from specialists to newcomers in dynamical systems and mathematical neuroscience.}

\section{Introduction}
\begin{figure}[t!]
   \begin{center}\includegraphics[width=1\linewidth]{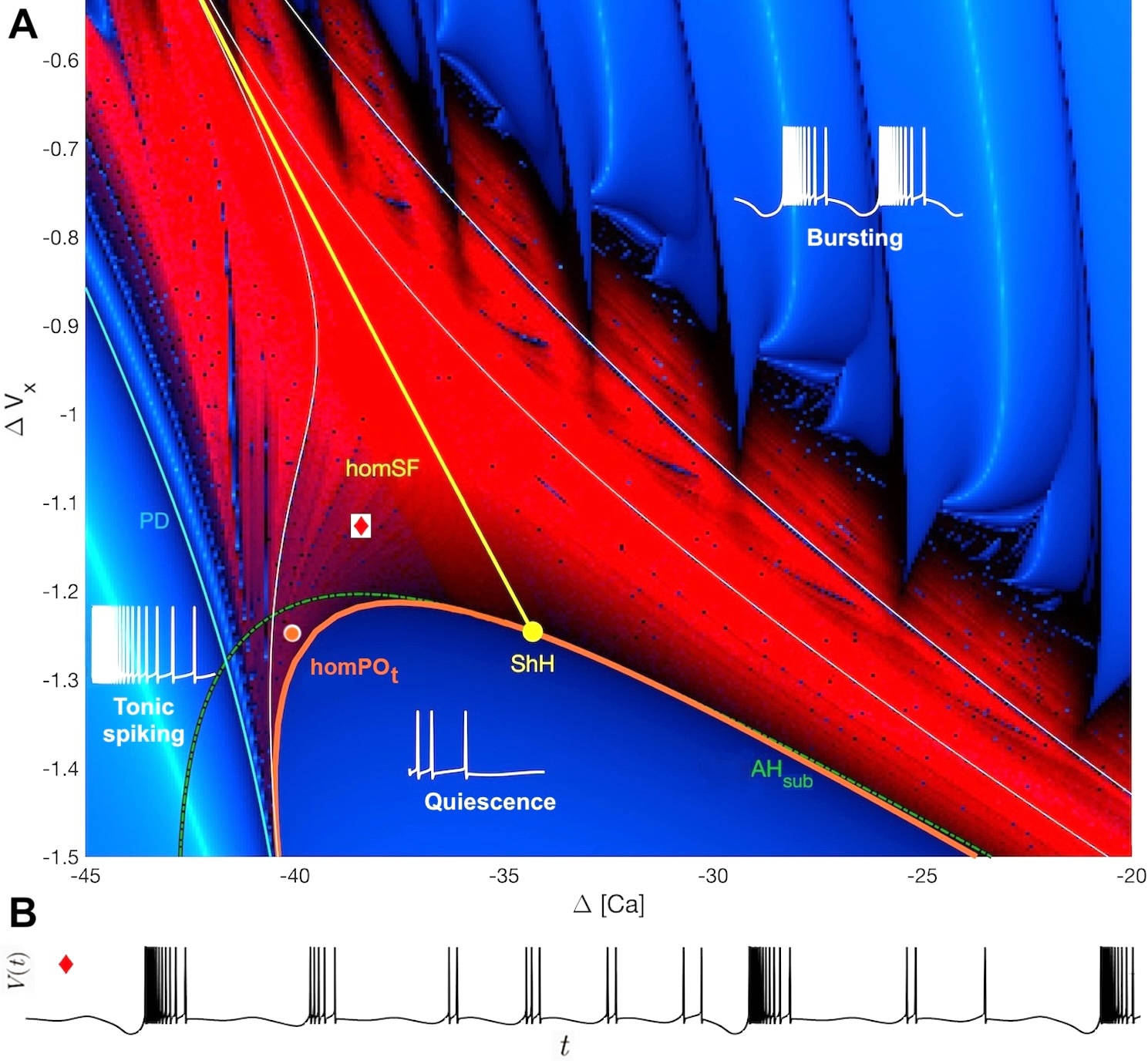}
       \caption{(A) Close-up of the widespread chaotic region (red) in a bi-parametric sweep of $\rm \Delta [Ca]$ and $\rm \Delta V_x$. The colormap illustrates the Lyapunov spectrum, with each RGB component representing a different Lyapunov exponent. The red channel corresponds to positive largest Lyapunov exponents (LLE). Green and blue both correspond to negative second Lyapunov exponents, but at different magnitudes - green indicates values closer to zero, and blue represents more negative values. Chaotic regions appear red, while non-chaotic regions are blue. Representative voltage traces are overlaid on the parameter space: tonic spiking (left), a hyperpolarized quiescent state (bottom), and regular bursting (top-right). The chaotic region is centered around the bifurcation curve $\rm homSF$ (yellow line) corresponding to a Shilnikov saddle-focus. This curve extends from the codimension-2 Shilnikov-Hopf (ShH) bifurcation point located on the subcritical Andronov-Hopf ($\rm AH_{sub}$) bifurcation line (green dashed line). The $\rm homPO_t$ curve (orange line), which separates the chaotic region from the quiescent region below, corresponds to a non-transverse homoclinic to a saddle periodic orbit (PO). In the region marked by the orange dot, situated between the $\rm AH_{\rm sub}$ curve and the $\rm homPO_t$ curve, the model exhibits bistability between chaotic and hyperpolarized quiescent attractors (see Fig.~\ref{fig15}). Thin white lines within the chaotic region indicate saddle-node bifurcations of periodic orbits, adjacent to stability windows. (B) Chaotic voltage trace corresponding to the parameter values at the red diamond in panel~A, showing the irregular bursting dynamics observed within the chaotic region.}\label{fig1}
  \end{center}
\end{figure}

Chaos in bursting neuronal systems can be compared to a kaleidoscope, where a single object fractures into intricate, shifting patterns -- each rotation revealing new layers of complexity \cite{ch1,ch2,ch3,ch4,ch5,ch6,ch7,ch8,ch9,ch10,ch11,ch12}. This analogy reflects the central focus of our study: an unusually broad region of chaos in a conductance-based neuronal model of the Hodgkin-Huxley type. In our model, the topological structure of the system amplifies chaotic dynamics, where a small seed of chaos is repeatedly fragmented and replicated, leading to widespread complexity. This amplification is driven by the underlying geometry of spike generation, which propagates and multiplies the seeds of chaos across a broad parameter space. We were intrigued by the vastness of this chaotic region, which far exceeds the narrow bands of chaos typically observed in both biologically realistic and phenomenological models during transitions, such as spike-adding. This finding challenges the prevailing view that chaotic behavior is confined to narrow parameter ranges. By studying a region where chaotic dynamics are naturally widespread and robust, we aim to uncover new insights into the mechanisms enabling flexible and adaptable oscillatory patterns, which may align more closely with the variability observed in biological neurons. In this paper, we provide a detailed exploration of this expansive chaotic region, emphasizing the diverse pathways leading to chaos across its boundaries.

One of the key features of neuronal systems, including isolated neurons in both normal and dysfunctional states, is the occurrence of self-sustained oscillations, regular or irregular, with varying recurrence times~\cite{sakurai2015phylogenetic,sakurai2016central,gamma,gatech1,gatech2}. These oscillations are typically stable in the Lyapunov sense~\cite{lyapunov}, whereas chaotic oscillations are inherently unstable~\cite{book}. Current mathematical models of neurons, however, often struggle to reproduce the desired neuronal qualities, such as variability and flexibility~\cite{Shilnikov2008a}, due to the deterministic nature of the modeling approach. Introducing noise into the model can address this limitation, particularly when the system is positioned close to semi-global bifurcations  underlying transitions between activity types such as quiescent states (stable equilibria), 
tonic spiking  and bursting oscillations (stable periodic orbits) -- where the dynamics become sensitive to small perturbations~\cite{noise1}. In contrast, far from bifurcations, mathematical models with exponentially stable (i.e., structurally stable) solutions remain robust against noise and other small perturbations, persistently exhibiting the same activity patterns~\cite{noise2,Mmo2005}.

The rigidity of most neuronal models arises from their classification as slow-fast systems with significant time-scale disparities between variables. In such models, typically consisting of two slow and one fast variable (or the reverse configuration), variability in the dynamics is only observed within narrow intervals corresponding to rapid transitions in parameter space~\cite{noise2}. Although noise can artificially widen these transition intervals, this approach remains a synthetic approximation of the fluidity observed in biological systems. The absence of well-explored neuronal models that exhibit flexible or widespread chaotic oscillations highlights the importance of our study. By developing a comprehensive framework for studying these chaotic oscillations, our work lays the groundwork for in-depth studies of other neural models that more accurately capture the adaptability and variability of neuronal activity. These insights could prove invaluable in both theoretical neuroscience and in the design of artificial systems that emulate biological flexibility.

The key to generating flexible self-oscillations in this neuron model lies in the natural widening and overlapping of spike-adding transitions as they approach a Shilnikov saddle-focus bifurcation in the parameter space. This bifurcation is well known in the dynamical systems community~\cite{LP1,LP2,LP3,sfbif,Sciheritage,LPbook17}, but less familiar in the context of neuronal modeling~\cite{feudel2000homoclinic,descroches2013}. The homoclinic bifurcation culminates in a more exotic phenomenon, the cod-2 Shilnikov-Hopf (ShH) or Belyakov-I point~\cite{belyakovshh}, where the Shilnikov homoclinic saddle-focus meets the subcritical Andronov-Hopf ($\rm AH_{sub}$) bifurcation in the parameter space of the model. Globally, the homoclinic structure interacts with the spiking manifold $\rm M_{PO}$, which expands, coils, and folds the flow, generating a topological Smale horseshoe in the phase space of the neuron model.

The interaction of the homoclinic structure with the tonic spiking manifold constitutes a rich bifurcation and dynamical ``Klondike'' in the model under investigation: its codimension-1 (cod-1) and codimension-2 (cod-2) bifurcations far exceed the typical range. These include other bifurcations of equilibria with characteristic exponents $(0,0)$ (the Bogdanov-Takens point) and $(0,\, \pm i \omega)$ (fold-Hopf due to Gavrilov and Guckenheimer), sub- and super-critical Andronov-Hopf bifurcations, the cod-2 Bautin point, multiple homoclinic bifurcations~\cite{Shilnikov2004a,Shilnikov2012}, another cod-2 Belyakov-II bifurcation describing a homoclinic saddle to saddle-focus transition~\cite{belyakovssf}, the blue-sky catastrophe~\cite{Mmo2005,Shilnikov2005,blue,blue1,blue2}, period-doubling cascades, likely a torus bifurcation with 1:2-resonance, and various hereto-dimensional cycles, among others.

To explore this wealth of dynamical phenomena, we employ a comprehensive set of topological and computational tools, including slow-fast decomposition, MATCONT parameter continuation, symbolic dynamics, conditional block entropy, templates, Lyapunov exponents, and one-dimensional return maps.

This paper is structured as follows: it begins with a description of the model, followed by theoretical background to provide context and foundational concepts. The main results start with the interpretation of three two-dimensional bifurcation diagrams, followed by illustrations of the homoclinic structure arising from the saddle-focus. Next, we illustrate the topology of the pseudo-attractor using templates. Finally, we employ one-dimensional maps to investigate the routes to chaos through the Shilnikov-Hopf bifurcation, the degeneration of spike-adding transitions, and homoclinic tangencies. The paper concludes with an in-depth discussion of these results, their implications, and future research directions.

\section{Model Description}
The model central to this study is a conductance-based neuronal model consisting of five ODEs~\cite{plant75,plant76,plant81} of the Hodgkin-Huxley type. The original purpose of this Plant model was to demonstrate the mechanism of endogenous bursting recorded in R15 cells from sea slugs of the genus {\em Aplysia}\cite{RP1985,LL88}. Over the decades, the model has been adapted for various modeling and theoretical studies due to the flexibility of its slow subsystem\cite{RL87,CCB91,butera98}.

We refer to this adapted version as the swim inter-neuron (SiN) model, as substantial changes were made, including the introduction of two new control parameters, to tailor it for studies of swim central pattern generators in two specific sea slugs~\cite{Alacam2015,pairing}. This SiN model is a slow-fast system with three fast variables governing spike generation, and two slow variables on the timescale of endogenous bursting. These slow variables, $x$ and $\rm [Ca]$, create a hysteresis loop that can be manipulated by altering their voltage response properties. This manipulation reveals a range of dynamic behaviors, including bursting, tonic-spiking, quiescence, and ultimately the onset of chaos, which is the primary focus of this paper.

\begin{figure}[t!]
  \begin{center}
  \includegraphics[width=.85\linewidth]{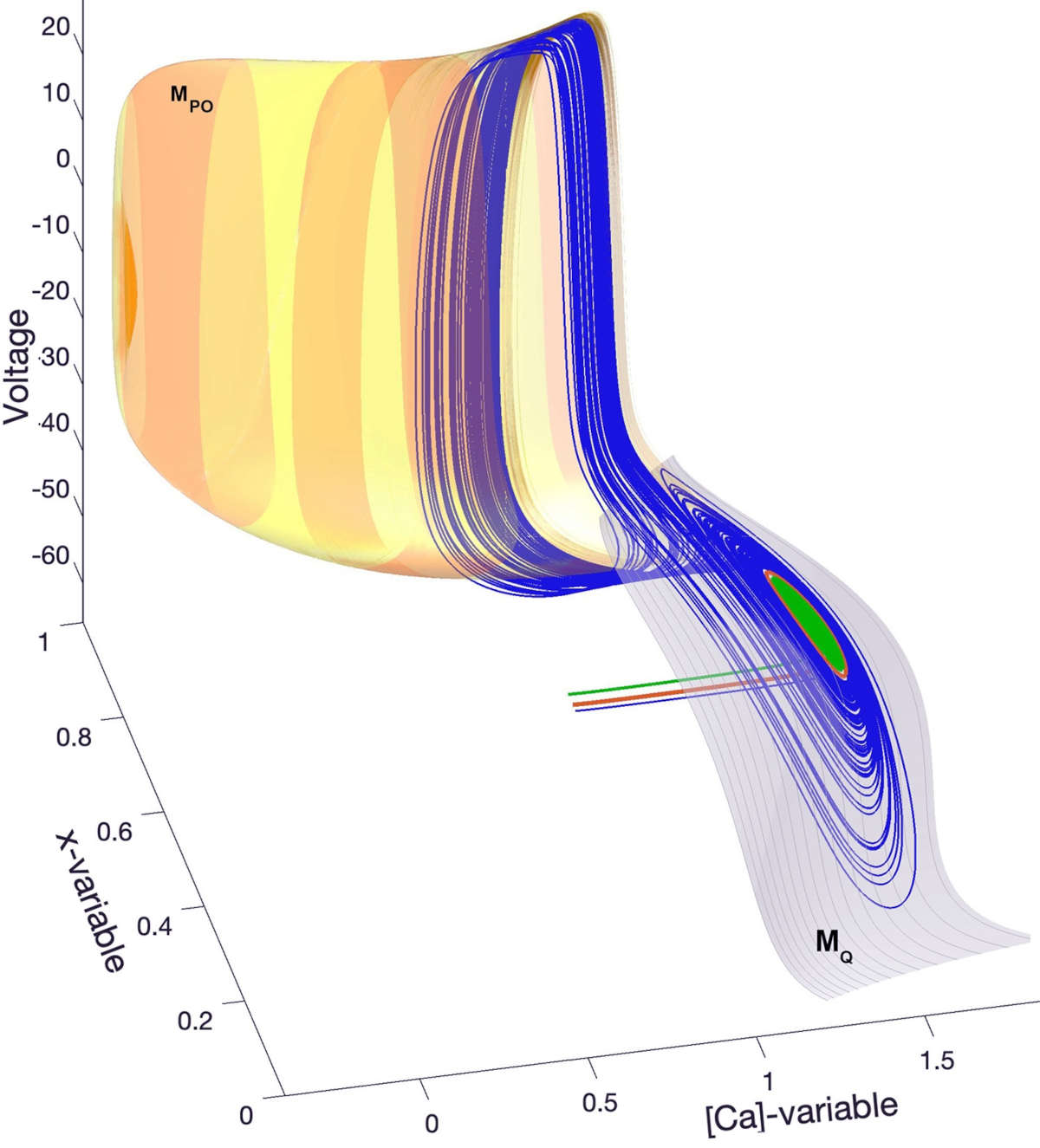}
  \end{center}
  \caption{3D $\left ( V,\, [{\rm Ca}],\, x \right )$-phase space projection of the two critical manifolds that globally determine the slow-fast dynamics in the phase space of the SiN model at $\Delta \rm [Ca]=-4.6$ and $\rm \Delta V_x = -2$. A 2D spiking manifold $\rm M_{PO}$ (yellow) is foliated by fast, round periodic orbits of the fast subsystem and features a characteristic fold. A 2D slow-motion (grey) manifold $\rm M_Q$, also referred here to as the ``dune,'' consists of equilibrium states in the fast subsystem. Three orbits are shown here: chaotic bursting (in blue) switching back and forth between $\rm M_{PO}$ and $\rm M_Q$, a green orbit converging to a stable equilibrium state, and a red saddle periodic orbit whose stable manifold locally separates two coexisting attractors, see Figs.~\ref{fig4}A and B below.}\label{fig2} 
\end{figure}

The fast subsystem, derived from the original Hodgkin-Huxley model, represents the mechanism responsible for generating fast spikes. It includes three currents: an inward sodium and calcium current, $I_{I}$, with dynamic inactivation gated by the $h(t)$ variable and instantaneous activation, $m_{\infty}(V)$; an outward potassium current, $I_K$, gated by the probability $n(t)$; and an instantaneous leak current, $I_{leak}$. We also treat the voltage $V(t)$ as a fast variable because it responds quickly to perturbations, although it plays a crucial role in the slow subsystem as well. The equations governing the fast subsystem are given by:
\begin{gather}
C_{m} {V}^\prime= - I_{I} - I_{K} - I_{leak} - I_{T} - I_{KCa}, \label{eq:voltage} ~\\
h^\prime = \frac{h_{\infty}(V)-h}{ {\tau_{h}(V)} }, \qquad n^\prime = \frac{n_{\infty}(V)-n}{\tau_{n}(V)} \label{eq:hngating}
\end{gather}
with the membrane capacitance $C_m=1$, and the currents defined as follows:
\begin{align}
&I_{I} = g_{I}\,h\, m^{3}_{\infty}(V)\,(V-E_{I}), \label{eq:icurrent}\\ 
&I_{K} = g_{K} \, n^{4}\, (V-E_{K}), \label{eq:kcurrent} \\ 
&I_{leak} = g_{L} (V-E_{L}).
\end{align}

Gradual spike frequency adaptation and post-inhibitory rebound in the SiN model are due to the slow dynamics of two currents: the TTX-resistant inward sodium and calcium current ($I_{\rm T}$) and the outward calcium-sensitive potassium current ($I_{\rm KCa}$). These are defined as:
\begin{align} 
I_{T} &= g_{T} x (V-E_{I}), \label{eq:ttxcurrent} \\
I_{KCa} &= g_{KCa}\frac{[{\rm Ca}]}{0.5+[{\rm Ca}]}(V-E_{K}), \label{eq:cacurrent} 
\end{align}
where the dynamic variables, calcium concentration $[{\rm Ca}]$ and the voltage-gated probability $x(t)$, evolve according to:
\begin{align}
x^\prime &= \frac{1}{\tau_x}\left [ \frac{1}{1+e^{-0.15(V+50- \Delta_{V_x})}}-x \right ], &\tau_x &\gg 1, \label{SlowSub1} ~\\
[{\rm Ca}]^\prime  &=  \rho \left ( K_{c}\, x\, (E_{Ca}-V + \Delta {\rm [Ca]})-[{\rm Ca}] \right  ), &\rho &\ll 1.  \label{SlowSub2}  
\end{align}

Here, $\rm \Delta V_x$ and $\Delta {\rm [Ca]}$ are bifurcation parameters introduced to control the slow dynamics of the SiN model. Our previous research provides a detailed description of the properties of the SiN model as they relate to our neuronal network modeling work~\cite{pairing}. Further information on the parameters, their biological interpretations, and the corresponding activation functions is given in Appendix III below.

\begin{figure}[t!]
  \begin{center}
  \includegraphics[width=1.0\linewidth]{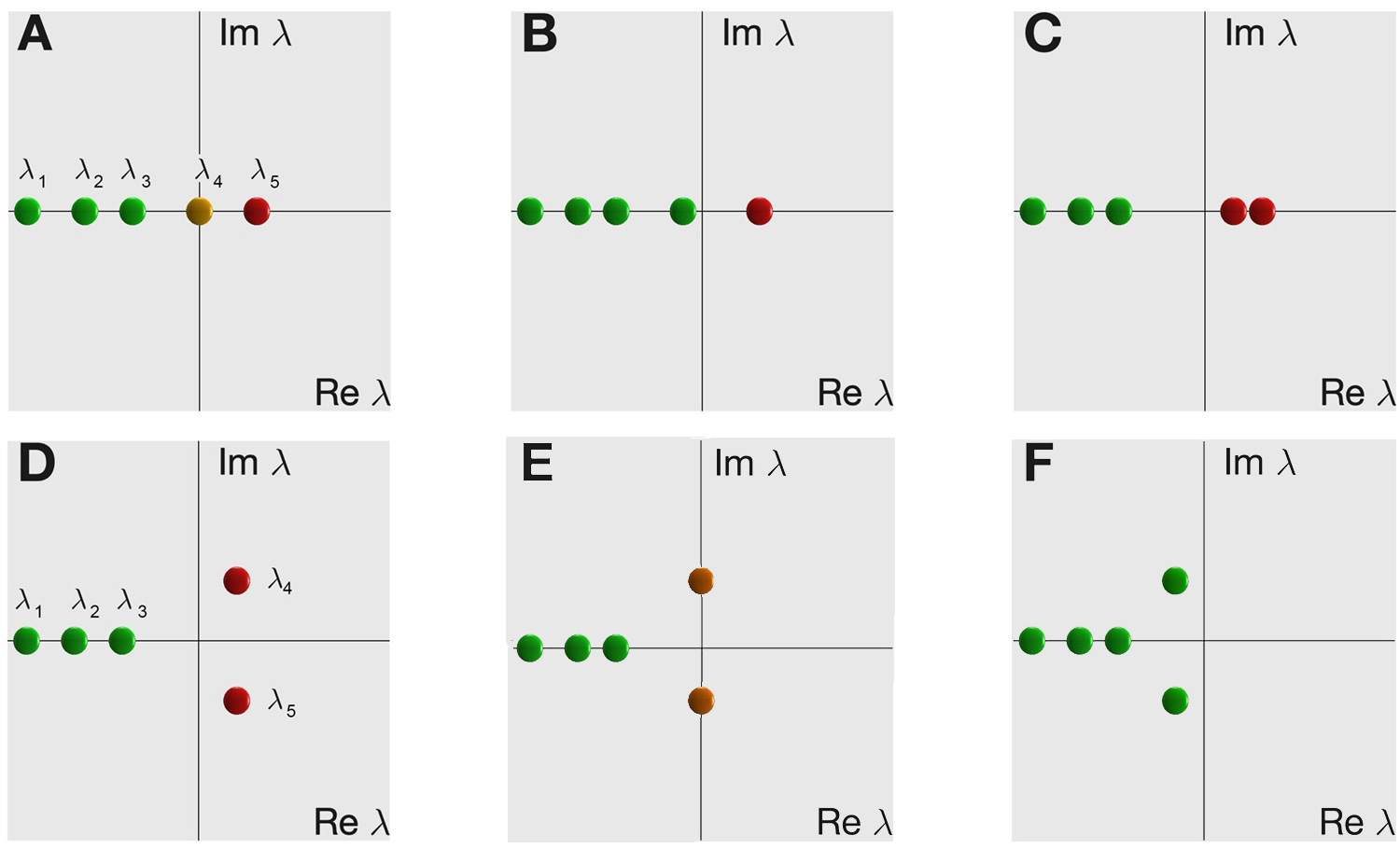}
    \caption{Qualitative stability diagrams depicting the positions of characteristic exponents in the complex plane for the two equilibrium states of the SiN model. (A) A saddle-saddle bifurcation of an equilibrium state with characteristic exponents $\lambda_{1,2,3} < \lambda_4 = 0 < \lambda_5$ gives rise to the onset of two saddles: an upper saddle of the topological (4,1)-type (B) and a lower saddle of the (3,2)-type (C) on the 2D slow-motion manifold (the dune) in the 5D phase space of the SiN model. (D) After its positive eigenvalues merge and form a complex conjugate pair next to the imaginary axis in the right open half-plane, the lower saddle transforms into the Shilnikov saddle-focus. (E-F) An Andronov-Hopf bifurcation (sub or super-critical) makes the saddle-focus a stable focus after its complex-conjugate eigenvalues cross over the over imaginary axis to the left open half-plane.}\label{fig3}
  \end{center}
\end{figure}

The first bifurcation parameter, $\rm \Delta V_x$, represents a deviation from the voltage value of $-50$ mV at which the TTX-resistant $\rm Na^{+}$-$\rm Ca^{2+}$ current becomes half-activated (see Eq.(\ref{SlowSub1})). At this point, the corresponding activation function, $x_{\infty}(V)=1/ \left (1+e^{-0.15(V+50-\Delta V_{x})} \right )$, reaches a value of $1/2$. The second bifurcation parameter, $\Delta {\rm [Ca]}$, introduced in Eq.(\ref{SlowSub2}), shifts the calcium reversal potential from its hypothetically high value of $+140$ mV, a level that cannot be experimentally validated due to the excessive current that would deplete the cell.

A segment of the $\left( \Delta {\rm [Ca]},\Delta V_x \right)$ bifurcation diagram of the SiN model is shown in Fig.~\ref{fig1}. A high-level observation reveals that the diagram can be locally divided into three peripheral regions corresponding to distinct regular dynamics: quiescent, tonic-spiking, and bursting activity. The center of the diagram is dominated by what we term informally as the ``Bermuda triangle'' of chaos, where the model exhibits irregular chaotic bursting characterized by a variable, non-repeating number of spikes in each successive burst.

The slow subsystem can be readily visualized, which is crucial for interpreting 2D and 3D (including voltage) illustrations of the slow phase space. Figure~\ref{fig2} shows a three-dimensional projection of the phase space of the SiN model, highlighting two distinct two-dimensional slow or critical manifolds. Trajectories of the model typically converge to one of these two structures, which together form the slow manifold.

For slow subsystem states where the fast subsystem contains a stable equilibrium, this equilibrium lies upon a slow manifold $\rm M_Q$ appearing like a bent sheet (shown as the grey surface in Fig.~\ref{fig2}). This slow manifold is referred to as the ``quiescent manifold'' or informally as the ``dune.'' This quiescent manifold corresponds to the region where the fast dynamics are subdued, allowing the slow variables to dominate. When the fast subsystem instead contains a stable periodic orbit (PO) representing tonic-spiking (TS) fast oscillations, the slow spiking manifold $\rm M_{PO}$ takes on a cylindrical shape (depicted with yellow and orange stripes). This cylinder is constructed by treating the state variable $\rm [Ca]$ as a parameter and allowing other state variables to converge to the stable PO. As $\rm [Ca]$ is varied, the stable spiking PO traces out $\rm M_{PO}$ in the phase space.

The dune $\rm M_Q$ wraps over the spiking manifold $\rm M_{PO}$ and repeatedly coils around before returning to its own underside, where tonic-spiking activity in the fast subsystem collapses into a hyper-polarized steady state. Bursting trajectories alternate between $\rm M_Q$ and $\rm M_{PO}$, while tonic-spiking trajectories remain as stable orbits on $\rm M_{PO}$. Quiescent trajectories and subthreshold oscillations are confined to $\rm M_Q$.

\section{Background and conceptual framework}

\begin{figure*}[t!]
  \begin{center}
\includegraphics[width=.75\linewidth]{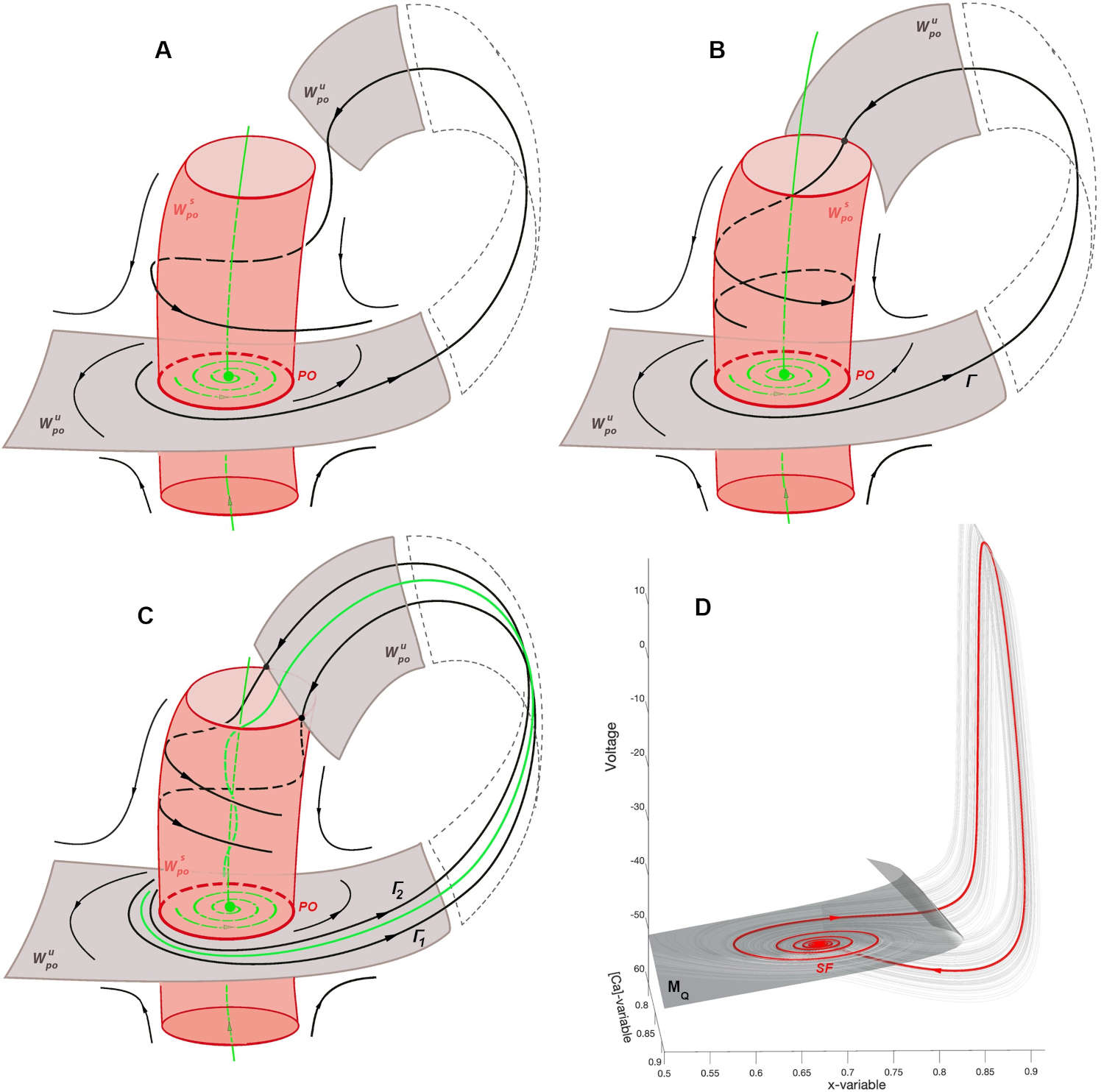}
\caption{A sketch of the unfolding of the Shilnikov-Hopf (the ShH point in Fig.~\ref{fig1}) or Belyakov-I cod-2 bifurcation in a 3D phase space, where a saddle periodic orbit (a red ring labelled by $\rm PO$) with transversely intersecting 2D stable (light-red) and unstable (grey) manifolds, $W^s_{\rm PO}$ and  $W^u_{\rm PO}$, collapses into the Shilnikov saddle-focus through a sub-critical AH-bifurcation.  The manifolds do not cross initially  in panel~A, then form a tangency along a homoclinic orbit $\Gamma$ in (B).  In (D) the manifolds cross transversally along two homoclinic obits $\Gamma_1$ and $\Gamma_2$ allowing an access to the stable equilibrium state (green dot) within the region bounded by $W_{\rm PO}^s$ (a light-red cylinder). (D) 3D phase projection depicting a Shilnikov saddle-focus with a (red) homoclinic orbit spiraling away and returning to it on the dune $\rm M_Q$ (bent grey surface) from below; this homoclinic bifurcation occurs along the corresponding curve $\rm homSF$ in the parameter plane (see Fig.~\ref{fig1} above) of the SiN model, at $\rm (\Delta [Ca], \Delta V_x)$-parameter value $(-35.98, -1.1)$.} \label{fig4}
  \end{center}
\end{figure*}

This section provides the necessary background and selectively introduces theoretical insights into the role of homoclinic bifurcations of two saddle equilibria in the onset of chaos in this neuron model. Table \ref{table1} consolidates the abbreviations for the bifurcation curves used in this study.

We first examine the local bifurcations of the equilibria. Two equilibrium states persist throughout the chaotic region, playing a central role in its dynamics. Both equilibria are key to the development and resolution of chaotic behavior in the system. These equilibria emerge simultaneously through a saddle-saddle (SS) bifurcation and will be referred to as the upper and lower saddles, with the upper saddle characterized by higher $\rm [Ca]$- and $x$-coordinates.

\subsection{Saddle-saddle bifurcation}
In a saddle-node bifurcation, one characteristic exponent is zero, and all others are either positive or negative. In contrast, a saddle-saddle bifurcation also has one zero characteristic exponent, but the remaining exponents are split between positive and negative values, placing the equilibrium in a saddle state; see Ref.~\cite{SS} for further details on this bifurcation and its non-local implications.

This saddle-saddle (SS) bifurcation occurs on the far left of the bifurcation diagram in Fig.~\ref{fig1}, corresponding to small values of the $\rm \Delta [Ca]$ parameter, and is further illustrated in the bi-parametric sweeps shown in Fig.~\ref{fig5}. The first three panels of Fig.~\ref{fig3} qualitatively illustrate the SS-bifurcation of a non-hyperbolic equilibrium state with characteristic exponents $\lambda_{1,2,3} < 0 = \lambda_{4} < \lambda_{5}$ (Fig.~\ref{fig3}A). This bifurcation produces two saddles of two different topological types.

The upper saddle has 4D stable and 1D unstable manifolds (Fig.~\ref{fig3}B), while the lower saddle has 3D stable and 2D unstable manifolds, corresponding to its characteristic exponents (Fig.~\ref{fig3}C). Note that the characteristic exponents $\lambda_{4,5}$ arise from the slow subsystem and thus remain small, located near the origin in the complex plane.

The upper saddle, with its 1D unstable manifold, marks the termination of the spiking manifold $\rm M_{PO}$, and its unstable manifold forms the outer boundary of the attractor. The homoclinic bifurcations of the upper saddle are crucial in determining the positions of spike-adding transitions in the parameter plane. The lower saddle, with its 2D unstable manifold, transitions into a saddle-focus at low $\rm \Delta V_x$-values.

The homoclinic bifurcation curve of the lower saddle, labeled by homSF, runs through the center of the chaotic region and terminates in a Shilnikov-Hopf (ShH) cod-2 bifurcation. The unfolding of this bifurcation defines the lower boundary of the chaotic region, labeled homPO$\rm_t$, as shown in the bifurcation diagram in Fig.~\ref{fig1}.

\subsection{Andronov-Hopf bifurcations and a Bautin cod-2 point}
 
Andronov-Hopf (AH) bifurcations can be classified into two types: sub-critical and super-critical. An AH bifurcation occurs when a pair of complex-conjugate characteristic exponents, denoted $\lambda_{4,5} = \alpha \pm i \omega$, cross the imaginary axis in the complex plane, as illustrated in Figs.~\ref{fig3}D-F. The type of the AH bifurcation is determined by the sign of the leading Lyapunov coefficient\cite{book}.

If the Lyapunov coefficient is positive, the AH-bifurcation is a sub-critical one. In this case, a stable focus (with characteristic exponents as shown in Fig.~\ref{fig3}F) becomes a saddle-focus with 2D unstable and 3D stable manifolds (see Fig.~\ref{fig3}D), as a saddle periodic orbit collapses into it. If the Lyapunov coefficient is negative, the AH-bifurcation is a super-critical one, through which the stable focus loses its stability and becomes a saddle-focus, accompanied by the emergence of a stable periodic orbit.

The cod-2 Bautin bifurcation occurs when the Lyapunov coefficient changes sign at the AH bifurcation (for $\left ( \Delta {\rm [Ca]},\Delta V_x \right )=(27.41,\,-2.7)$ in the SiN model), marking the transition between sub-critical and super-critical AH behaviors. The corresponding Bautin point (BP) separates two branches of the AH curve in the parameter plane: sub-critical (dashed green, AH$_{\rm sub}$) and super-critical (solid green, AH$_{\rm super}$), as shown in Figs.~\ref{fig1} and \ref{fig5}.

The unfolding of the Bautin bifurcation includes an additional curve, homSN$_{\rm PO}$, (shown in blue), corresponding to a saddle-node periodic orbit, where the saddle and stable periodic orbits merge and annihilate. The bottom three panels of Fig.~\ref{fig3} show how the characteristic exponents of the lower saddle evolve through the Andronov-Hopf bifurcation (on the green line in Fig.~\ref{fig1}). Two real positive exponents become complex conjugate (Fig.~\ref{fig3}D), cross the imaginary axis (Fig.~\ref{fig3}E), and move into the left half-plane in the complex plane (Fig.~\ref{fig3}F).

\subsection{Blue-sky catastrophe}

The blue-sky catastrophe is a saddle-node periodic orbit with a structurally stable homoclinic connection in three or more dimensions\cite{book}. The key feature of this non-local bifurcation is that it leads to the onset of a stable periodic orbit with an infinitely long period and spatial extent. This bifurcation is common in slow-fast systems\cite{Mmo2005,blue}, where it describes transitions between tonic-spiking and bursting oscillations in neuronal systems~\cite{Shilnikov2005,Shilnikov2008a}. The bifurcation occurs in the given neural model on the $\rm homSN_{PO}$ curve in the bifurcation diagrams presented in Fig.~\ref{fig5}.

\subsection{Shilnikov saddle-focus}

A homoclinic saddle-focus was discovered and analyzed by L.P. Shilnikov~\cite{LP1,LP2,LP3} in the general case for $\mathbb{R}^n$. He demonstrated that when a complex-conjugate pair of characteristic exponents, $\lambda_{4,5}=\alpha \pm i \omega$, is closest to the imaginary axis in the complex plane (as shown in Fig.~\ref{fig3}D), the neighborhood of the homoclinic orbit of the saddle-focus contains countably many periodic orbits. This phenomenon is a prerequisite for the onset of spiral chaos\cite{sfbif}, which is directly associated with the Shilnikov saddle-focus~\cite{book}.

Figure~\ref{fig4}D illustrates such a homoclinic orbit in the phase space projection of the SiN model. This bifurcation occurs along the curve labeled by homSF in the diagrams shown in Figs.~\ref{fig1} and \ref{fig5}.

\subsection{Shilnikov-Hopf or Belyakov-I point of cod-2}

The lower boundary of chaotic behavior can be understood through the unfolding of the Belyakov type-I bifurcation, also known as the Shilnikov-Hopf (ShH) bifurcation, occurring at $\left ( \Delta {\rm [Ca]},\,\Delta V_x \right )=(-34.32,\,-1.246)$. This cod-2 bifurcation occurs when a homoclinic orbit connects to a saddle-focus undergoing a subcritical Andronov-Hopf bifurcation. L.A. Belyakov \cite{belyakovshh} was the first to examine the case of a weak Shilnikov saddle-focus in $\mathbb{R}^n$, characterized by a pair of purely imaginary characteristic exponents (see Fig.~\ref{fig3}E). In his analysis, Belyakov assumed that the Andronov-Hopf bifurcation was subcritical, a key condition for this scenario. This implies that, prior to becoming a saddle-focus, the system possesses a stable focus enclosed by an unstable periodic orbit on a local two-dimensional center manifold or by a saddle periodic orbit (PO) with homoclinic trajectories in the full phase space. The key stages of the ShH unfolding are illustrated in Fig.~\ref{fig4}.

The local unfolding of the ShH point includes three cod-1 bifurcation curves: $\rm AH_{sub}$, representing the subcritical Andronov-Hopf bifurcation; $\rm homSF$, representing the homoclinic saddle-focus bifurcation; and $\rm homPO_t$ representing homoclinic tangencies to a saddle periodic orbit surrounding the stable focus. The $\rm AH$ curve divides the unfolding parameter plane, placing the homoclinic saddle-focus curve on one side and the curve of homoclinic tangencies on the other. The homoclinic orbits to periodic orbits form an ``umbrella'' shape in the parameter plane, featuring a tangency at the cod-2 ShH bifurcation point. Along this umbrella, the unstable and stable manifolds of the saddle periodic orbits form a tangent intersection, representing the boundary of possible chaos (see Fig.~\ref{fig4}B). The region between the AH curve and the umbrella of critical tangencies contains no homoclinic orbits but may still exhibit chaos, as well as stability windows (see Fig.~\ref{fig4}A).

On the opposite side of the umbrella, there are two intersections between the stable and unstable manifolds. The segment of the unstable manifold between these intersections belongs to the stable basin that emerges through the subcritical AH bifurcation (see Fig.~\ref{fig4}C). In this region, the system exhibits transient chaotic behavior, but trajectories ultimately settle into the stable focus. The bifurcation unfolding is depicted in the bifurcation diagrams in Figs.~\ref{fig1} and \ref{fig5}, and further illustrated with one-dimensional return maps in the results section.

The stages of this global bifurcation unfolding involving a saddle periodic orbit with structurally stable and tangent homoclinics in the 3D case are illustrated in Fig.~\ref{fig4}A-C. In a three-dimensional phase space, a pair of two-dimensional invariant manifolds typically either do not intersect, as depicted in panel~A, or intersect transversally, as shown in panel~C of Fig.~\ref{fig4}. In the transverse case, intersections occur along two homoclinic trajectories denoted by $\Gamma_{1,2}$ in Fig.~\ref{fig4}C, which are bi-asymptotic to the saddle PO as $t \to \pm \infty$. Between these cases, there is a moment when the manifolds become tangent along a single, non-transverse homoclinic orbit $\Gamma$, resulting from the merging of $\Gamma_{1}$ and $\Gamma_{2}$, as shown in Fig.~\ref{fig4}B.
 
The occurrence of transverse homoclinic orbits to the saddle PO implies that, nearby, there exists a hyperbolic subset. This subset includes countably many longer saddle periodic and homoclinic orbits, along with a continuum of Poisson stable trajectories~\cite{LP4,LP5,LP6}; for more information see the collections~\cite{book,Sciheritage} and references therein. When the saddle PO collapses into a stable focus at the cod-2 ShH point, it becomes a Shilnikov saddle-focus with a homoclinic orbit, as shown in the phase space projection in Fig.~\ref{fig4}D.

\subsection{Belyakov-II: saddle to saddle-focus transition of cod-2}

\begin{figure*}[t!]
  \begin{center}
  \includegraphics[width=0.67\linewidth]{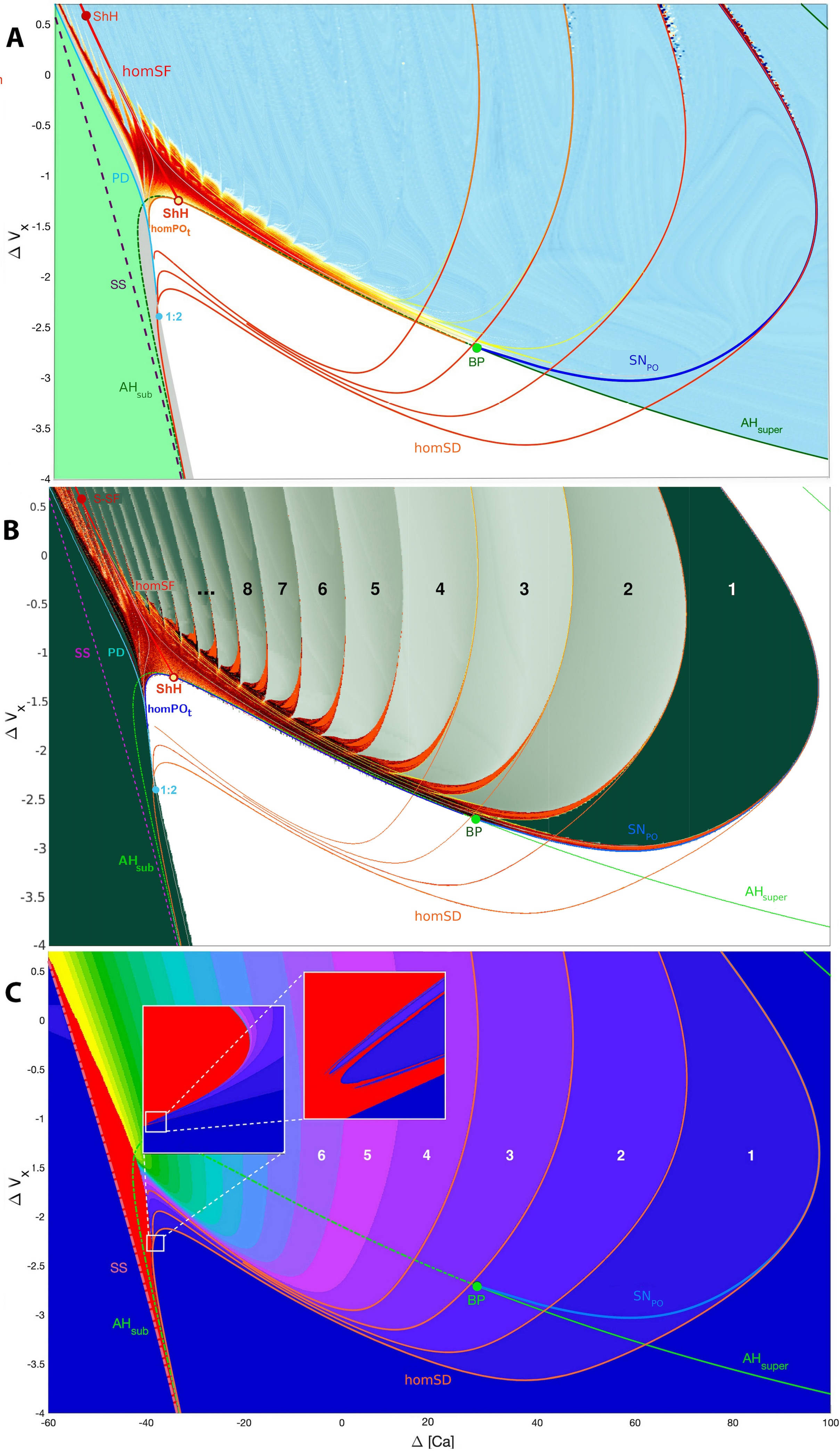} 
  \caption{Bi-parametric, $\left(\Delta [{\rm Ca}],\,\rm  \Delta V_x \right)$-sweeps of the SiN model revealing its dynamics with the aid of three different approaches. Overlaid bifurcation curves are computed using the MATCONT parameter continuation package. (A) Lyapunov exponent sweep, with intensified chaos indicated with the yellow and red colors. (B) Colormap corresponds to the inter-spike-interval (ISI) variance for periodic dynamics (green), and to Lempel-Ziv complexity for chaotic dynamics (red). (C) The symbolic sweep detects homoclinic bifurcations of the upper saddle, with inset panels illustrating the homoclinic U-shaped curves. Numbers $1, 2, \dots$ in (B)-(C) indicate number of spikes per burst in those regions.}\label{fig5}
  \end{center}
\end{figure*}

On the far side of the chaotic region along the bifurcation homSF curve, the seed of chaos begins with a saddle to saddle-focus (S--SF) transition, occurring at $\left ( \Delta {\rm [Ca]},\Delta V_x \right )=(-53.49,\,0.58)$; see Figs.~\ref{fig1} and \ref{fig5}. This cod-2 bifurcation was separately investigated by L.A.~Belyakov in Ref.~\cite{belyakovssf} for a 3D system, and to the best of our knowledge his result has not yet been generalized to the $n$-dimensional case, especially for the case when the saddle has two, not just one positive characteristic exponent, as in our system.  

The transition from saddle to saddle-focus occurs when its two positive real exponents, $\lambda_{4,5}$ merge and become a complex conjugate pair (see Figs.~\ref{fig3}C and D) on the homoclinic bifurcation curve. As long as $\lambda_{4,5}$ remain real and small, the homoclinic bifurcation of such a saddle gives rise to the emergence of a single saddle periodic orbit in the phase space. When the characteristic exponents $\lambda_{4,5}$ form a complex conjugate pair, there is a locus of bifurcation curves branching out of the cod-2 S--SF point and therefore the local dynamics become of complex nature. Belyakov showed that countably many secondary homoclinic orbits (homoclinics for short) could arise on one side of the transition, distinguished by the number of turns around the saddle-focus. Additionally, countably many saddle-node of periodic orbits emerge through this bifurcation, which give rise to the onset stability windows, which is a typical phenomenon associated with the Shilnikov saddle-focus at least in the 3D case. These two families of curves coincide at the cod-2 transition point S--SF where the characteristic exponents merge and split along the primary homoclinic curve.

In the case of the SiN model, this transition occurs in a region dominated by stable bursting orbits, so while complex dynamics may yet exist, they are not easily observable until they are revealed as the globally stable orbits collapse.

\section{Bifurcation sweeps}

The triangular-shaped region of widespread chaos in the parameter plane is located at the intersection of the bursting, quiescent, and tonic spiking regions, as shown in Fig.~\ref{fig1}. Three $\left (\Delta [{\rm Ca}],\,\rm  \Delta V_x \right )$ bifurcation diagrams presented in Fig.~\ref{fig5} illustrate the boundaries of this chaotic region, revealed through bi-parametric sweeps superimposed with key bifurcation curves continued using the parameter continuation package MATCONT~\cite{matcont,matcont1}. In these diagrams, the bursting region is at the top-right, the quiescent region is at the bottom, and the tonic-spiking region is on the left.

The first of these bifurcation diagrams, Fig.~\ref{fig5}A, represents a bi-parametric sweep using the largest Lyapunov exponents (LLEs), which are canonical indicators of chaos\cite{lyapunov}. In this diagram, a triangular-shaped region of chaos is indicated by red color, representing the most positive LLE values, while yellow indicates smaller positive LLEs. The light blue color stands for the periodic bursting region with zero LLE. The tonic-spiking and quiescent regions (where the LLE is negative) are colored light green and white, respectively, for visual clarity. We will describe and discuss selected bifurcation curves superimposed on the diagram below. 

The Lyapunov spectra from the 1000x1000-point scan in Fig.~\ref{fig5}A are calculated using the DynamicalSystems.jl package\cite{datserisdynamicalsystemsjl,diffeqjl}. Each trajectory runs for a fixed time interval $0 \leq t \leq 10^6 + 10^4$, discarding an initial transient of length $10^4$; the integration tolerance is set to $10^{-6}$ and the time step is fixed at $0.1$.

The second bifurcation diagram in Fig.~\ref{fig5}B illustrates the spike-adding process through the bursting region from right to left in various shades of green. Here, the color scheme shows the inter-spike interval (ISI) variance for non-chaotic regions, with dark green corresponding to high ISI variance and white corresponding to a lack of spikes observed after an initial transient. The sharp boundaries between shades of green fall precisely on spike adding transitions, which become blurred approaching the large chaotic region. Additional details on the route to chaos through a spike adding cascade will be provided later in the text. In the large chaotic region, the color of shading represents the Lempel-Ziv complexity~\cite{lz76} computed from a particular kneading sequence, indicating the topological complexity of trajectories. A detailed explanation of this procedure can be found in section~\ref{section:topologicalorganization}.

The third bifurcation diagram in Fig.~\ref{fig5}C shows the homoclinic bifurcation curves of the upper saddle, detected using a newly developed symbolic representation technique. The shaded colors represent the spike count associated with the unstable manifold of the saddle, with the boundaries between colors marking homoclinic bifurcations of the upper saddle in the parameter plane. Two inset panels show how these secondary curves bend into U-shaped forms\cite{rossler2020,xing21,sharkov2024}, except for the primary curve, which completes the unfolding of a cod-2 Bogdanov-Takens bifurcation point located outside of the given parameter frame. This unfolding includes three curves labeled $\rm SS$, $\rm AH$, and $\rm homSD$, which originate from the bifurcation point (see their descriptions in Table I below).

The applied symbolic technique is justified by the observation that the largest subthreshold oscillation reaches a sharp voltage maximum in the $([{\rm Ca}],\, x)$ projection of the phase space of the SiN model. This occurs as the trajectory passes near the stable manifold of the upper saddle. If the voltage maximum exceeds this point, the trajectory spikes and loops around the spiking manifold $\rm M_{PO}$, as shown in Fig.~\ref{fig2}. The largest subthreshold oscillations follow the lower branch of the unstable manifold of the saddle, while spikes (with the highest $\rm [Ca]$ value) follow the upper branch. Thus, a homoclinic orbit to the upper saddle lies between successive spike-adding events along its unstable manifold. The spike-adding cascade is illustrated by the numbers $1, 2, \dots$, in Fig.~\ref{fig5}C, corresponding to the number of spikes per periodic burst in the color-mapped regions of this bi-parametric sweep.

\begin{figure*}[t!]
  \begin{center}
    \includegraphics[width=.95\linewidth]{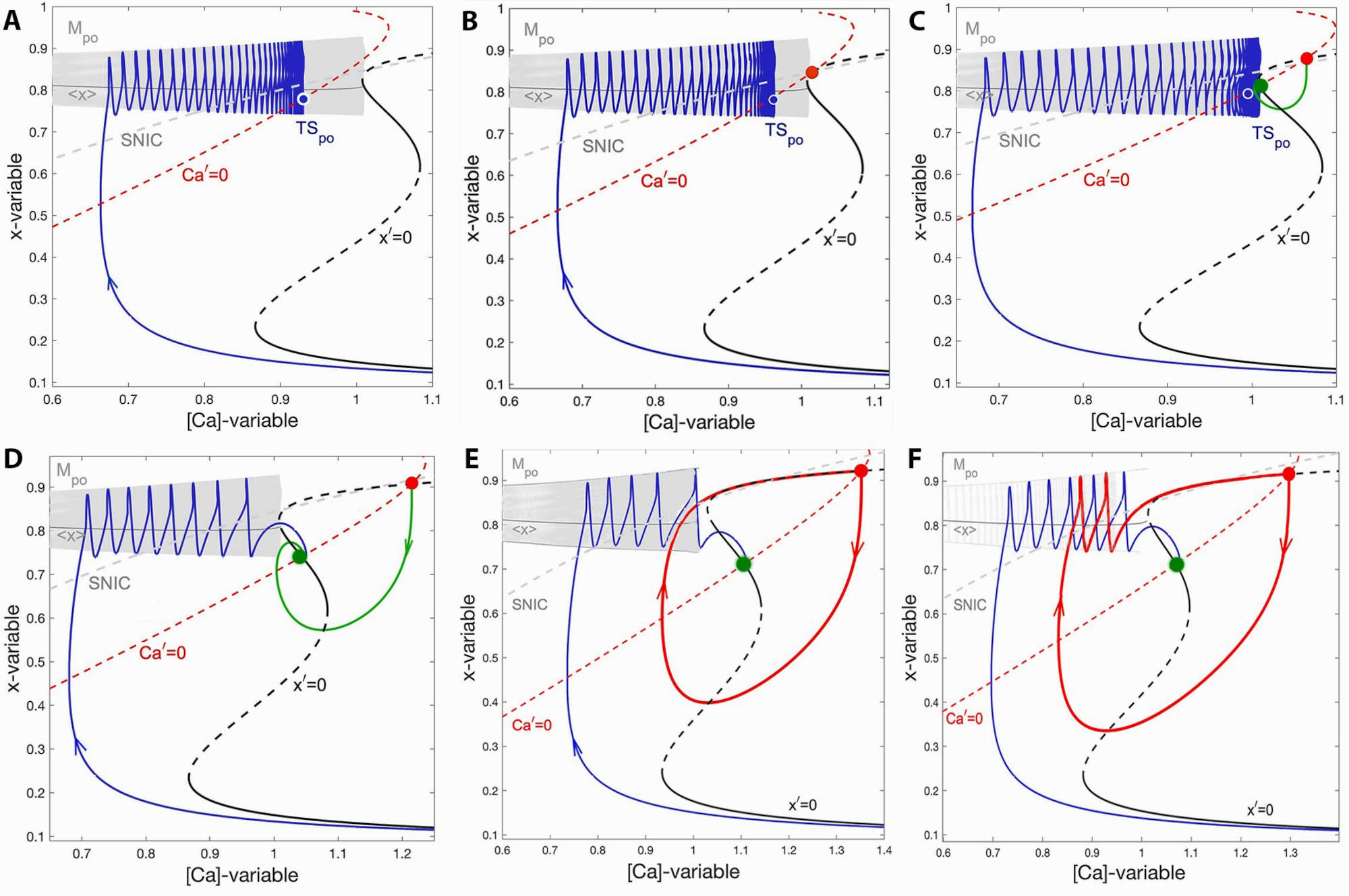}
    \caption{Six panels of the 2D ($[{\rm Ca}],\, x$)-phase projections depicting several key stages of transformative bifurcations underlying simple dynamics of the SiN model. (A) A stable tonic-spiking PO labeled as $\rm TS_{\rm po}$ on the manifold $\rm M_{PO}$ is the only attractor in the phase space at $([{\rm Ca}],\, x) = (-44.0, -2.7)$. (B) At the moment of a saddle-saddle (SS) bifurcation at $( [{\rm Ca}],\, x) = (-40.4, -2.7)$, occurring when the calcium nullcline $\rm Ca' = 0$ touches the $x$ nullcline $x' = 0$ above its top fold, two saddles emerge, with eigenvalues shown in Figs. \ref{fig3}B and C, respectively. (C) Increasing the $\rm [Ca]$ value to $-36.95$ shifts the nullcline $[{\rm Ca}]' = 0$ below the top fold on the nullcline $x' = 0$, thus making the lower saddle a stable focus through a subcritical AH bifurcation (see eigenvalues in Figs. \ref{fig3}D-F). Here, the system becomes bistable as there are two attractors: the tonic-spiking PO $\rm TS_{PO}$ and a newborn quiescent state (green dot). (D) Representative dynamics in the quiescent region of the parameter space (see Fig.~\ref{fig1}) where a typical orbit converges to a stable hyperpolarized equilibrium (green dot) at $([{\rm Ca}],\, x) = (-20.0, -2.7)$ or at $([{\rm Ca}],\, x) = (10.0, -2.7)$. (E) At $([{\rm Ca}],\, x) = (-40.122, -1.6379)$, a primary homoclinic orbit to the upper saddle (red dot) gives rise to a saddle PO. This PO collapses into the stable equilibrium (green dot) to make it a saddle-focus through a subcritical AH bifurcation above the AH curve in the parameter space. (F) A subsequent 2-spike homoclinic orbit (in red) to the upper saddle occurring at $([{\rm Ca}],\, x) = (-10.22, -2.77)$.}\label{fig6}
  \end{center}
\end{figure*}

Based on these observations, a three-step computational algorithm was applied to a regular grid over the $\rm (\Delta [Ca],\, \Delta V_x)$-parameter plane to detect the corresponding bifurcation curves:
\begin{enumerate} 
    \item Select an initial condition just beneath the saddle along its 1D unstable separatrix (in the $V$ and $x$ directions). This point is calculated from a small eigenvector associated with the positive eigenvalue. If the eigenvector points upward, it is multiplied by $-\varepsilon$ instead of $\varepsilon$ for scaling. 
    \item Evolve the trajectory and count the number of voltage maxima whenever the voltage exceeds that of the saddle.
    \item Terminate the trajectory if any voltage maximum is less than that of the saddle.
\end{enumerate}

The spike counts, corresponding to voltage maxima along the lower 1D separatrix $\Gamma_{\rm SD}^-$ of the upper saddle equilibrium $\rm SD$, are plotted as a heat map in Fig.~\ref{fig5}C. There are two sets of homoclinic orbits to the upper saddle. The first set includes orbits with an initial spike that follow the upper branch $\Gamma_{\rm SD}^+$ of the unstable manifold of $\rm SD$ up and over the spiking manifold $\rm M_{PO}$ before initially reaching the dune $\rm M_Q$. The second set consists of orbits that follow the lower branch $\Gamma_{\rm SD}^-$ of the unstable manifold of $\rm SD$ directly to the dune $\rm M_Q$. Only the trajectories from the latter set, which follow the lower branch $\Gamma_{\rm SD}^-$, are represented in Fig.~\ref{fig5}C. For some parameter values, the two branches of the unstable manifold nearly merge, while for others, they diverge significantly in the phase space.

This omission does not imply that homoclinic orbits along the upper branch of the unstable manifold are less relevant to the system. A complete analysis of the homoclinic structures of the upper saddle could provide key insights for a deeper topological investigation, but we believe such an analysis would detract from the main objectives of this study.

All three bi-parametric sweeps in Fig.~\ref{fig5} are superimposed with several selected bifurcation curves detected by the parameter continuation package MATCONT, whose descriptions are provided in Table~\ref{table1}. The SS curve gives rise to two saddle equilibria to the right of the curve. Following the SS curve downward, outside the displayed parameter range, the SS curve hits a cod-2 Bogdanov-Takens (BT) bifurcation at $\rm (\Delta [Ca],\, \Delta V_x) = (-10.2, -10.8)$. The AH curve finally folds into a cod-2 zero-Hopf point at $\rm (\Delta [Ca],\, \Delta V_x) = (-4.61, -13.8)$. 

Below the $\rm AH_{sub}$ bifurcation curve, the equilibrium state corresponding to hyperpolarized quiescence behavior becomes stable as the system transitions into a resting state. As $\rm AH_{sub}$ is followed upward, it intersects with the curve of primary homoclinics to the saddle-focus, $\rm homSF$ (red). This intersection occurs at the Shilnikov-Hopf (ShH) point, a critical center of the global unfolding. The AH bifurcation transitions from subcritical to a supercritical at the Bautin point, marking a shift in the behavior of the SiN model. To the right of the Bautin point (BP), the $\rm AH_{sup}$ (solid green) curve gives rise to subthreshold oscillations, which occur above the $\rm AH_{sup}$ curve in the parameter plane.

The origin of the first bursting trajectory, a 1-spike burster, is also associated with the Bautin point. This first spiking orbit is born through a saddle-node of periodic orbits, $\rm homSN_{PO}$ (blue), where stable and saddle POs merge and annihilate. This saddle-node has a homoclinic attachment, a feature characteristic of a blue-sky catastrophe, where stable periodic orbits of increasingly long periods emerge. The $\rm homSN_{PO}$ curve extends upward from near the Bautin point and continues beyond the top of the frames in Fig.~\ref{fig5}. A more detailed unfolding of the Bautin point and its relationship to the onset of bursting, including a cusp of periodic orbits, is explored in the 1D maps section below. Each spike-adding transition is accompanied by several additional saddle-node bifurcations of periodic orbits, as discussed in the spike-adding subsection.

While the quiescent equilibrium begins along the AH curve, this bifurcation curve alone does not bound the chaotic region in the parameter space. The stable manifold of the saddle periodic orbit (PO), which emerges from $\rm AH_{sub}$, creates a barrier that prevents chaotic trajectories from reaching the stable equilibrium. This barrier is destroyed when the unstable and stable manifolds of the saddle PO form a non-transverse homoclinic trajectory in the curve $\rm homPO_t$ (orange in Fig.~\ref{fig5}A, dark-blue in Fig.~\ref{fig5}B) in the bifurcation diagram. This homoclinic tangency forms the lower boundary of the chaotic region, originating at the cod-2 Shilnikov-Hopf (ShH) point and extending outward below the AH curve in both directions. The points on the $\rm homPO_t$ curve were identified by manually recording and analyzing homoclinic trajectories from 1D maps. More details on the computation and the association of specific homoclinics with this homoclinic tangency can be found in the 1D maps section.

The spike-adding sequence is punctuated by homoclinic orbits to the upper saddle, occurring on the curve labeled $\rm homSD$ in the $\rm (\Delta [Ca], \, \Delta V_x)$-parameter place. These homoclinic bifurcations play a central role in organizing the spike-adding process, as discussed below. Four such homoclinic bifurcations were computed using MATCONT, confirming the bifurcation diagram shown in Fig.~\ref{fig5}C.

The equilibria that arise through the saddle-saddle bifurcation (SS) are the same as those involved in the AH-bifurcation, as well as in both sets of homoclinic orbits to saddle equilibria occurring on the curves $\rm homSF$ and $\rm homSN_{PO}$ in the parameter space. The changes in the characteristic exponents of these equilibria as they split into upper and lower saddles are illustrated in Fig.~\ref{fig3}.
 
 \begin{table}[t!]
\centering
\begin{tabular}{l|l}
 Labels & Description of bifurcation curves in Fig.~\ref{fig5} \\
\hline \hline
 SS & Saddle-saddle bifurcation giving rise to 2 saddles \\
 \hline  $\rm{AH_{sub/sup}}$ & Andronov-Hopf bifurcations: sub- and super-critical \\
 \hline  BP & cod-2 Bautin point between $\rm AH_{sub}$ and $\rm AH_{sup}$ \\
 \hline  PD & First period-doubling bifurcation of a stable TS-orbit \\
 \hline  ShH & cod-2 Shilnikov-Hopf bifurcation of the saddle-focus \\
 \hline  homSF & Shilnikov bifurcation of the lower saddle-focus \\
 \hline  S--SF & Cod-2 transition between a saddle-focus and a saddle \\
 \hline $\rm homSN_{PO}$  & Saddle-node PO with a homoclinic orbit \\
 & (blue-sky catastrophe) \\
 \hline $\rm homPO_t$ & Tangent homoclinic to a saddle PO \\
 & Forms the lower boundary of the chaotic region. \\
 \hline $\rm homSD$ & Homoclinic orbit to the upper saddle
   \end{tabular}
\caption{Abbreviations used for bifurcations throughout the text.}
\label{table1}
\end{table}

\subsection{$\rm \Delta [Ca]$-routes across the bifurcation diagram}

\begin{figure}[t!]
  \begin{center}
  \includegraphics[width=0.85\linewidth]{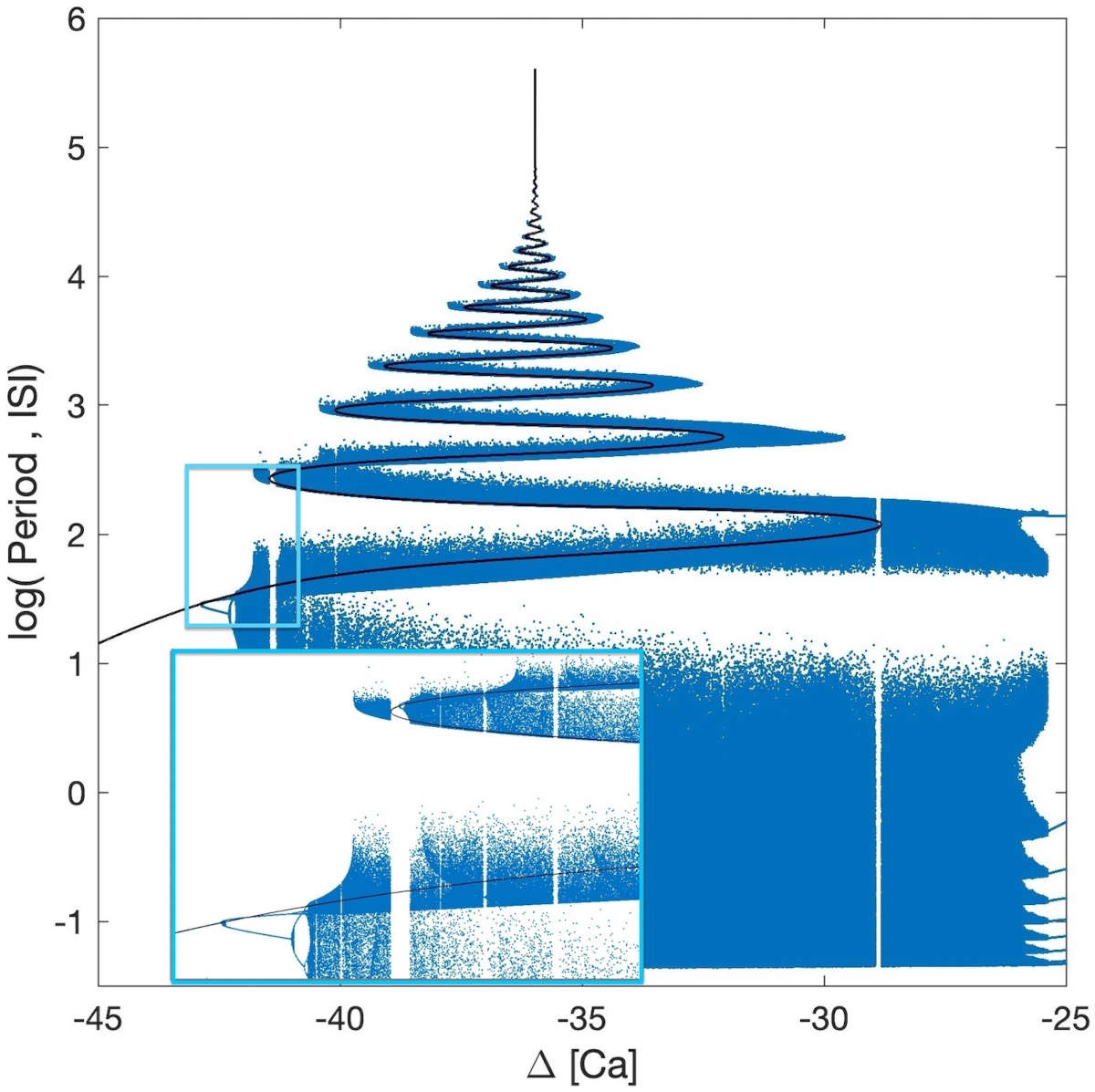}\\
    \caption{Superimposed are a one-parameter $\rm \Delta [Ca]$ sweep (in blue) of the inter-spike intervals (ISIs) in bursting voltage traces, courtesy of Dr. A.~Neiman, and the MATCONT parameter continuation (black curve), representing the period of the tonic-spiking periodic orbit on a logarithmic scale. This periodic orbit becomes a homoclinic orbit to the Shilnikov saddle-focus around $\Delta \rm [Ca]=-35.98$, at the level $\rm \Delta V_x =-1.1$, transverse to the bifurcation curve homSF in the two-parameter sweep in Fig.~\ref{fig1}. There is a close agreement between the positions of the turning points, corresponding to saddle-node bifurcations at the edges of prominent stability windows within the ISI diagram. An inset at bottom left shows a few initial stages of period-doubling bifurcations near each fold---a feature of the Shilnikov saddle-focus.}\label{fig7}
  \end{center}
\end{figure}

\begin{figure}[t!]
 \begin{center}
  \includegraphics[width=1.0\linewidth]{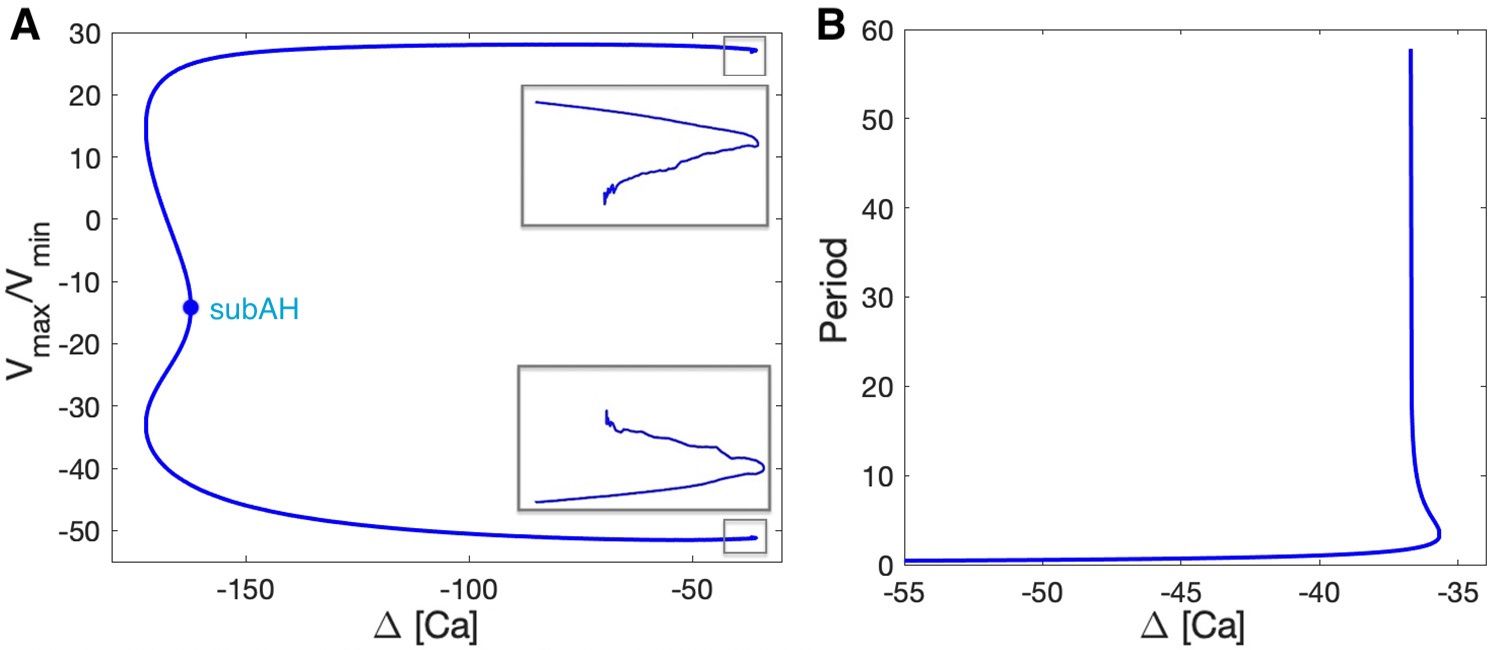}\\
    \caption{
(A) A MATCONT-based one-parameter $\Delta \rm Ca$ sweep of the spiking manifold $\rm M_{\rm PO}$ begins with a subcritical AH bifurcation of the depolarized equilibrium (EQ) in the SiN model at the level $\rm \Delta V_x=-3.5$. Next, $\rm M_{\rm PO}$ folds back with a characteristic bend corresponding to a saddle-node bifurcation of two tonic-spiking periodic orbits (TS POs), one stable and one unstable (saddle-type) (see Fig.~\ref{fig2}). Closer to $\Delta \rm Ca \simeq -37$, the outer stable section merges with an inner unstable one, characterized by a saddle PO that emerges from the primary homoclinic bifurcation (the line $\rm homSD$ in Fig.~\ref{fig5}). Emergent slow quasi-periodic oscillations around this small fold on $\rm M_{\rm PO}$ are a prerequisite for the formation of a saddle (canard) torus in the phase space of the model. (B) Period of the tonic-spiking PO is plotted against the $\Delta \rm Ca$ parameter. The turning point near $\Delta \rm Ca \simeq -37$ and the increasing period indicate, respectively, the folds on the manifold $\rm M_{\rm PO}$ shown in panel~A, and that the inner unstable PO becomes a flat homoclinic orbit to the upper saddle (see Fig.~\ref{fig6}E).}\label{fig8}
  \end{center}
\end{figure}

In this subsection, we discuss the sequences of bifurcations that occur along a three different pathways across the bifurcation diagrams in Fig.~\ref{fig5} for a fixed $\rm \Delta V_x$-parameter as $\rm \Delta [Ca]$ is increased. The first route illustrates the transition from tonic-spiking to quiescence, avoiding the bursting region in Fig. \ref{fig6}. This occurs at $\rm \Delta V_x = -2.7$. This existence of this route was the original purpose for introducing the $\rm \Delta V_x$ and $\rm \Delta [Ca]$ parameters in modeling studies. The second route at Fig \ref{fig7} shows the homoclinic structure, period-doubling bifurcations, and stability windows throughout the center of the chaotic region at $\rm \Delta V_x = -1.1$. The third route in Fig.~\ref{fig8} illustrates the possible existence of a small torus canard (of under-determined stability) at $\rm \Delta V_x = -3.5$.

Figure~\ref{fig6} shows six inset panels, each representing a snapshot of a pathway across the parameter plane with $\rm \Delta V_x = -2.7$. Each panel shows the phase space projected onto the slow variables. Across all panels, the red dashed line $\rm [Ca]' = 0$ represents the $\rm [Ca]$ nullcline, and the black dashed line $x' = 0$ is the nullcline for $x$. Blue trajectories begin nearby the lower stable branch of the x nullcline, and flow onto the tonic spiking manifold. The SNIC line (dotted grey) corresponds to a saddle-node on the invariant-circle, which is a homoclinic bifurcation to a saddle-node that separates tonic spiking from bursting in the fast subsystem. The grey shaded region shows a projection of the spiking manifold $\rm M_{\rm PO}$. This manifold was calculated by continuation of the periodic orbit in the fast sub-system with the MATCONT package. The average, or ``gravity center'' of this manifold is roughly approximated with a trajectory where the timescale of $x$ is slowed down. This approximation is illustrated with a solid grey line labelled by $\langle x \rangle$, and it can be functionally interpreted as an average $x$ nullcline. The approximation of the manifold and its average, $\langle x \rangle$, were calculated with $\rm \Delta V_x = -2.7$. For details on the calculation of the nullclines, see Appendix~I below.

Figure~\ref{fig6}A depicts a stable tonic-spiking periodic orbit labeled by $\rm TS_{PO}$ on the slow-motion manifold $\rm M_{PO}$ on the left from the SS curve. The average of this stable limit cycle is indicated by a blue dot, which lies on the $\rm [Ca]$ nullcline. This point is located near an intersection of the average $\langle x \rangle$ curve and the $\rm [Ca]$ nullcline. The deviation between the intersection and the blue dot is due to approximation error.

In Fig.~\ref{fig6}B, as the $\rm \Delta V_x$ parameter is increased, the $\rm [Ca]$ nullcline $\rm [Ca]' = 0$ tilts to the right, causing the $x$ and $\rm [Ca]$ nullclines to meet tangentially. This corresponds to a saddle-saddle bifurcation of equilibria with three negative, one positive and one zero characteristic exponents (see Fig.~\ref{fig3}A). This saddle point with a zero characteristic exponent is indicated by the red dot.

With a further increase in $\rm \Delta \rm [Ca]$, the saddle-saddle decouples into two saddles, upper and lower, with characteristic exponents as shown in Figs.~\ref{fig3}B and C, respectively. In a very small interval thereafter, the lower saddle undergoes a subcritical AH bifurcation to become a stable focus as the $\rm [Ca]$ nullcline intersects the knee point of the x nullcline between the stable and unstable branches as the complex-conjugate characteristic exponents of the saddle move leftward across the imaginary axis in the complex plane (see Figs.~\ref{fig3}D-F). 

A snapshot of this case where the lower saddle-focus becomes a stable focus is pictured in Fig. \ref{fig6}C. The SiN model thereafter becomes bi-stable within a narrow band in the parameter space, with two co-existing attractors: the stable TS PO and the hyper-polarized quiescent state. The quiescent state is illustrated with a green dot, to which the unstable manifold of the upper saddle converges, illustrated by a green trajectory.

At larger values of $\rm \Delta [Ca]$, the stable quiescent state dominates the dynamics of the SiN model (see Fig.~\ref{fig6}D) after the stable tonic-spiking periodic orbit (TS PO) disappears from the $\rm M_{PO}$ manifold, leaving only a monostable equilibrium. The mechanism behind the termination of the tonic-spiking orbit is discussed in detail below, along the sweep at $\rm \Delta V_x = -3.5$.

On the right side of the quiescent region, the pathway at the level $\rm \Delta V_x = -2.7$  through the parameter plane intersects a series of homoclinic curves to the upper saddle, $\rm homSD$. The trajectories that emerge from these homoclinic bifurcations ultimately become bursting patterns (see also the spike-adding 1D maps below). Each homoclinic orbit has a distinct number of spikes in its separatrix loop. The 0-spike and 2-spike cases are illustrated in Figs.~\ref{fig6}E and F, respectively, with red trajectories indicating these unstable homoclinic orbits. A detailed visualization of this fractal homoclinic structure in parameter space is provided in Fig.~\ref{fig5}C.

During the breakdown of tonic spiking, a critical issue arises with the localization of the slow-motion manifold when the tonic-spiking orbit that defines it loses normal hyperbolicity, i.e., when it approaches a bifurcation. This disrupts the entire framework of slow-fast dissection, where slow-motion or critical manifolds are approximated by sweeping through attractors of the fast subsystem while treating the slow variables as parameters. According to the slow-fast dissection, the manifold $\rm M_{PO}$ should terminate at the top knee point where the nullcline $x' = 0$ intersects the SNIC. However, this does not hold true in the full system, which is why the slow-fast approach breaks down near bifurcations of the stable tonic-spiking orbit. Such bifurcations typically occur when one or more Floquet multipliers of the periodic orbit cross the unit circle outward, causing convergence to or divergence from the manifold to slow down, aligning with the timescale of the slow subsystem. This can result from various bifurcations: a saddle-node bifurcation with a multiplier of +1, a period-doubling bifurcation with a multiplier of -1, or a torus bifurcation where a pair of multipliers crosses the unit circle. Additionally, this loss of normal hyperbolicity can occur when a periodic orbit approaches a homoclinic bifurcation to a saddle or saddle-node~\cite{book}. 

Following the $\rm \Delta [Ca]$-sweep at the pathway $\rm \Delta V_x = -1.1$, the stable tonic-spiking (TS) orbit loses its stability through a period-doubling bifurcation cascade as the $\rm \Delta [Ca]$ parameter increases. The resulting chaotic pseudo-attractor occupies the red triangular-shaped region, interspersed with numerous stability windows, as shown in Figs.~\ref{fig1} and \ref{fig5}A,B. Figure~\ref{fig7} illustrates how the original TS orbit transitions into a homoclinic orbit of the saddle-focus along the homSF curve in the bifurcation diagram. This transition suggests that the terminal orbits, which initially converge to the cylindrical manifold $\rm M_{PO}$, begin to resemble a ``treble clef'' as they make a large global turn followed by smaller oscillations around the saddle-focus.

Figure~\ref{fig7} superimposes a continuation diagram (in black) of the original TS orbit, in which its period is plotted against the continuation parameter $\rm \Delta [Ca]$, atop a uni-parametric $\rm \Delta [Ca]$-sweep of inter-spike intervals (ISIs) (in blue). Both representations align exceptionally well, showing that the period of the tonic-spiking periodic orbit increases logarithmically ($-\log \varepsilon$, $0 < \varepsilon \ll 1$) as the homoclinic bifurcation is approached, as indicated by the large periods in the middle of the panel. The ``stairway to heaven'' pattern, characterized by the tightening undulations of both curves, is a hallmark of the homoclinic saddle-focus bifurcation \cite{book,gl_sparrow,huber}. At each switchback, where the curves reverse in the $\rm \Delta [Ca]$ direction, a fold occurs, resulting in the emergence of a stable and unstable pair of periodic orbits. As the $\rm \Delta [Ca]$-parameter increases, each successive fold introduces orbits with an additional turn around the saddle-focus, accumulating to an infinite number of turns as its period becomes arbitrarily large near the homSF curve in the parameter space. For long periodic trajectories, the distribution of ISIs will be finite after discarding transients. For chaotic trajectories, however, the distribution of ISIs becomes dense.

The inset panel in Fig.~\ref{fig7} provides a close-up view of the first two folds, each followed by miniature period-doubling cascades. The 1D sweep also reveals stability windows (white stripes) adjacent to saddle-node (SN) bifurcations, cutting through all rungs of the ladder. These stability windows are characteristic of the period-doubling route to chaos and are visible as blue strips within the red chaotic region in the bi-parameter sweep in Fig.~\ref{fig1}, and as white strips in the sweeps shown in Figs.~\ref{fig5}A and B.

The construction of the sweep in Fig.~\ref{fig7} is based on measuring the inter-spike intervals (ISI) of a long trajectory. At each value of $\rm \Delta \rm [Ca]$, a single long trajectory is run, and the ISIs are recorded. The final state of the current trajectory is used as the initial condition for the trajectory at the next $\rm \Delta \rm [Ca]$ value. The $\rm \Delta [Ca]$-sweep at $\rm \Delta V_x = -3.5$ illustrates the fate of the tonic spiking orbit below the 2:1 resonance point with two-parameter continuation curves in Fig.~\ref{fig8}.


Note that $\rm \Delta [Ca]$-sweeps at larger $\rm \Delta V_x$ values, such as at the level $\rm \Delta V_x = -1.1$ throughout the chaotic region, all demonstrate that the stability loss and termination of the tonic-spiking orbit on the manifold $\rm M_{\rm PO}$ typically involves a cascade of period-doubling bifurcation (see Fig.~\ref{fig7}). However, at certain values of the $\Delta V_x$ parameter, the tonic-spiking orbit undergoes a 2:1 resonance transition, where the period-doubling (PD) curve ends in Fig.~\ref{fig5} at $\rm \Delta V_x = -2.6$. This cod-2 point was identified using MATCONT. The 2:1 resonance arises from a strongly reciprocal interaction between the tonic-spiking orbit on the spiking manifold $\rm M_{PO}$ and the unstable limit cycle emerging through the $\rm AH_{sub}$-bifurcation on the underlying quiescent manifold $\rm M_{Q}$. There are several possibilities for the unfolding of this 2:1 resonance point, where the TS orbit has a pair of Floquet multipliers equal to -1. One scenario involves the emergence of an invariant torus, which may be stable or unstable\cite{book}. While we present some evidence suggesting the possible existence of such a torus, a thorough investigation was not pursued, as any invariant tori do not appear to play a significant role in the onset of chaos in the SiN model.

The $\rm \Delta [Ca]$ sweep at $\rm \Delta V_x = -2.7$, discussed above, occurs below the 2:1 resonance point on the curve PD. To provide more details, we specifically investigate the fate of the tonic-spiking orbit for this case in a $\rm \Delta [Ca]$ sweep at the lower pathway $\rm \Delta V_x = -3.5$ in Fig.~\ref{fig8}. Its panel~A shows the maximum and minimum voltage values of the tonic-spiking orbit as $\rm \Delta [Ca]$ is varied. At low $\rm \Delta [Ca]$-values (far to the left of the parameter ranges shown in the figures of this paper), the tonic-spiking orbit is born through a saddle-node of periodic orbits (SNPO). This appears as a fold on the left side of panel~A. The outer stable orbit becomes the tonic-spiking orbit, while the inner unstable orbit contracts and terminates in a subcritical AH bifurcation. Note that this AH bifurcation occurs on a completely different orbit from the one marked by the $\rm AH_{\rm sub}$ lines in Figs.~\ref{fig1} and \ref{fig5}, and occurs at different $\rm \Delta [Ca]$ values.

On the right side, the continuation of the tonic-spiking orbit turns back on itself before terminating. The tiny fold, visible in the zoomed-in inset panels, may suggest the existence of an invariant torus. The unstable inner branch represents an unstable periodic orbit, which emerges through this terminating fold on the far right of Fig.~\ref{fig8}A.

Figure~\ref{fig8}B shows that the period of this orbit can grow arbitrarily large after the fold, suggesting that the unstable orbit eventually becomes a homoclinic trajectory. The numerical difficulty of continuing the unstable orbit beyond the final fold is evident from the roughness of the continuation curve in panel~A.

\subsection{Saddle-focus homoclinics near homSF}
\begin{figure}[t!]
  \begin{center}
  \includegraphics[width=1\linewidth]{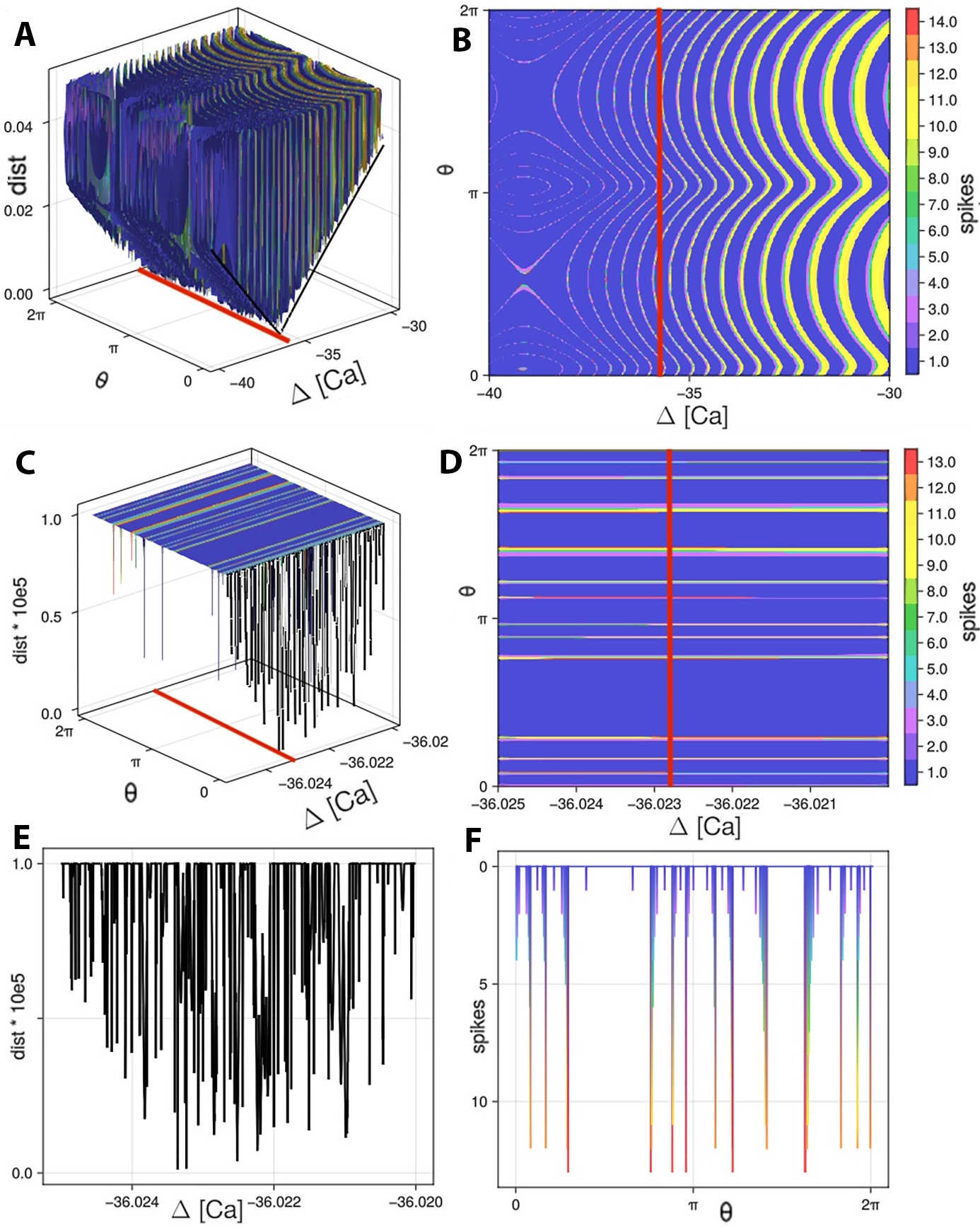}\\
    \caption{A variety of homoclinic orbits with different numbers of spikes originate from different positions on the unstable manifold of the saddle-focus. The 2D $\rm (\Delta [Ca],\, \theta)$ sweep is conducted at $\rm \Delta V_x = -1.1$, where the angular variable $\theta$ parameterizes a small circle around the saddle-focus. Trajectories flow for a single global turn, and subsequently their distances to the saddle-focus in both the $\rm [Ca]$ and $x$ directions are measured. Panel~A shows a sharp V-shape where the distance is minimized at the apex. The corresponding spike counts are more clearly seen from the top view in panel~B. Panels C and D show a zoomed-in view of the apex. Thick red lines mark the values where homoclinic bifurcations occur. The black curves show the minimum distance over all $\theta$ values, and panel~E displays this information as well. Panel~F shows that the apex of the V-shape comprises a variety of different spike numbers, corresponding to the color bars to the right.}\label{fig9}
  \end{center}
\end{figure}

As previously noted, the curve of the primary homoclinic orbit to the Shilnikov saddle-focus cuts through the center of the chaotic region. Continuation analysis has identified this curve, shown in yellow in the parameter diagram in Fig.~\ref{fig1} and in red in Fig.~\ref{fig5}. The curve originates at the cod-2 ShH point and extends upward through the cod-2 S--SF point, where the saddle-focus transforms into a saddle with all real characteristic exponents (see Figs.~\ref{fig3}C and D), continuing beyond the depicted frame. In phase space, the homoclinic orbit (see Fig.~\ref{fig4}D) forms when the saddle-focus intersects the reinsertion loop, where both slow-motion manifolds, $\rm M_Q$ and $\rm M_{PO}$, merge. To the right of the homoclinic curve, the saddle-focus rests on the slow-motion dune, while to the left, it resides on the fast spiking manifold $\rm M_{PO}$.

In the SiN model, homoclinic orbits with different numbers of spikes can exist along distinct trajectories of the unstable manifold of the saddle-focus. These homoclinics are not secondary (or tertiary, etc.) in the conventional sense\cite{book}, as they only return to a neighborhood of the stable manifold once. Instead, these homoclinic orbits are tightly clustered due to the strong contraction toward the spiking manifold, a global feature of the flow. Each primary homoclinic orbit induces secondary and higher-order homoclinics, with varying numbers of turns around the unstable manifold of the saddle-focus after each global excursion. These homoclinic structures can interact significantly and even involve multiple bursts. For instance, a secondary homoclinic orbit may exhibit three spikes during its first global excursion, followed by a large number of turns, say seven, around the saddle-focus, and finally four spikes on its second global excursion before returning to the saddle-focus.

Figure~\ref{fig9} presents a numerical investigation into the nature of these multiple homoclinics in the parameter plane. A one-dimensional cut is taken through the parameter space with the $\rm \Delta V_{x}$ parameter fixed at -1.1, while varying the $\Delta {\rm{Ca}}$ parameter. Since the unstable manifold of the saddle-focus is two-dimensional, initial conditions are sampled from a small circle around the saddle-focus. This circle, $\mathbf{x}_0$, is constructed using the unstable eigenvectors $\mathbf{u}_1$ and $\mathbf{u}_2$, and is parameterized by the variable $\theta$ as follows
\begin{equation}
	{\mathbf x}_0(\theta) = (\epsilon ({\mathbf u_1}\,\cos(\theta) + \mathbf{u}_2\,\sin(\theta)).
\end{equation}

Next, a 1000x1000 grid is sampled over the $\Delta {\rm{Ca}}$-parameter and the angular variable $0 < \theta < 2\pi$. Each initial condition is numerically integrated and events are captured at voltage maxima. Voltage maxima above the value of the upper saddle are classified as spikes, while those below are considered subthreshold oscillations. To focus on primary homoclinics, the trajectory is allowed to continue until at least one spike occurs, excluding subthreshold oscillations within the dune. Once a spike has occurred, the trajectory is terminated at the next subthreshold oscillation. The number of spikes and the $L^2$-distance, $(x^2+\rm{[Ca]}^2)^{1/2}$, from the saddle-focus in the $x$- and $\rm [Ca]$ coordinates are recorded. For visual clarity, distances in Fig.~\ref{fig9} are clipped at a maximum value.

The results from this scan reveal a single ``V''-shaped pattern in $\rm \Delta [Ca]$ (Fig.~\ref{fig9}A), indicating that the homoclinic bifurcation curves corresponding to different numbers of spikes are packed extremely closely together if they differ at all. Figure~\ref{fig9}B illustrates how the spike count of the first burst varies with $\theta$- and $\rm \Delta [Ca]$ parameters.

To further investigate this peak, we repeated the scan over successively smaller ranges of $\rm \Delta [Ca]$ values around the minimum. Figures~\ref{fig9}C-F show the results of the high-resolution scan with 2000x2000 points. Figure~\ref{fig9}E provides a clear view of the minimum across all $\theta$ values at each $\Delta {\rm{Ca}}$ value, where a singular minimum is still observed. These figures demonstrate that nearly homoclinic trajectories exhibit a full range of spike counts depending on the initial $\theta$ value. This investigation confirms that homoclinic orbits with different spike counts exist, and that their bifurcation curves are tightly packed together in parameter space. An orbit returning to the neighborhood of the saddle-focus can, in principle, pass through any $\theta$ value with a small perturbation as it exits along the unstable manifold. This underscores the complexity of the higher-order homoclinic structure, as the saddle-focus interacts globally with the tight coiling of trajectories around the spiking manifold.

It is surprising for a saddle-focus to produce multiple primary homoclinic tangencies at the same parameter value, as this implies that multiple tangencies exist between its unstable and stable manifolds. Typically, if the manifolds of tangency move independently in parameter space, this phenomenon would only occur at points of high codimension. We hypothesize that, due to the coiled fast-slow nature of the inertial manifold, tangencies may be packed together extremely closely. This is difficult to demonstrate in numerical simulations, especially in such a fast-slow system. 

\section{Topological organization of the SiN model}\label{section:topologicalorganization}

\subsection{Attractor structure}

Although chaos is known to exist near the Shilnikov saddle-focus homoclinic, the relationship between this localized chaos and the broader chaotic behavior across the parameter region remains unclear. To better understand the chaotic behavior across the entire region, we examine the topology of the model by constructing a topological template.
This involves projecting the chaotic attractor onto a two-dimensional branched manifold called a \emph{template}.
The template serves as a powerful analytical tool, enabling us to apply techniques from symbolic dynamics to the study of complexity in the SiN model.
In order for the Birman-Williams projection\cite{birmanwilliams} of the SiN system to project the phase space onto a two-dimensional template, the Hausdorff dimension of the $\omega$-limit set (i.e., the attractor) must be less than three.
To verify the applicability of the template approach, we calculated the Lyapunov spectrum of trajectories on the attractor at many parameter values using ChaosTools.jl.\cite{datserisdynamicalsystemsjl}

Let $\lambda_1, \lambda_2, \ldots, \lambda_5$ represent the Lyapunov spectrum.
In every case we examined, the Lyapunov dimension
\begin{equation}
	\sup_{u \in K}\ \dim_L(\varphi^t, u) = \sup_{u \in K}\ 2 + \frac{\lambda_1(u)}{|\lambda_3(u)|}
\end{equation}
of the flow $\varphi^t$ with respect to the globally attracting set K is less than three, typically only slightly above two due to the extreme timescale separation between slow and fast subsystems.
In accordance with the Kaplan-Yorke conjecture\cite{kaplanyorke}, the attractor can therefore be expected to embed in a three-dimensional manifold, the attractor itself having Hausdorff dimension less than three.
Carrying out numerical simulations of many trajectories near the attractor, we find that a projection to the three variables, $[{\rm Ca}],\,x$ and $V$, sufficiently approximates the embedding to straightforwardly deduce the topological structure of the attractor.
Our motivation for choosing this projection is that the slow dynamics in the $\rm [Ca]$ and $x$ variables exhibit expansion associated with the outward-spiraling saddle-focus equilibrium, while the voltage $V$-value is of physiological interest and is empirically observable in physical cells.
Selecting a three-dimensional inertial manifold projection of phase space follows the method described by Gilmore and Lefranc for studying strange attractors in dissipative systems with Lyapunov dimension less then three~\cite{gilmorelefranc}, avoiding the need to carry out the intractable task of explicitly performing the Birman-Williams projection.

\begin{figure}[t!]
  \begin{center}\includegraphics[width=1.0\linewidth]{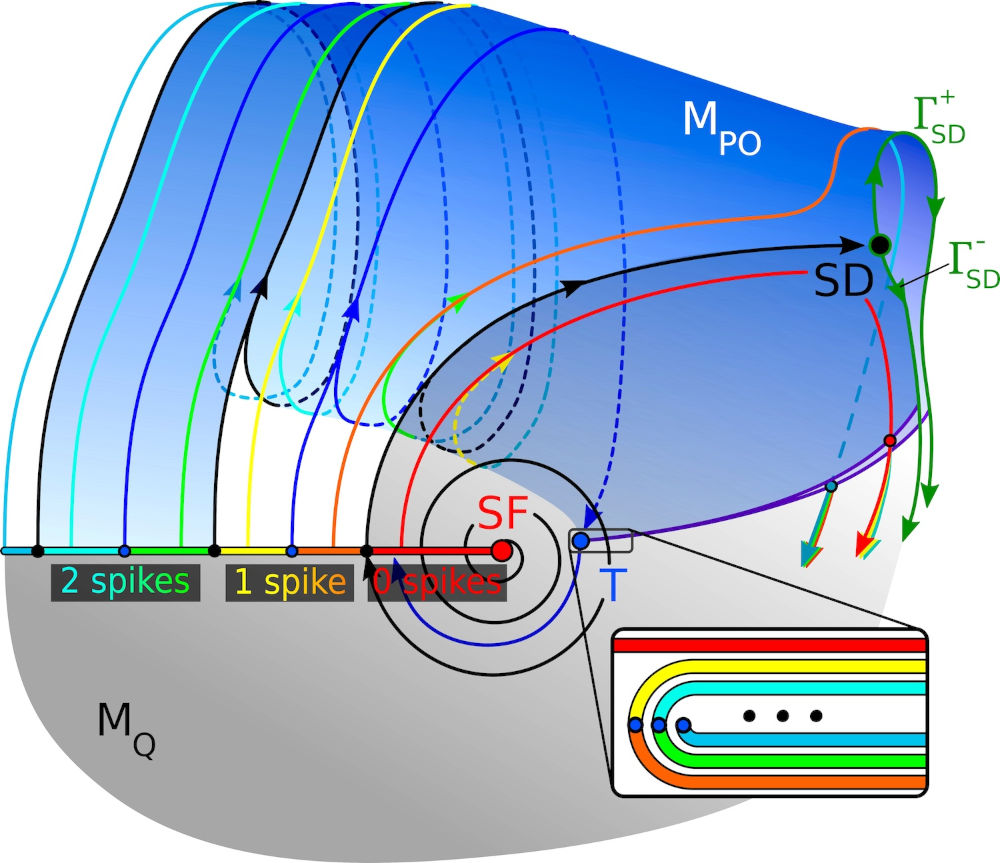}
    \caption{
    An illustration of how layering occurs in phase space. At the bottom left, a line on the dune $\rm M_Q$ is indicated, terminating at the saddle-focus (SF) marked by a red dot. Each color on the interval corresponds to a different layer, after a single pass through $\rm M_{PO}$. The final layering is shown via the nested Smale horseshoe at the bottom right. The blue surface represents the trajectories that spike twice as they pass through $\rm M_{PO}$ before returning to $\rm M_Q$. The heteroclinic connection from the unstable manifold of the saddle-focus (red dot SF) to the upper saddle (black dot SD) is shown as a black line. The reinsertion loop, where the rainbow section returns to the dune, is colored purple and spans from the unstable manifold $\Gamma_{\rm SD}^- \cup \Gamma_{\rm SD}^+$ of the upper saddle $\rm SD$ to its furthest extent at the point marked $\rm T$ in blue. $\rm T$ is the return point on $\rm M_Q$ closest to SF. This point occurs at the crease of the spike-count bands, marked by blue dots.
The flow onto each layer is highlighted by single trajectories with corresponding coloring. These trajectories are chosen to overlap as they coil onto $\rm M_{PO}$, illustrating how the folding occurs.
}\label{fig10}
  \end{center}
\end{figure}

Figure~\ref{fig10} illustrates the general structure of the chaotic attractor in a 3D phase space.
To arrive at this understanding of the topological structure of the phase space, we iteratively refined this illustration by the use of four tools:
\begin{itemize}
	\item MATCONT was used to locate equilibria, periodic orbits, homoclinic structures, and bifurcations of each.
	\item Trajectories were numerically integrated via DifferentialEquations.jl\cite{differentialequationsjl} to interactively visualize flows in the phase space.
	\item 1-dimensional return maps were calculated to study the variety of bursting behaviors exhibited on the attractor.
	\item Paper models were constructed by hand to understand how spiking trajectories layer atop one another.
\end{itemize}

The horizontal direction corresponds to the $\rm [Ca]$ variable, while the vertical direction corresponds to the $x$ variable.
At the bottom, trajectories on the dune $\rm M_Q$ flow slowly and laminarly clockwise around the saddle-focus equilibrium $\rm SF$.
At the top, trajectories rapidly wind around the topological cylinder $\rm M_{PO}$ as they gradually progress to the right.
The slow-motion manifolds $\rm M_Q$ and $\rm M_{PO}$ come together in two locations: at left, the rainbow line labeled with spike counts greater than $0$ constitutes the region of transition from quiescent or refractory behavior on $\rm M_Q$ to fast spiking activity on the spiking manifold $\rm M_{PO}$; at right, the purple ``reinsertion loop'' constitutes the cessation of spiking activity on $\rm M_{PO}$ and a transition to slow subthreshold activity on $\rm M_Q$.
Note that only those curves with arrowheads are trajectories in the phase space; for example, the reinsertion loop is not a trajectory.

The farther from the saddle-focus $\rm SF$ a trajectory is when it transitions from $\rm M_Q$ to $\rm M_{PO}$, the greater the number of spikes observed in a burst---i.e., the number of times the trajectory winds around $\rm M_{PO}$ before returning to $\rm M_Q$ via the reinjection loop.
The rainbow line at left in Fig.~\ref{fig10} is partitioned by segments of roughly equal length, each corresponding to a different number of spikes.
We refer to these segments as ``spike-count bands.''
The points on the boundary of each spike-count band belong to the stable manifold of the upper saddle equilibrium $\rm SD$; therefore the trajectories flowing forward from any two such boundary points are strongly contracted against one another against the spiking manifold $\rm M_{PO}$.
In fact, the entirety of a spike-count band will flow forward and layer underneath each band associated with a lesser number of spikes, pressing against $\rm M_{PO}$ so quickly that the flows of any two bands cannot be visually distinguished in numerical integrations.
This outside-to-inside scrolling structure of the attractor is common in neuronal models\cite{gilmorepei} and has also been observed in laser models\cite{laserswissroll}, obtaining a moniker of ``g\^ateau roul\'e'' or ``Swiss-roll attractor.''

When transitioning from spiking activity on $\rm M_{PO}$ to subthreshold activity on $\rm M_Q$, trajectories emanating from a given spike-count band are folded and creased about the tangency point $\rm T$ as they are pressed rapidly against the underside of $\rm M_Q$ at the reinsertion loop.
Inside the reinsertion loop, the forward flows of the countably many spike-count bands are tightly layered; the magnified callout diagram at bottom right in Fig.~\ref{fig10} shows this layering.
The folded spike-count bands subsequently flow clockwise along $\rm M_Q$ back to the rainbow line at left, producing several nested Smale horseshoes.

Notably, the attractor is bounded on $\rm M_Q$ by the forward flow of the tangency point $\rm T$ and whichever of the branches $\Gamma_{\rm SD}^\pm$ of the 1-dimensional unstable manifold of the upper saddle equilibrium $\rm SD$ lies furthest from the saddle-focus $\rm SF$.
After an initial transient burst, all trajectories starting on the rainbow line at left in Fig.~\ref{fig10} will begin winding around $\rm M_{PO}$ between $\Gamma_{\rm SD}^\pm$ and the forward flow of $\rm T$ upon subsequent bursts.
This implies that the quantity of spikes per burst that may be observed in the SiN model for a given parameter value---after relaxation---are limited to the range between that of the first bursts observed in the forward flow of $\rm T$ and a neighborhood of $\rm SD$.
Hence the placement of the reinsertion loop on $\rm M_Q$ is the main qualitative feature of the attractor in the SiN model which varies as parameters are manipulated.

\begin{figure}[h!]
  \begin{center}\includegraphics[width=.85\linewidth]{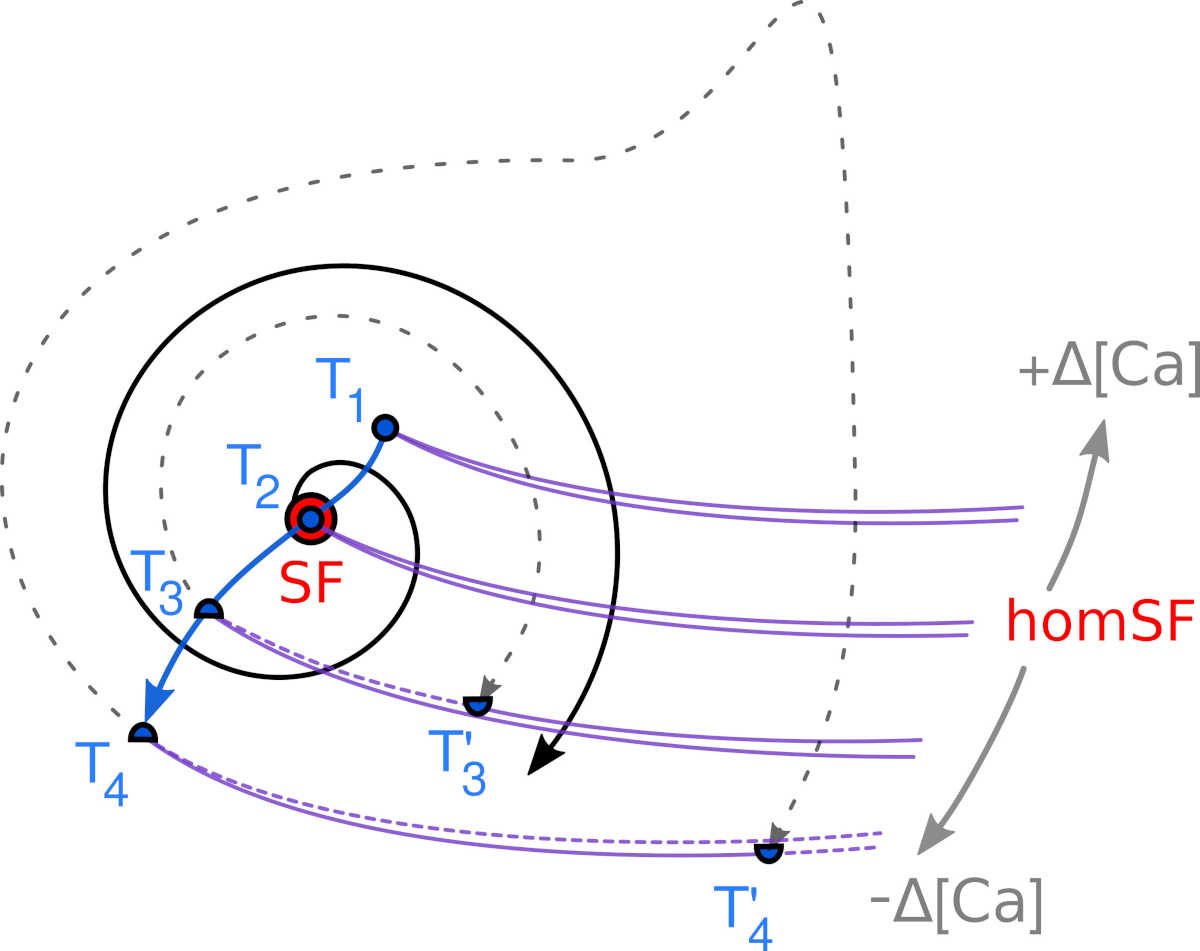}
    \caption{
    	Illustration of the reinsertion loop passing through the saddle-focus (SF) in phase space as the parameters cross the homSF bifurcation curve.
    	The purple curves indicate the location of the reinsertion loop---joining the spiking manifold $\rm M_{PO}$ to the dune $\rm M_Q$---at different values of the parameter $\rm \Delta [Ca]$.
    	The tangency point on the reinsertion loop, labeled T$_i$, is marked by blue dots (or upper half-dots for $\rm T_{3, 4}$).
    	$\rm T_1$ is the location of the tangency point when $\rm \Delta [Ca]$ is selected to the right of the homSF bifurcation curve, while $\rm T_2$ occurs at the homSF bifurcation.
    	The positions of $\rm T_3$ and $\rm T_4$ occur to the left of homSF, where the tangency point flows to the left initially..
    	Trajectories just to the left of $\rm T_{3, 4}$ will go for a long excursion before reinsertion to $\rm M_Q$, whereas this was not the case for $\rm T_{1, 2}$.
    	This has the effect of ``erasing'' a part of the original reinsertion loop, illustrated by dashed segments between $\rm T_{3, 4}$ and $\rm T_{3, 4}'$.
    	For $\rm T_4$, the trajectory spikes before its first return to the reinsertion loop on $\rm M_Q$, completely occluding the orientation-preserving branch from accessibility by any trajectory in the attractor.
	}\label{fig11}
  \end{center}
\end{figure}

Upon the $\rm homSF$ bifurcation curve, there exists a trajectory homoclinic to the saddle-focus equilibrium $\rm SF$.
Thus the tangency point $\rm T$ of the reinsertion loop must lie exactly on the saddle-focus $\rm SF$ at these parameter values.
Fig.~\ref{fig11} shows the movement of the reinsertion loop in a parameter neighborhood of the $\rm homSF$ curve.
As $\rm \Delta [Ca]$ is increased, $\rm T$ moves upward in the $x$ variable, while decreasing $\rm \Delta [Ca]$ moves $\rm T$ downward.
When $\rm \Delta [Ca]$ is decreased from the parameter value lying on $\rm homSF$, $\rm T$ lies below $\rm SF$ and $\rm SF$ is an unstable focus on $\rm M_Q$.
Hence to the left of the $\rm homSF$ curve in the $\rm (\Delta [Ca], \Delta V_x)$-parameter plane, no trajectory starting on the reinsertion loop can return to the segment of the upper half of the reinsertion loop between the tangency point $\rm T$ and the first return $\rm T'$ of the tangency point to the reinsertion loop.
If the trajectory starting at $\rm T$ exhibits one or more spikes before returning to the reinsertion loop, the entire top half of the reinsertion loop becomes inaccessible and the bottom half will become inaccessible on the right side between $\rm T'$ and the upper saddle equilibrium $\rm SD$.
This occlusion is indicated by a dashed segment between $\rm T_{3, 4}$ and $\rm T_{3, 4}'$ on the reinsertion loops in Fig.~\ref{fig11}.
Later, we will show that this occlusion causes one-dimensional return maps for the SiN model to become discontinuous.

\begin{figure}[h!]
  \begin{center}\includegraphics[width=.6 \linewidth]{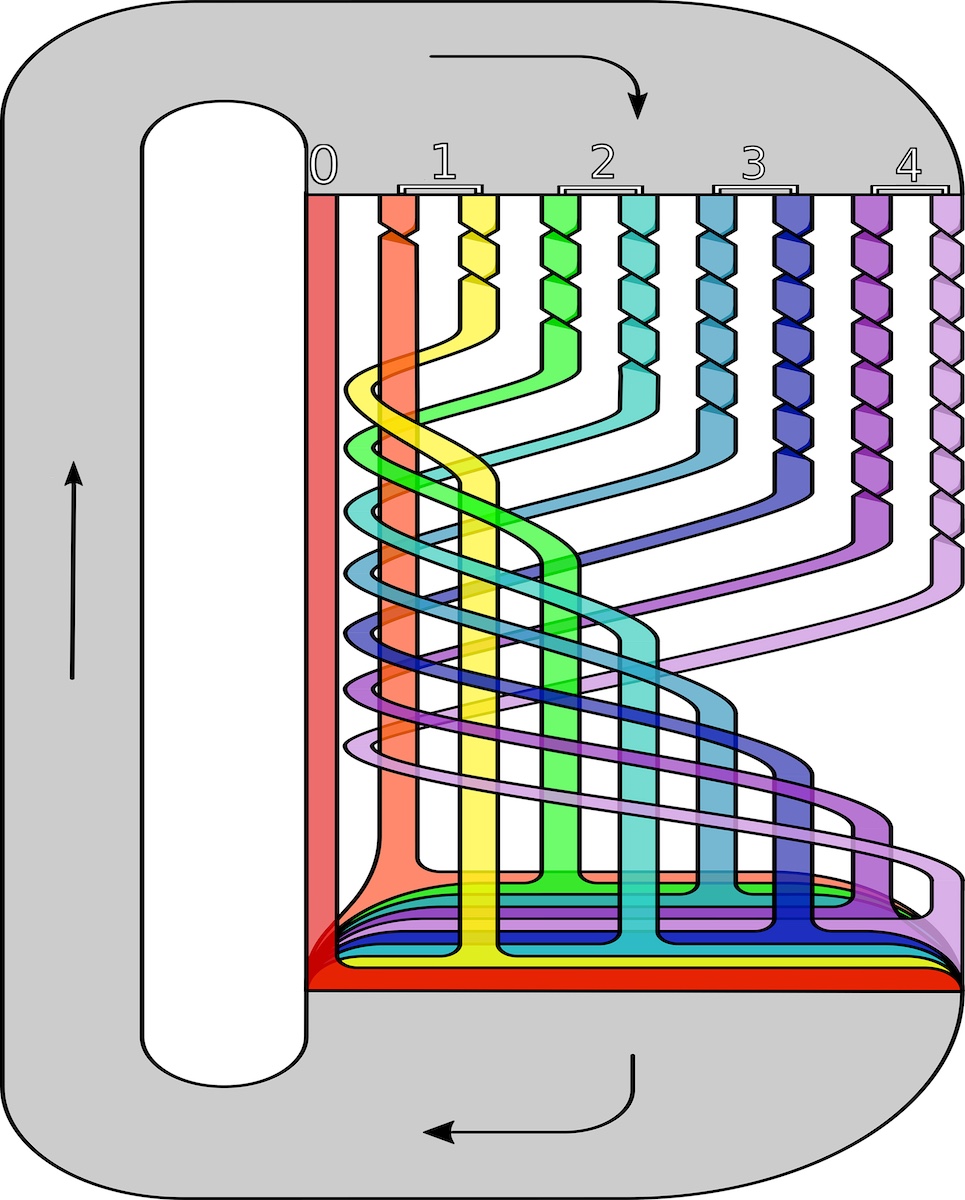}
    \caption{Template describing the dynamics of the SiN model including trajectories with up to four revolutions around the spiking manifold $\rm M_{PO}$. Labels on the upper branch line correspond to the number of times trajectories wind around $\rm M_{PO}$ as they traverse the Swiss roll; i.e., the number of spikes that occur in a burst. Lower spike counts are found on the outer layers of the reinsertion back onto the dune $\rm M_Q$, while higher spike counts are nested deeper within the folded sheets---this is reflected in the order in which branches insert into the lower branch line. This template represents the fully expansive case where the attractor in phase space includes a trajectory homoclinic to the upper saddle equilibrium, spiking four times, as well as a trajectory homoclinic to the saddle-focus equilibrium for each spike count.}\label{fig12}
  \end{center}
\end{figure}

Figure~\ref{fig12} is the template obtained by the Birman-Williams projection of the phase space when there are at most four spikes observed in the attractor.
The template comprises a grey branch representing $\rm M_Q$ and a rainbow-colored array of branches associated to the flow from the rainbow line of Figure~\ref{fig10} around the spiking manifold $\rm M_{PO}$ to the reinsertion loop during a burst.
The branches are conjoined at two branch lines: the upper branch line is associated with the rainbow line at left in Fig.~\ref{fig10} and the lower branch line is associated with the reinsertion loop at right in Fig.~\ref{fig10}.
This simplified representation of the dynamics as a semiflow on a branched-manifold template makes it much easier to understand how trajectories on the attractor intertwine with one another, constrained to flowing clockwise around the template along its many branches.

\subsection{General template}
The template of Fig.~\ref{fig10} can be greatly simplified. In parameter regions where the attractor contains bursts of up to $n$ spikes, the template requires $2n+1$ branches. No template of this kind can serve as a general model for bursting dynamics, as permitting an unbounded quantity of spikes per burst would make the template non-compact. We introduce a recurrent branch, accommodating an arbitrary number of spikes per burst in a single template. The branches of this general template induce a partition for generating symbolic representations of the dynamics across parameter space. The general template is related to the original through a \emph{template inflation} mapping the original template  diffeomorphically onto a subset of the general template and preserving the structure of symbolic itineraries\cite{ghristholmes}; however, we omit the details of this inflation.

\begin{figure}[h!]
  \begin{center}\includegraphics[width=.7\linewidth]{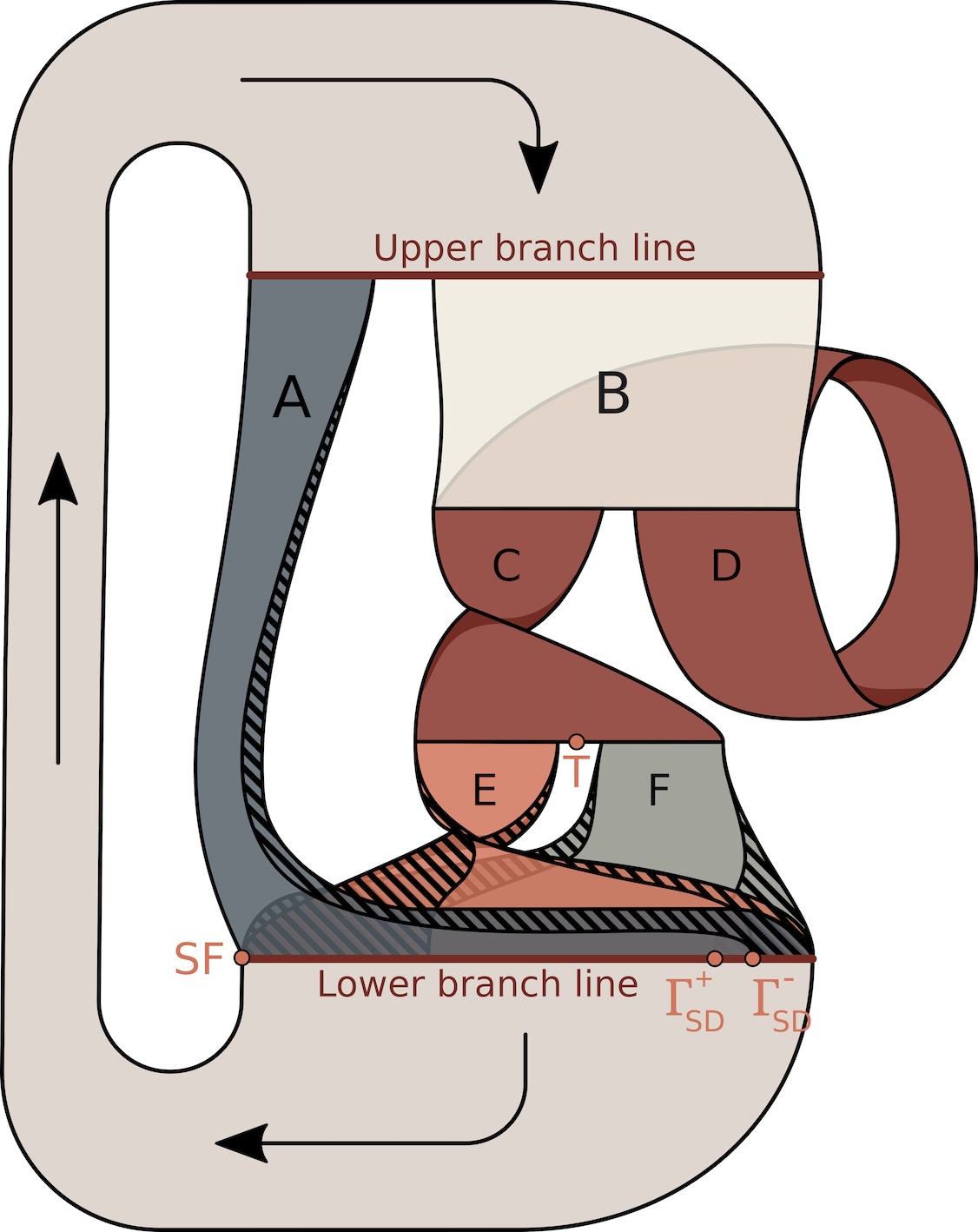}
    \caption{The general template may accommodate any number of spikes. The branch labeled A corresponds to sub-threshold oscillations, while B wraps over $\rm M_{PO}$.  Branch D continues to spike, while C leaves the spiking manifold and reinserts to $\rm M_Q$. Branch E corresponds to the orientation-preserving front segment of the reinsertion loop, while F corresponds to the orientation-reversing back segment. The points SF and SD mark the saddle-focus and upper saddle, respectively. The label T indicates the tangency point. Striped shading represents subsets of the template that may be inaccessible for particular parameter values.}\label{fig13}
  \end{center}
\end{figure}

Figure~\ref{fig13} illustrates the general template, which is reduced to six branches. The split between the A and B branches corresponds to either completing a subthreshold oscillation on the dune $\rm M_Q$ or moving onto the spiking manifold $\rm M_{PO}$.

In the original template, each successive branch has an additional half-twist. This twisting action can be seen in the upper right portion of the template in Fig. \ref{fig12}. In the general template, corresponding twists are accomplished during recurrent passes through the branch labeled by D. Each pass through D represents an additional full twist corresponding to an extra spike. The branch C represents leaving the spiking manifold $\rm M_{PO}$. This branch is not strictly necessary; the general template could be simplified by eliminating the branch C and consolidating the two branch lines at the bottom of branches $\rm B$ and $\rm C$ into a single branch line from which three branches (F, E, and D) emanate. This simplification would result in a non-binary tree structure. Structuring the template as a binary tree simplifies symbolic analysis. Additionally, the half twist of branch C makes the template easier for humans to understand by preserving a clear correspondence between branches and recognizable features of phase space. 

The branches E and F serve to distinguish the orientation-preserving branches associated to trajectories in the full phase space which reinsert to the dune $\rm M_Q$ on the front side of the spiking manifold $\rm M_{PO}$, and the orientation-reversing branches associated to trajectories which reinsert behind $\rm M_{PO}$ respectively. The splitting point on the branch line from which E and F flow corresponds to the tangency point of the reinsertion loop where creasing of a spike band occurs, labeled $\rm T$ in Fig.~\ref{fig10}. The lower branch line stretches between the saddle-focus ($\rm SF$) and whichever of the branches $\Gamma_{\rm SD}^\pm$ of the unstable manifold of the upper saddle ($\rm SD$) lies furthest from the saddle-focus equilibrium $\rm SF$.

The symbolic itinerary of any trajectory is directly associated to the sequence of traversed branches in the associated semiflow on the template of Fig.~\ref{fig13}.
We use double brackets $\llbracket x \rrbracket$ to denote the symbolic itinerary of a trajectory or point $x$ in the phase space of the SiN model, Birman-Williams projected onto the template of Fig.~\ref{fig13}, becoming a semiflow from an initial condition located on the upper branch line.

From a small set of symbolic itineraries associated with specific initial conditions, we can determine how to shade the template for a particular parameter value and directly calculate the topological entropy of the flow in the phase space of the SiN model.
The remainder of this section specifies concrete computational procedures to calculate an efficient proxy for the topological entropy, measuring the complexity of the dynamics for a given parameter value.

\begin{figure}[t!]
	\begin{center}
	\includegraphics[width=\linewidth]{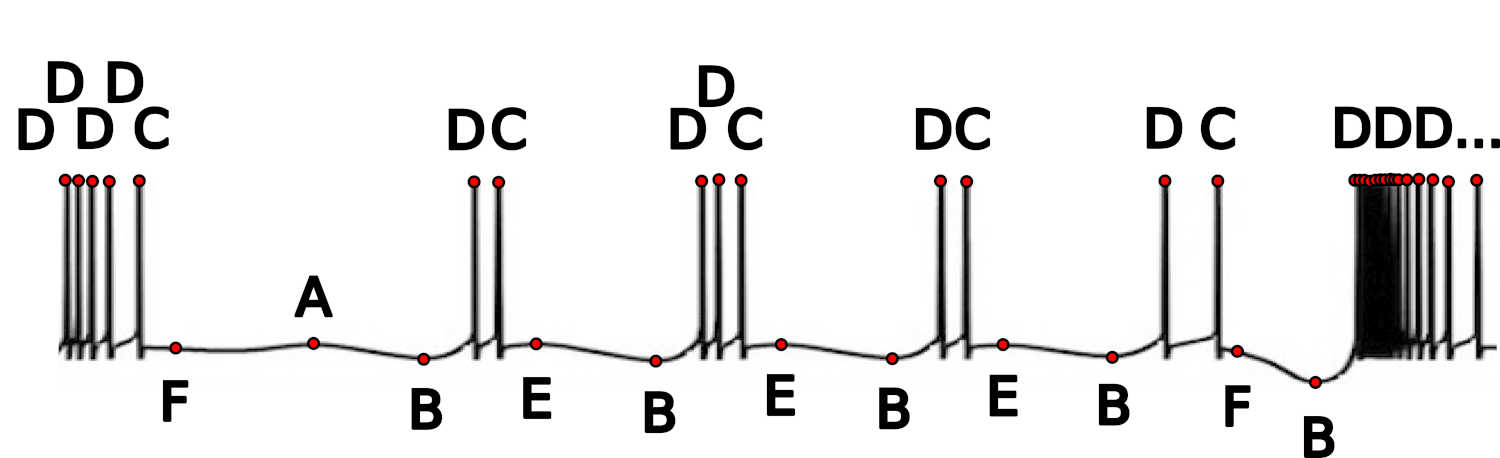}\\
	  \caption{Illustration of how symbols are generated from voltage traces. The symbols correspond to the sequence of traversed branches on the general template.}\label{fig14}
	\end{center}
\end{figure}

We compute the symbolic representation of a trajectory in the SiN model directly from a voltage time series and its time derivative.
Figure~\ref{fig14} provides an example of the qualitative features of a voltage trace corresponding to the sequence of branches traversed by a semiflow on the template of Fig.~\ref{fig13}.
The symbols in the itinerary correspond to the following features apparent in the voltage trace:
\begin{enumerate}
\item A: Subthreshold oscillation.
\item B: Final voltage minimum before a spike.
\item C: Final spike in a burst.
\item D: Non-final spike in a burst.
\item E: Post-burst repolarization, followed by depolarization.
\item F: Post-burst repolarization, followed by hyperpolarization.
\end{enumerate}

In effect, each burst with $n$ spikes produces a symbolic sequence with one $B$-symbol, $n-1$ many D-symbols, one C-symbol, and either an E or F, all in that order. The final symbol is determined by whether the voltage value increases (E) or decreases (F) immediately upon the end of the burst.
Between bursts, the only symbolic sequences which may occur are due to subthreshold oscillations: sequences [AA...AA] of any (potentially zero) length.

\subsection{Inadmissible regions of the template}

Although the compact template shown in Fig.~\ref{fig13} provides a general characterization of the trajectories of the SiN model in a neighborhood of its attractor, for a given parameter value of the system it may be the case that certain trajectories on the general template are inadmissible.

For example, in the region of stable 2-spike bursting, all but finitely many itineraries must exhibit eventually periodic behavior, converging under left-shifts to the periodic itinerary [BDCE BDCE...].

To establish which which portions of the general template are accessible for a given point in the parameter plane, we construct an ordering of trajectories such that trajectories which come first in the ordering always fall to the left on the upper branch line compared to trajectories on the right. For convenience we construct this ordering by assigning each sequence of symbols to a number in the interval $[0, 1]$. We refer to this number as an ``itinerary address.'' 

\begin{figure}[h!]
	\begin{center}
	\includegraphics[width=.9\linewidth]{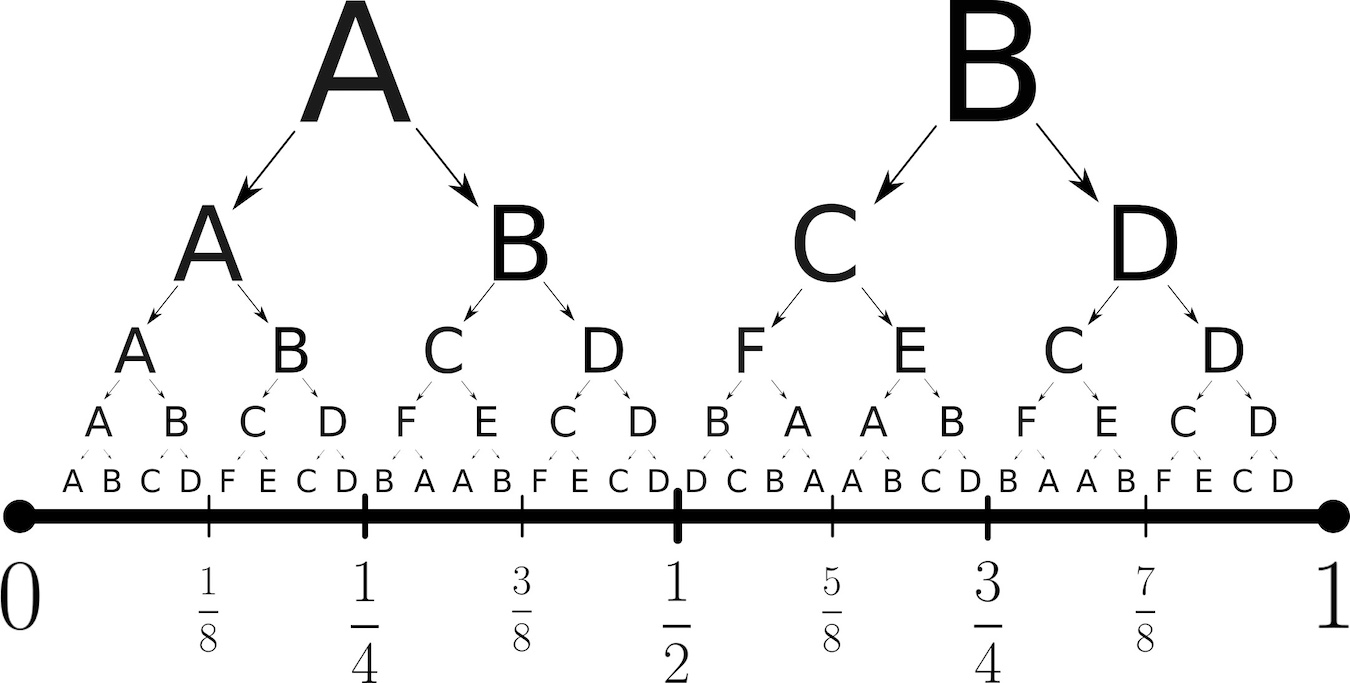}\\
	\caption{Diagram showing how symbolic sequences correspond to addresses in the interval $[0, 1]$.  The interval is successively bisected by the splitting points of the template. At each binary choice as a trajectory flows through the template, the left branch is assigned to the lesser subinterval, and the right is assigned to the greater subinterval.}\label{fig15}
	\end{center}
\end{figure}

We assign the value $0$ to the saddle-focus line and $1$ to the upper saddle point. To determine the coordinate of a trajectory, we successively bisect the interval $\left(0, 1\right)$, choosing the left or right subinterval at each step based on the branch the trajectory takes in the binary tree of possible paths. A truncated binary tree of possible branches and its relationship to the itinerary addresses between $\left(0, 1\right)$ is shown in Fig.~\ref{fig15}.

The assignment of symbols to subintervals depends on the accumulated orientation reversals in a subsequence. The accumulated orientation of a subsequence changes sign every time an orientation reversing branch (C or E) is encountered. Therefore when we bisect an interval, we assign the left subinterval first in alphabetical order if the total number of orientation flips in the subsequence prior to that symbol is even. We assign in reverse alphabetical order when the total number of orientation flips is odd. Trajectories which reach splitting points are identified with the midpoints. This addressing scheme agrees with the \emph{signed lexicographical ordering} ($\succ$) on itineraries.

For instance, consider a trajectory that passes over branch A, revolves clockwise around the template upon the unlabeled branch, traverses B, spikes twice by traversing D and C once each, and finally returns to the dune, $\rm M_Q$, through E. To construct the corresponding addressing subinterval of $[0, 1]$, we follow the itinerary taken by the trajectory:

\begin{itemize}
\item A: Select the left half, $\left(0, \frac{1}{2}\right)$.
\item B: Select the right half, $\left(\frac{1}{4}, \frac{1}{2}\right)$.
\item D: Select the right half, $\left(\frac{3}{8}, \frac{1}{2}\right)$.
\item C: Select the left half, $\left(\frac{3}{8}, \frac{5}{16}\right)$. Orientation flipped so that subsequent symbols are treated in reverse lexicographical ordering.
\item E: Select the right half, $\left(\frac{9}{32}, \frac{5}{16}\right)$. Orientation flipped so that subsequent symbols are treated in lexicographical ordering again.
\end{itemize}

Similarly, ignoring the first partial burst in Fig.~\ref{fig14} corresponding to the symbolic subitinerary [DDDDCF], one can determine from the itinerary corresponding to the depicted voltage trace that this trajectory would have crossed the upper branch line of the template at an address between $\frac{1826240}{2^{22}} \approx 0.4354095$ and $\frac{1826241}{2^{22}} \approx 0.4354098$.

For the purpose of determining the shading of the template, it suffices to calculate the addresses of the forward-iterated boundaries of certain branches. By convention, we refer to these boundaries as the ``right'' and ``left'' boundaries with respect to the semiflow through the template. For a branch with zero twists, this means the ``left'' boundary is on the right side of the illustration in Fig.~\ref{fig13}. For example, the left boundary of the A branch appears on the right-hand side of the branch as depicted, while the left boundary of the E branch appears at the top-right and bottom-left sides of that branch in the figure.

The particular branch boundaries which are shaded are chosen to constrain trajectories on the template in accordance with the topology of the attractor in phase space. There are generally two different ways in which the behavior of trajectories on the attractor is constrained: 1) The maximum number of spikes per burst for any trajectory is finite, bounded by the number of spikes observed in the first burst of the maximal itinerary $\max\,\{\llbracket \Gamma_{\rm SD}^- \rrbracket, \llbracket \Gamma_{\rm SD}^+ \rrbracket\}$, associated to a branch of the unstable manifold of the upper saddle, under the signed lexicographical ordering. 2) Trajectories in the attractor cannot reach the saddle-focus in parameter regions that are far from the saddle-focus homoclinic bifurcation curve, $\rm homSF$; this is governed by the placement of the tangency point $\rm T$ of the reinsertion loop on the dune $\rm M_Q$ in relation to the saddle-focus equilibrium $\rm SF$.

The maximum number of spikes per burst is limited by the two outgoing separatrices ($\Gamma_{\rm SD}^-$ and $\Gamma_{\rm SD}^+$) of the upper saddle $\rm SD$.
Every branch boundary which terminates at the right side (away from the SF) of the lower branch line is affected: these are the left side of A, the right side of E, and the left side of F.
The shading of left boundary of A and the right boundary of E are associated with the itinerary of $\Gamma_{\rm SD}^-$, $\llbracket \Gamma_{\rm SD}^- \rrbracket$, while the shading on the left boundary of F is associated with $\llbracket \Gamma_{\rm SD}^+ \rrbracket$.
There is no way to know \emph{a priori} which of $\llbracket \Gamma_{\rm SD}^- \rrbracket$ or $\llbracket \Gamma_{\rm SD}^+ \rrbracket$ is lesser; these may appear in either order, or in the same place, on the lower branch line.

To represent in the context of semiflows on the template the inaccessibility of parts of the dune $\rm M_Q$ near the saddle-focus $\rm SF$ by the attractor, we shade branches E and F at those boundaries which terminate at the lower branch line at the point labeled $\rm SF$. On the right side of the saddle-focus homoclinic curve $\rm homSF$ in the parameter plane, the entire upper and lower halves of the reinsertion loop are fully accessible to trajectories in the attractor, as shown by the purple curve touching the point labeled as $\rm T_1$ in Fig.~\ref{fig11}. In this case, the shadings for the left boundary of E and the right boundary of F agree; they correspond to the same itinerary address, calculated from the forward iterates of the tangency point $\rm T_1$. Upon passing over $\rm homSF$ from right to left, the two boundaries separate. Just to the left of the $\rm homSF$ curve, the shading at the right boundary of branch F still corresponds to $\rm \llbracket T_3 \rrbracket$, but the shading at the left boundary of branch E corresponds to $\rm \llbracket T_3' \rrbracket$, which is the same as the itinerary $\rm \llbracket T_3 \rrbracket$ without the first burst---the initial subitinerary of $\rm \llbracket T_3 \rrbracket$ is removed up to and including the first occurrence of branch A, E, or F to obtain the itinerary $\rm \llbracket T_3' \rrbracket$. In practice, this alleviates the requirement of computing an additional trajectory to determine the shading of the template.

For parameter values far to the left of $\rm homSF$, branch E becomes entirely inadmissible, so we shade the entirety of E so that the left half of branch C is inadmissible (one can also explicitly shade the left side of branch C so that only neighborhoods of boundaries of branches are shaded). In this case, the entire upper half of the reinsertion loop (in this case touching $\rm T_4$ in Fig.~\ref{fig11}) becomes inaccessible. The bottom half of the reinsertion loop is occluded to the right of $\rm T_4'$, so that no spiking trajectory returns to the dune $\rm M_Q$ further away from the saddle-focus $\rm SF$ than this point. Hence the left side of branch F must be shaded at least up to the address associated to the itinerary $\rm \llbracket T_4' \rrbracket$. However, recall that the shading at the left boundary of F is also associated with $\llbracket \Gamma_{\rm SD}^+ \rrbracket$. This conflict is reconciled by shading the left boundary of F at the address of the lesser of $\llbracket \Gamma_{\rm SD}^+ \rrbracket$ and $\rm \llbracket T_4' \rrbracket$ under the signed lexicographical ordering.

\subsection{Quantifying complexity}\label{lzcomplexity}

So far, we have described the paths under the semiflow on the template of Fig.~\ref{fig13} associated to trajectories in the phase space of the SiN model by encoding these trajectories as symbolic itineraries on branch labels from the template.
An equivalent, yet more concise symbolic encoding of a bursting trajectory will be referred to as a \emph{signed spike-count sequence} (SSCS). This encoding applies for any itinerary beginning with A or B, so that initial conditions on the upper branch line are assigned a unique SSCS according to their forward semiflows.

The number $0$ is assigned to a subthreshold oscillation which flows through branch A. For a burst, the total number of spikes is recorded, corresponding to the number of occurrences of symbols C and D. If the burst ends with an E-symbol, the number of spikes observed during the burst is committed to the SSCS; this value is positive, associated with the fact that the branches traversed preserve orientation in the signed lexicographical ordering. On the other hand, if the burst ends with an F-symbol, the negation of the number of spikes occurring in the burst is committed to the SSCS; the fact that this value is negative indicates that the branches traversed are in total orientation-reversing.
Then, again omitting the initial partial bursting subitinerary in Fig.~\ref{fig14} corresponding to the symbolic subsequence [DDDDCF] as well as [BDDD] for the incomplete final burst, the symbolic sequence [ABDCEBDDCEBDCEBDCF] corresponds to the SSCS [0, 2, 3, 2, -2]; the remainder of the itinerary for this trajectory, if integration and encoding were continued indefinitely, would produce an infinite SSCS.
Of course, the SSCS encoding scheme cannot represent trajectories which converge to a tonic-spiking regime with itinerary [DDD...].

\begin{figure*}[t!]
	\begin{center}
	\includegraphics[width=\linewidth]{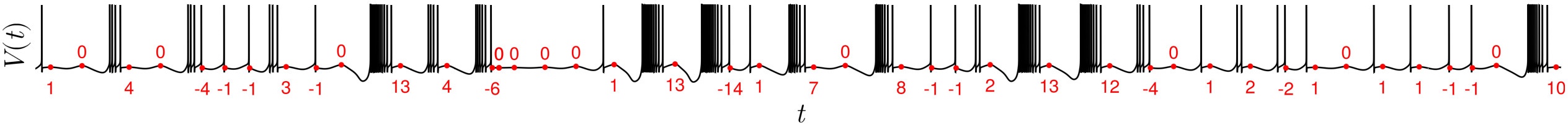}\\
	  \caption{Signed spike-count sequence (SSCS) in red corresponding to a long, chaotic voltage time series.
	  Depicted in black is the $V$ time series of the $\Gamma_{\rm SD}^-$ trajectory at value $(-38.285, -0.9)$ in the $\rm (\Delta [Ca], \Delta V_x)$-parameter plane.}\label{fig16}
	\end{center}
\end{figure*}

Figure~\ref{fig16} depicts a long, chaotic voltage time series having SSCS [1, 0, 4, 1, 0, 1, -1, 1, 0, 2, -2, -1, -14, ...].
The events detected in the time series $V(t)$ corresponding to the recognition of a new signed spike count in the SSCS are marked in red.

Appendix IV describes the encoding algorithm for converting voltage time series data from numerically integrated trajectories to SSCSs.

Because we are able to determine the symbolic itineraries associated with each of the critical trajectories $\Gamma_{\rm SD}^\pm$ and those flowing forward from T and T$'$---each corresponding to the itinerary of a critical point of the one-dimensional return map described in the next section---it is possible to calculate the topological entropy of the SiN model for a particular parameter value from the SSCSs associated with each of these trajectories; this approach to quantifying the dynamical complexity of the model comes from the Milnor-Thurston kneading theory\cite{milnorthurston}.
However, because the upper saddle equilibrium is highly unstable and the branches $\Gamma_{\rm SD}^\pm$ of its unstable manifold each tend to traverse the entire attractor (including near the tangency point T), it is possible to simplify the calculation of the topological entropy by eliminating its dependence on the symbolic itinerary associated to initial conditions $\rm T$ and $\rm T'$.

While it is possible to perform this simplified calculation explicitly, to do so would extend beyond the scope of this paper; we will instead use an auxiliary method to illustrate the dynamical complexity throughout the $\rm (\Delta [Ca], \Delta V_x)$-parameter plane.
This auxiliary method will be to consider the entropy rates of symbolic itineraries $\rm \llbracket T \rrbracket, \llbracket T' \rrbracket, \llbracket \Gamma_{\rm SD}^- \rrbracket, \llbracket \Gamma_{\rm SD}^+ \rrbracket$ of critical trajectories in the phase space associated to the shading in the template of Fig.~\ref{fig13}; there is not a direct correspondence between these individual entropy rates and topological entropy of the SiN model, but they should be highly correlated.
Because the trajectory $\Gamma_{\rm SD}^-$ tends to lie close to $\Gamma_{\rm SD}^+$ by virtue of both the fast-slow timescale separation and the observation that $\Gamma_{\rm SD}^-$ tends to eventually traverse most of the attractor in the phase space, we consider only the symbolic itinerary of the trajectory $\Gamma_{\rm SD}^-$.
As symbolic itineraries are in one-to-one correspondence with SSCSs, we suggest that the entropy rate of the SSCS associated with $\Gamma_{\rm SD}^-$ is an effective measure of the complexity of the dynamics.

Past literature has indicated the suitability of the Lempel-Ziv 1976 (LZ76) complexity\cite{lz76} in approximating the entropy rate of neuronal voltage traces\cite{lzentropyratespikingneuron}, although such studies considered the entropy rate of sequences encoding spiking behavior in fixed time bins in lieu of a generating partition for neuronal bursting dynamics.
To estimate the entropy rate of the SSCS, we calculate the LZ76 complexity of the SSCS of $\Gamma_{\rm SD}^-$ in parameter scans.

Figure~\ref{fig5}B shows in red the LZ76 complexity (above a certain threshold value) of the SSCS associated with $\Gamma_{\rm SD}^-$ in a scan over the $\rm (\Delta [Ca], \Delta V_x)$-parameter plane.
These LZ76 complexities appear to highlight the parameter regions which have positive leading Lyapunov exponent, as seen in Fig.~\ref{fig5}A; thus the LZ76 complexities of these SSCSs effectively locate chaotic parameter values, with some false positives between successive period-doubling curves near spike-adding transitions.
Many of the stability windows observed in the Lyapunov scan can be clearly observed in the LZ76 scan as well.
This close agreement between the Lyapunov and LZ76 scans serves as evidence that the SSCS symbolic encoding of trajectories is an effective representation of dynamical behavior, further supporting the use of templates in finding a partition for the phase space of the SiN model.

\section{1D return maps of an interval}

\begin{figure*}[t!]
  \begin{center}
  \includegraphics[width=.9\linewidth]{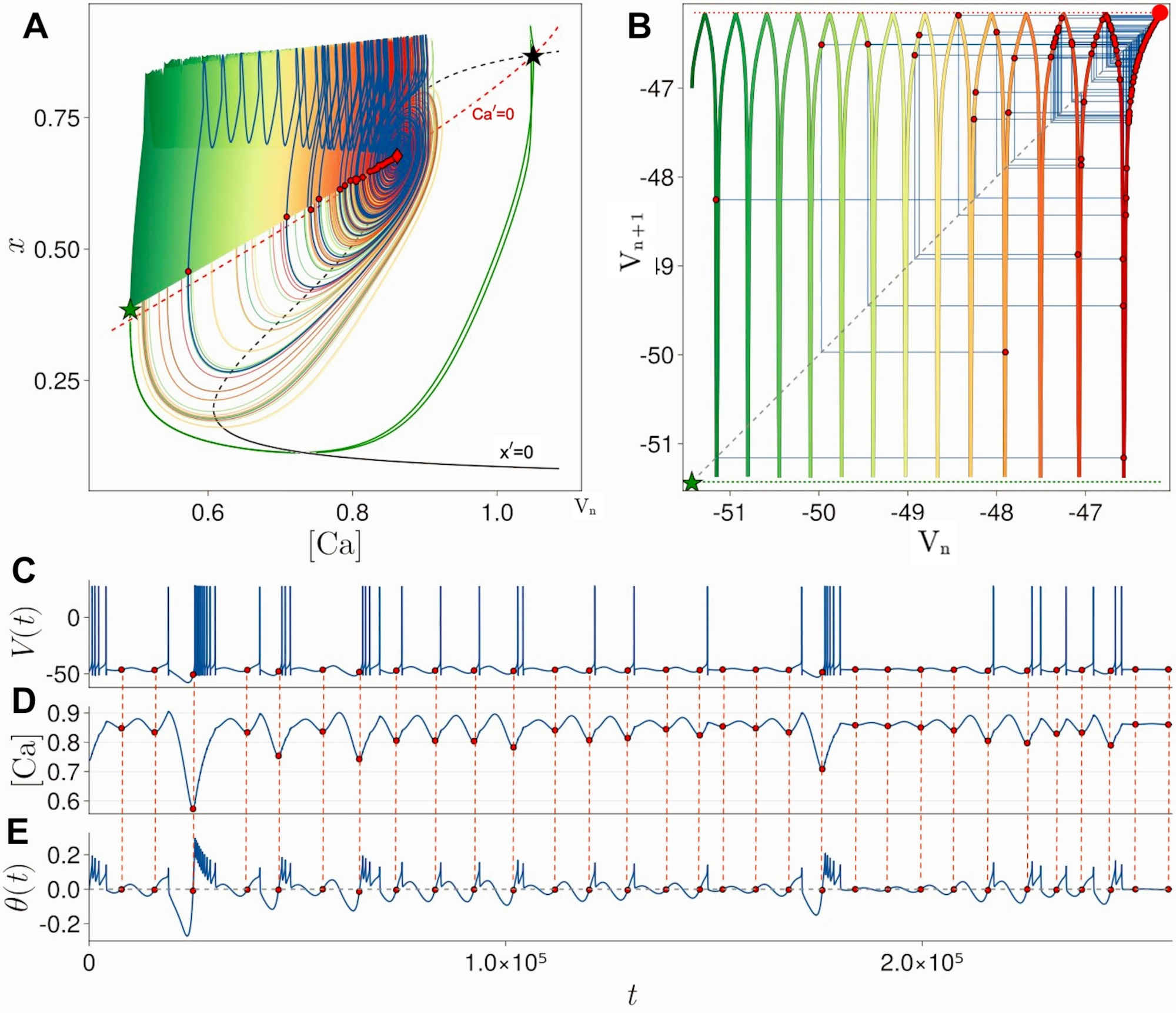}
  \caption{Panel A shows the construction of a typical 1D return map of an interval of phase-space trajectories that originate and terminate on a cross-section drawn between the origin and the position of the saddle-focus (red diamond at the crossing of the $x$ and the $\rm [Ca]$ nullclines) in the slow $\rm ([Ca],\,x)$ projection. The range of the map is determined by the one-dimensional outgoing separatrices from the upper saddle (marked with a black star). Panel~B depicts the one-dimensional return map $f: V_n \to V_{n+1}$ of an interval with 13 self-similar structures for $\Delta \rm V_x = -1.08$ and $\Delta \rm{[Ca]} = -36.3$, populated with chaotic transient iterates (red dots). The crossings of the dashed $45^\circ$ line with the map graph indicate (mostly unstable) fixed points. The coloring of the 1D map corresponds to the coloring of the phase-space trajectories in Panel~A. Three panels below present the corresponding time series of a chaotic trajectory: voltage trace in Panel~C, $\rm [Ca]$ trace in Panel~D, and the angle $\theta$ trace, used for calculating return events on the cross-section, in Panel~E. Zeroes of the $\theta$-variable closely approximate the position of the $\rm [Ca]$ nullcline, as marked by the red dots where the trajectory intersects the cross-section from below in the projection shown in Panel~A.}
  \label{fig17}
  \end{center}
\end{figure*}

The onset of chaos to the right of the central homoclinic curve $\rm homSF$ (see Figs.~\ref{fig1} and \ref{fig5}) can be analyzed using one-dimensional return maps that map a cross-section of the slow-motion dune back onto itself. The construction of these maps is detailed in Fig.~\ref{fig17}A. The chosen cross-section is placed along a straight line connecting the origin to the saddle-focus (lower saddle), which sufficiently approximates the $\rm [Ca]$ nullcline. Note that using calcium minima, determined by the conditions $\rm [Ca]' = 0$ and $\rm [Ca]'' > 0$, will yield similar return maps. Regardless of how the cross-section is selected, its essential property is that it should be nowhere tangent to the trajectories on the dune $\rm M_Q$.

We note that in Ref.~\cite{butera98}, which analyzed a multi-bursting neuronal model, a similar cross-section for generating 1D point-wise return maps was chosen along a line corresponding to a saddle-node bifurcation of the fast subsystem. This bifurcation line demarcates the boundary separating tonic-spiking from quiescent states in the 2D slow phase subspace.

Maps for this study are taken from a straight line connecting the origin, $\rm [Ca] = 0,\, x = 0$, to the saddle-focus. The construction process of the computed 1D return maps includes the following steps:
\begin{enumerate}
  \item Create a sample by linearly interpolating 500 $(\rm [Ca],\, x)$-pairs between the origin and the saddle-focus, $({\rm [Ca]_{SF}},\, x_{\rm SF})$.
  \item Estimate the fast subsystem variable initial conditions by solving for fast subsystem equilibrium. This places the initial conditions near $\rm M_{PO}$.
  \item Integrate trajectories from these initial points for at least 50ms and until they intersect they cross-section after one global excursion. 
  \item This crossing is computed with low-to-high zero-crossings of the the test function $\theta = \tan^{-1} \left( \frac{x(t)}{[\mathrm{Ca}](t)} \right) - \tan^{-1} \left( \frac{x_{SF}}{\rm [Ca]_{SF}} \right)$.
  \item Solve again for the fast subsystem equilibria at the final state of the $\rm [Ca]$- and $x$-variables. This resolves error that arises due to the fact that the initial conditions are not exactly on $\rm M_Q$.
\end{enumerate}

The return map plots are generated by plotting the initial condition voltage values against the return voltage values, yielding the map $f: V_n \mapsto V_{n+1}$. This choice of coordinate is chosen only because voltage is commonly observable in neural systems. In general, any coordinate, such as $\rm [Ca]$, will produce a topologically map, so long as it regularly parameterizes the return section as a one-dimensional submanifold of the dune.

The return maps have a characteristic structure consisting of a succession of ``arches'' (see Fig.~\ref{fig17}B), each corresponding to a burst with a distinct number of spikes. Depending on the parameter values, the chaotic bursting traces and the corresponding maps will have a smaller number of accessible arches. See the bifurcation diagram in Fig.~\ref{fig5}C, which is partitioned into zones indicated by the largest possible numbers of spikes per burst within.
  
The maximum of each arch occurs when trajectories are located closest to the saddle-focus along the return section. The minima between arches correspond to unidirectional heteroclinic connections from the saddle-focus to the upper saddle. The unstable manifold of the upper saddle forms the boundary of the dune in phase space and thus bounds from below the range of the return map as well. On the saddle-focus homoclinic curve, homSF, the saddle-focus lies on the reinsertion loop where the spiking manifold $\rm M_{PO}$ curls beneath and merges back onto the dune $\rm M_Q$. To the left of the homoclinic curve, the saddle-focus moves across the tangency point of the reinsertion loop further away from the spiking manifold, downward on the dune.

\begin{figure*}[t!]
  \begin{center}
  \includegraphics[width=.999\linewidth]{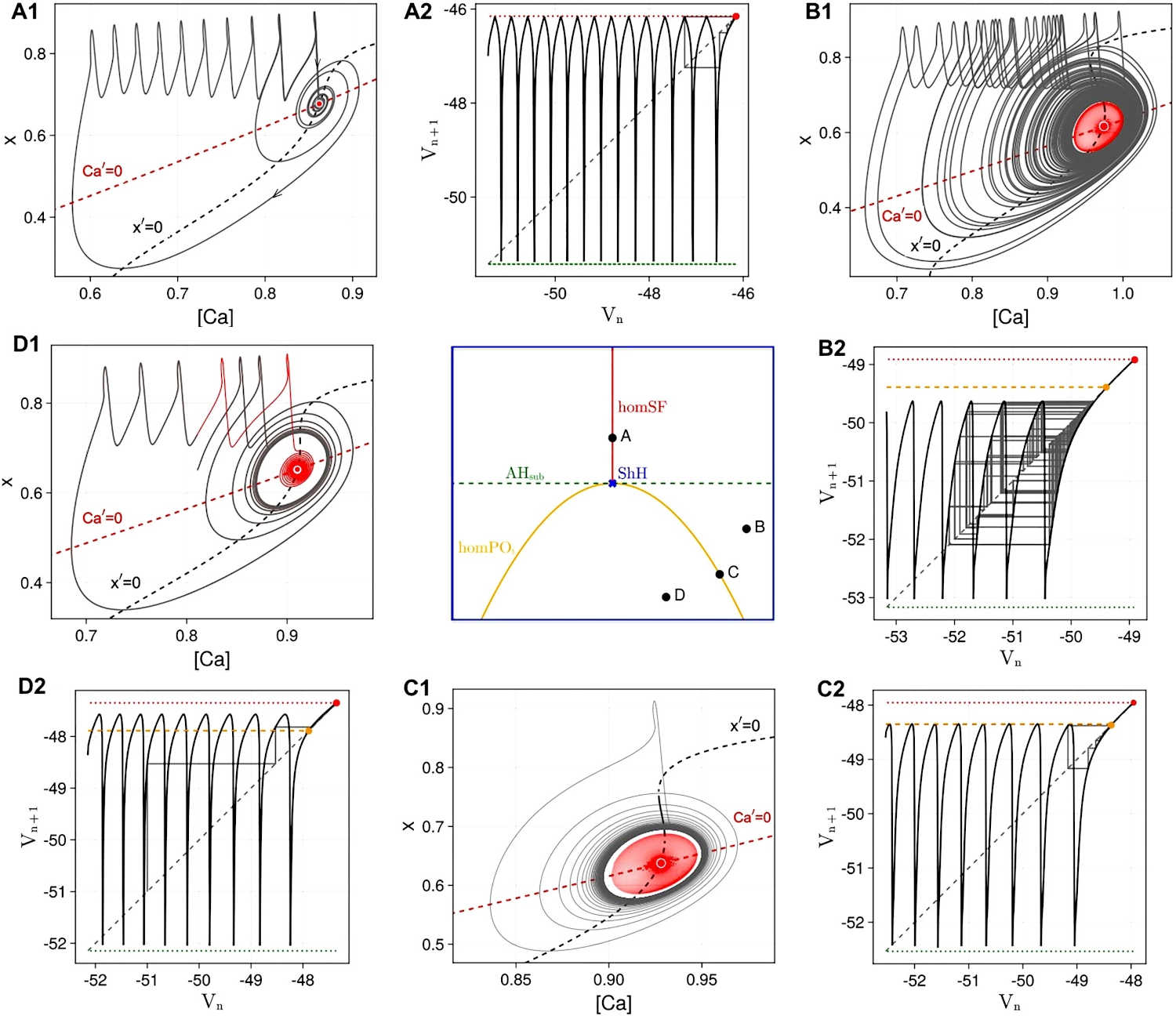}\\
    \caption{Unfolding of the cod-2 Shilnikov-Hopf bifurcation illustrated with selected phase space trajectories and corresponding one-dimensional (1D) return maps. The central panel is an unfolding sketch including three principal bifurcation curves --- AH, homSF, and $\rm homPO_t$ (for reference, see the sweep in Fig.~\ref{fig1}). The outer panels are organized in pairs labeled A1–D2; each pair corresponds to a single parameter value and is represented by the corresponding return map along with the $([{\rm Ca}],\, x)$-phase portrait. In the maps, the dotted red line marks the position of the fixed point (FP) (red dot), representing either the Shilnikov saddle-focus when the latter is unstable or a stable focus when the red FP is stable. The dotted orange line marks the location of the unstable FP (orange dot), representing the saddle periodic orbit (PO) emerging through the subcritical AH bifurcation. Sampled trajectories are shown in grey in the maps and phase portraits. For the phase space trajectories pictured in the corner panels, red portions of the trajectory indicate convergence to the stable focus. The $\rm [Ca]$ and $x$ nullclines in the phase planes are drawn as dotted red and black curves, respectively. Panels A1 and A2 illustrate the homoclinic orbit of the Shilnikov saddle-focus occurring on the bifurcation homSF curve. Panels B1 and B2 demonstrate the coexistence of a chaotic attractor and a stable equilibrium state, separated by the unstable PO (orange FP) that arises below the $\rm AH_{sub}$ curve in the parameter space. Panels C1 and C2 illustrate the onset of homoclinic tangencies on the $\rm homPO_t$ curve. Panels D1 and D2 show how chaos begins to vanish after the homoclinics to the saddle PO become transverse (above the orange FP), and transient trajectories find their way towards the stable focus (red dot).}\label{fig18}
  \end{center}
\end{figure*}

To illustrate how to interpret these maps, we consider Fig.~\ref{fig17}B as an example. This map represents the SiN model in the chaotic region, specifically on the homoclinic homSF curve. Starting from the right, the first fixed point (FP) on the map occurs where it becomes tangent to the $45$-degree line (bisectrix); by construction, this fixed point corresponds to the saddle-focus, depicted as a thick red diamond on the dune in Fig.~\ref{fig17}A. The red dotted line at the top of the frame represents the pre-images of the saddle-focus; thus, any point on the map that intersects this dotted red line is associated with a homoclinic orbit to the corresponding fixed point.The first monotonically increasing segment of the map, before the first local minimum on the left, corresponds to the oscillatory flow on the slow dune near the saddle-focus, terminating at the sub-threshold (i.e., zero-spike) heteroclinic connection to the upper saddle. The first arch to the left of the subthreshold monotonic section, characterized by a single local maximum, corresponds to regular or chaotic voltage traces with a single spike. The next arch to the left corresponds to two-spike bursting trajectories, and so forth, as reflected in the voltage trace shown in Fig.~\ref{fig17}C, where the number of spikes per burst varies unpredictably. At each local minimum (cusp-shaped) between arches, the map exhibits a discontinuity due to slight differences in how the upper and lower outgoing separatrices of the unstable manifold of the upper saddle return to the section. These two unstable separatrices in phase space are plotted as green curves in the $\left(\rm [Ca],\, x \right)$ phase-space projection in Fig.~\ref{fig17}B. The point where the two separatrices hit the return section is illustrated at the bottom of Fig.~\ref{fig17}B by a dotted green line marked with a green star. Note that at these parameter values within the chaotic region, the separation between the upper and lower branches is too small to notice, but this is not always the case. 
Each arch of the one-dimensional interval map is composed of two orientation-preserving branches: one monotonically increasing and one decreasing. Each branch corresponds to a separate strip of the topological template shown above and is twisted once for each spike (see Fig.~\ref{fig12}).

The sequence of repeated arches continues to the left of the plot frame in Fig.~\ref{fig17}B. However, these additional branches or arches may not be accessible because they extend beyond the reach of the unstable separatrices of the upper saddle. The point where the dotted green line meets the identity line marks the most extreme accessible value. Additional details on how to interpret key features of the map relevant to specific bifurcations are included in their respective sections below.

\subsection{Homoclinic tangencies of periodic orbits}

The lower boundary of the chaotic region in the $\rm \left( \Delta [Ca], \Delta V_x \right)$ parameter plane (Figs.~\ref{fig1} and \ref{fig5}) is delineated by the $\rm homPO_t$ curve. This curve signifies the onset of homoclinic tangencies, where the tubular stable manifold of a saddle periodic orbit --- emerging from the subcritical Andronov-Hopf bifurcation ($\rm AH_{sub}$) --- first contacts its unstable manifold (see Fig.~\ref{fig4}). Below the $\rm AH_{sub}$ bifurcation curve, a stable equilibrium state (focus) exists, which can coexist with chaotic dynamics above the $\rm homPO_t$ curve.

In this region, the stable manifold of the saddle periodic orbit separates the basin of attraction of the chaotic attractor from that of the stable focus within the five-dimensional phase space. The unstable manifold of this saddle periodic orbit extends into the chaotic attractor over time. As the saddle periodic orbit increases in size farther from the $\rm AH_{sub}$ curve, its stable manifold eventually intersects the basin of the chaotic attractor along the $\rm homPO_t$ curve. For parameter values below this curve, chaotic behavior becomes transient, and the system ultimately settles into the basin of attraction of the stable equilibrium state.

Figure~\ref{fig18} presents several panels that showcase matching return maps and corresponding phase portraits of the slow dynamics near the cod-2 Shilnikov-Hopf (ShH) point in the parameter space. In the one-dimensional maps, the height of each arch indicates the closest approach to either the saddle-focus or the stable fixed point (FP) --- depicted by the red dot --- that is accessible from the chaotic attractor.

Figures~\ref{fig18}A1 and A2 illustrate the formation of two homoclinic orbits associated with the Shilnikov saddle-focus at parameter values $\Delta {\rm [Ca]} = -36.6$ and $\Delta \rm V_x = -1.078$. In these panels, the unstable fixed point (FP), corresponding to the saddle-focus, is located in the upper corner of the one-dimensional map. The homoclinic orbits can be identified by the alignment of the peaks of the arches with the dotted red line, which represents the preimages of the saddle-focus.

As we move downwards along the red homSF curve in parameter space, the derivative of the map at the saddle-focus FP increases to +1. At this critical point, the map becomes tangent to the 45-degree identity line, marking the Shilnikov-Hopf (ShH) bifurcation. Below this tangency point, at the subcritical Andronov-Hopf bifurcation ($\rm AH_{sub}$), a new FP (orange dot) emerges in the map. This FP represents a saddle periodic orbit formed through the Andronov-Hopf bifurcation.

Directly below the ShH bifurcation point, the Andronov-Hopf bifurcation immediately suppresses chaotic behavior, leading to a stable, hyperpolarized state that dominates this region. On either side of the $\rm hom_{SF}$ curve, after transitioning through the $\rm AH_{sub}$ bifurcation, the saddle periodic orbit has room to grow before interacting with the chaotic attractor. In the maps, the available space for the saddle periodic orbit to expand without interfering with the chaotic attractor is reflected by the vertical separation between the peaks of the arches and the dotted red line, which indicates the height corresponding to the saddle-focus.

The pre-images of the emerging FP in the maps are highlighted by a dotted orange horizontal line, representing the stable manifold of the saddle periodic orbit. Intersections between the map and the dotted orange line indicate homoclinic orbits to the saddle periodic orbit (orange FP). If the critical points of the arches do not reach the height of the saddle-focus (unstable red FP), homoclinic connections to the saddle periodic orbit are absent, and the chaotic attractor persists. However, when the critical points of the arches surpass this threshold, the dynamics of the system change, allowing trajectories to escape the chaotic attractor and converge to the stable focus (stable FP), leading to the suppression of chaos.

Figures~\ref{fig18}B1 and B2 computed at $\rm \Delta [Ca] = -2.75$ and $\rm \Delta V_x = -2.04$ depict the region situated below the $\rm AH_{sub}$ curve and above the $\rm homPO_t$ curve, where bistability between chaotic dynamics and a quiescent stable focus occurs. In this region, the height of the fixed point (indicated by the orange dashed line) is above the peaks of the arches in the map, resulting in no interaction between the chaotic attractor and the stable focus.

Along the $\rm homPO_t$ curve, the peaks of the arches in the map become tangent to the dashed orange line, as illustrated in Figs.~\ref{fig18}C1 and C2. This represents the critical case where the stable and unstable manifolds of the saddle periodic orbit (PO) form a non-transverse intersection. Intersections between the map and the dashed orange line indicate homoclinic orbits to the saddle periodic orbit (orange fixed point). The example in Figs.~\ref{fig18}C1 and C2 is calculated at $\Delta {\rm [Ca]} = -1.69$ and $\Delta \rm V_x = -16.56$.

In Figs.~\ref{fig18}D1 and D2, computed at $\Delta {\rm [Ca]} = -1.54$ and $\Delta \rm V_x = 23.12$ respectively, each arch of the map forms two intersections with the dashed blue line, each corresponding to a transverse homoclinic orbit to the saddle periodic orbit. The portions of the map lying above the dashed line converge rapidly to the saddle-focus.

We emphasize the importance of these interval maps, such as in Fig.~\ref{fig18}C2, which were used to compute the edge-of-chaos boundary ($\rm homPO_t$) in the bifurcation diagrams shown in Figs.~\ref{fig1} and \ref{fig5}. Accurately determining this boundary would be impossible relying solely on trajectory computations in the five-dimensional phase space of the SiN model.
  
\begin{figure}[t!]
  \begin{center}
  \includegraphics[width=1.0\linewidth]{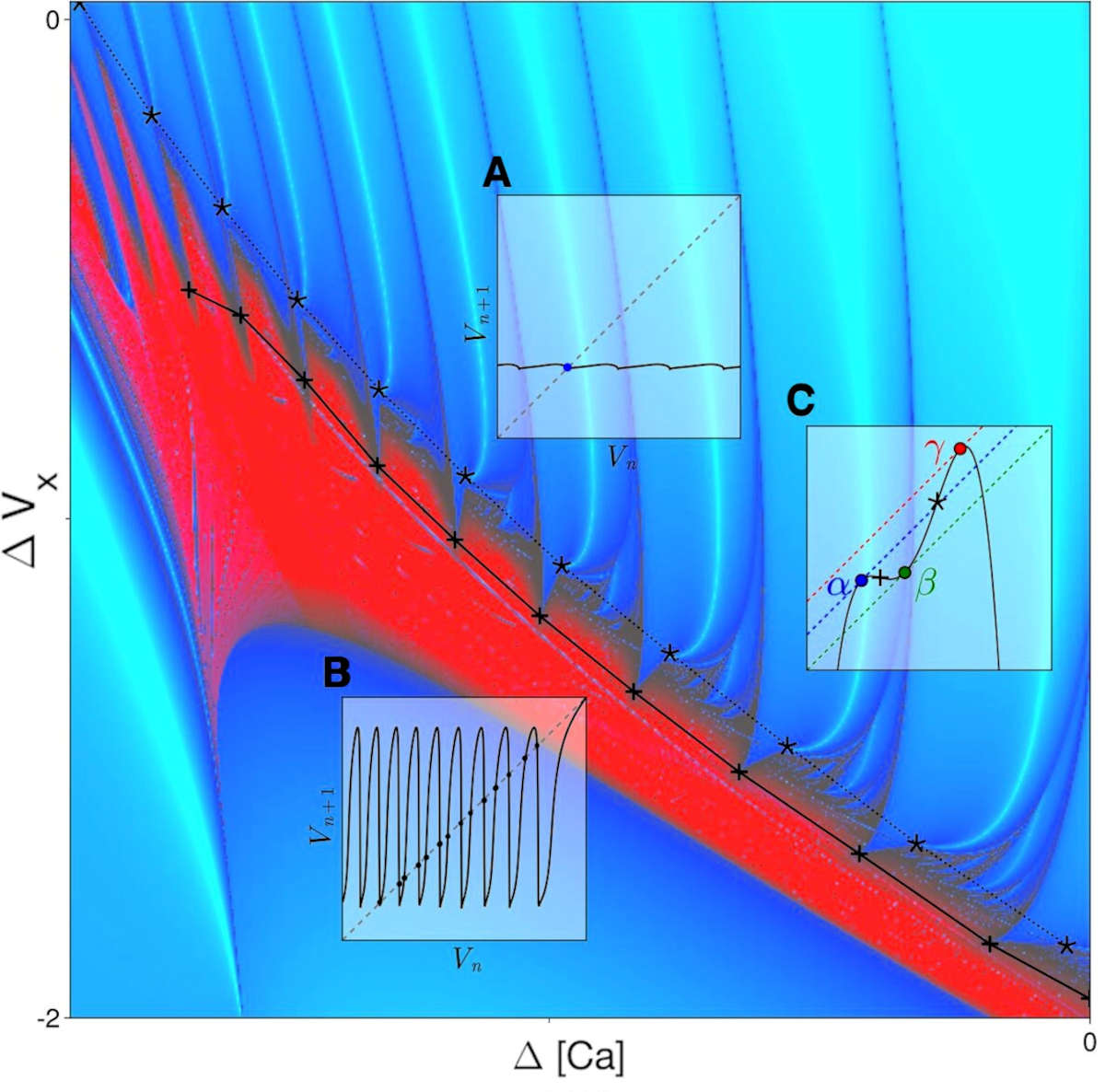}
  \caption{
    Fragment of the bifurcation diagram illustrating the transition between regular bursting (top blue region) and chaotic bursting dynamics (red chaos-land). A transition occurs when the arches of the 1D map increase in height, as illustrated by the contrast between inset panel~A (stable bursting) and panel~B (chaos). A transitory zone exists between these two cases, as represented in panel~C. Along the upper and lower boundaries of the transition region, cusp points of the map occur, marked by ``*'' and ``+'' respectively. These cusps occur when an inflection point on the map forms a cubic tangency with the $45^\circ$  bisectrix line.
    Inside the transition region, three points with derivative equal to +1 exist labeled by $\alpha, \beta,$ and $\gamma$ in inset panel~C. When these points touch the bisectrix, saddle-nodes occur, corresponding to saddle-node bifurcations of periodic orbits in the SiN model. The coloring of the Lyapunov spectrum is the same as in Fig.~\ref{fig1}.
    }\label{fig19}
  \end{center}
\end{figure}

\begin{figure*}[t!]
  \begin{center}
  \includegraphics[width=1.0\linewidth]{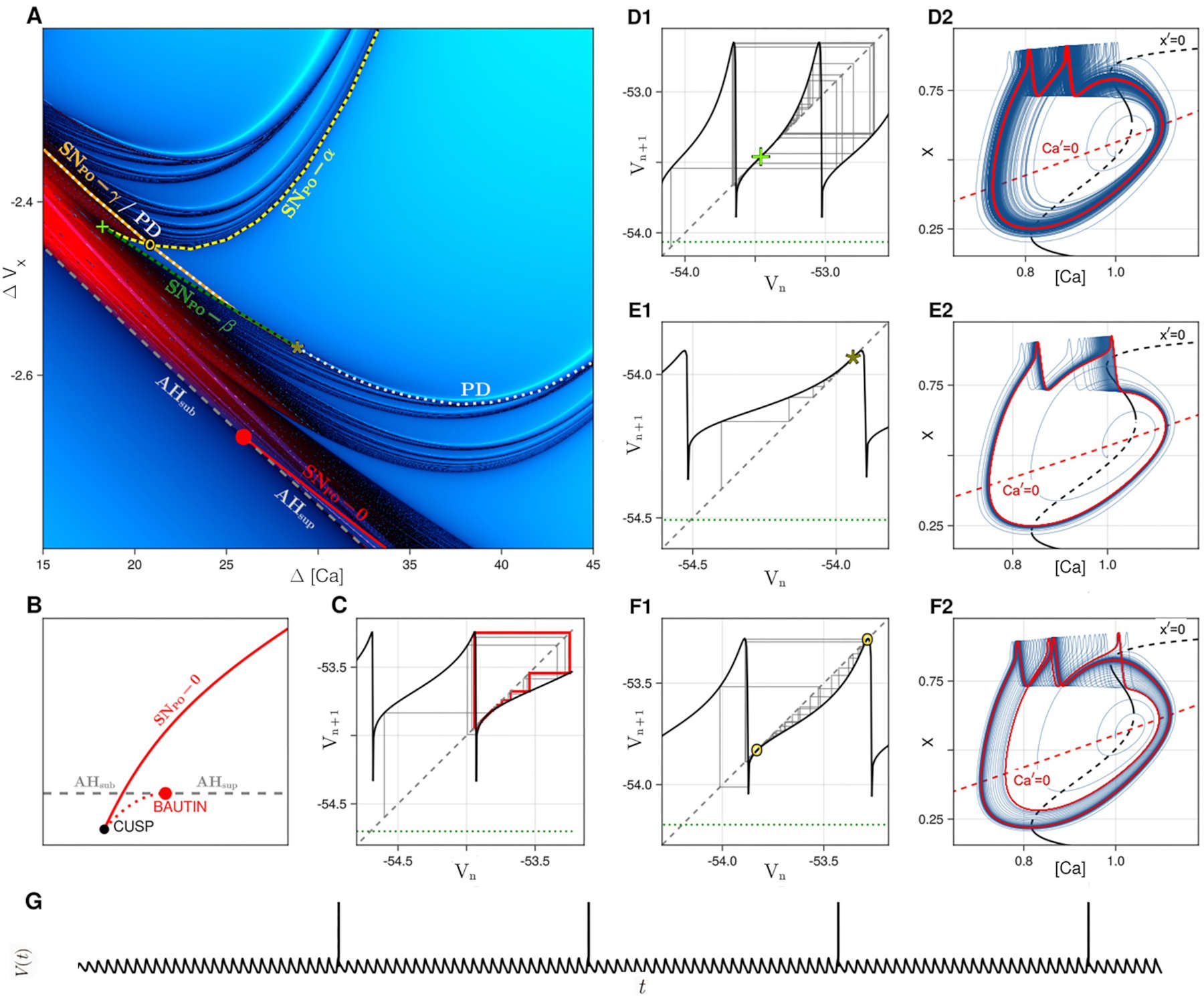}
    \caption{
    Exploration of the spike-adding route to chaos using 1D return maps and slow phase portraits. (A) A magnified two-parameter sweep displays the first and second largest Lyapunov exponents, overlaid with bifurcation curves derived from manual inspection of the maps. These curves highlight the loss of stability of the stable two-spike burster. The curves SN$\rm _{PO} - \alpha$ (yellow), SN$\rm _{PO} - \beta$ (green), and SN$\rm _{PO} - \gamma$ (orange) are shown. SN$\rm _{PO} - \gamma$ is closely shadowed by a period-doubling (PD) curve, which extends beyond the right edge of the parameter frame, marked by a dotted white line. SN$\rm _{PO} - \alpha$ and SN$\rm _{PO} - \beta$ coincide at the cusp point marked ``+'', located at ($\Delta \rm Ca,\, \rm \Delta V_x$) = (18.2, -2.42). This cusp is illustrated in panels D1 and D2: the map is shown in D1, and the slow phase plane in D2, with the cusp orbit highlighted by a red trajectory and transients in blue. SN$\rm _{PO} - \beta$ and SN$\rm _{PO} - \gamma$ coincide at the cusp point marked ``*'', located at ($\Delta \rm Ca,\,\Delta V_x$) = (28.9, -2.57). Panels E1 and E2 illustrate this cusp formation. The two SN$\rm _{PO}$ orbits meet at the point marked ``o'', with ($\Delta \rm Ca,\,\Delta V_x$) coordinates = (20.83, -2.45), where a blue-sky heteroclinic cycle exists. This case is illustrated in panels F1 and F2, with the blue-sky connections in blue and the two saddle-node orbits in red. (B) A sketch of a zoomed-in bifurcation diagram around the Bautin point at ($\Delta \rm Ca,\,\Delta V_x$) = (25.9, -2.67) (shown as a red dot), which gives rise to a saddle-node of periodic orbits curve that turns around at a nearby cusp. The red curve SN$\rm _{PO} - 0$ begins at this cusp and defines the loss of stability of subthreshold oscillations. Panel~C demonstrates the stability of the blue-sky orbit at SN$\rm _{PO} - \alpha$ via a long trajectory on the map. Panel~G shows a voltage trace near a blue-sky orbit exhibiting long episodes of subthreshold oscillations interrupted by single spikes.
    }\label{fig20}
  \end{center}
\end{figure*}
\subsection{Route from bursting to chaos through widening spike-adding transitions}

In the following, we analyze various global transitions that occur in the SiN model as the parameter value exits the wide chaotic region over its upper-right boundary toward regular bursting activity (see Fig.~\ref{fig19}). Along this boundary, multiple thin regions of spike-adding chaos expand and merge into one another. The nature of these transitions can be clearly understood through the use of one-dimensional (1D) interval maps. These maps reveal a transitional zone where the spike-adding transitions blend together, organized by a sequence of cod-2 cusp points of periodic orbits.

To understand the transition from bursting to chaos, we first explain the sequence of bifurcations that occur along each spike-adding transition, moving from right to left in the bifurcation diagram in Fig.~\ref{fig19}. The shape of the maps in the bursting region can be seen in Fig.~\ref{fig19}A. At high values of the $\Delta x$-parameter, the map is compressed vertically: the arches are flat, but the valleys between the arches remain steep, although they are shortened. As the $\Delta {\rm [Ca]}$-parameter decreases, the maps compress horizontally, which has the effect of ``pulling'' each arch rightwards across the 45-degree line (dotted grey line).

The sequence of bifurcations that occurs in spike adding from $n$ to $n+1$ spikes can be inferred from the sequence of geometric features of each arch that crosses the identity line as $\Delta {\rm [Ca]}$ decreases:
\begin{enumerate}
\item Homoclinic to the upper saddle: The minimum on the right side of the arch crosses the identity line, signifying a homoclinic orbit to the upper saddle. Only after this homoclinic bifurcation does the new arch become accessible, in the sense that the unstable manifold of the saddle-focus contains trajectories that spike an additional time at precisely this point.
\item Saddle-node of periodic orbits (SN$_{\rm PO}$): As the left side of the $n$-spike arch of the map becomes tangent to the $45$-degree line, the saddle periodic orbit---which previously emerged from the homoclinic orbit to the upper saddle---merges and annihilates with the stable bursting orbit.
This saddle-node of periodic orbits leaves no remaining stable periodic trajectories, leading to a brief window of chaos. The chaotic trajectories in this window will undergo a mix of both $n$-spike and $n+1$-spike excursions. Note that the derivative of the map at this tangency must equal one.
\item Period-doubling cascade: As the derivative of the map at its intersection with the 45-degree bisectrix decreases past $-1$, a period-doubling cascade occurs in reverse.
\item Stabilization of the bursting orbit: When the derivative equals $-1$, the new $n+1$-spike bursting orbit becomes stable and dominates the long-term dynamics.
\end{enumerate}

As the arches of the map grow taller, the expanding proportion of the map increases, making chaos more common toward the bottom of each spike-adding transition. This stretching of the map is illustrated by the contrast between Figs.~\ref{fig19}A (bursting) and \ref{fig19}B (chaos). When the arches of the map are either very short or very tall, there is precisely one point on each arch where $\frac{df}{dV_n} = 1$, and one point where $\frac{df}{dV_n} = -1$. These two points separate the contracting and expanding regions in the 1D interval map. As the map is perturbed by horizontal compression, the relative sizes of the expanding and contracting regions reveal why stretching the map spreads chaos in the parameter plane. When the map is short, the expanding regions are very small, leading to predominantly bursting behavior. When the map is tall, the expanding regions are very large, leading to predominantly chaotic dynamics.

A transitional zone exists between the short bursting and tall chaotic cases, illustrated by the region between the dotted and solid lines in the bifurcation diagram in Fig.~\ref{fig19}. In this zone, there are not one, but three points on each arch where $\frac{df}{dV_n} = 1$, and each arch gains an additional region of expansion. The transitional zone is interspersed with a large number of cod-1 $\rm SN_{PO}$ bifurcations that can be separated into three classes based on the order in which they occur in the domain of the map. These classes correspond to $\frac{df}{dV_n} = 1$ points on the map labeled as $\alpha$, $\beta$, and $\gamma$ in order of increasing $V_n$, as seen in Fig.~\ref{fig19}C. The corresponding $\rm SN_{PO}$s are labeled as $\rm SN_{PO}-\alpha$, $\rm SN_{PO}-\beta$, and $\rm SN_{PO}-\gamma$, respectively.

The boundaries of the transitional zone are defined when two $\frac{df}{dV_n} = 1$ points merge to form an inflection point with derivative one. The $\alpha$ and $\beta$ points come together to form the lower boundary toward chaos (solid line), and the $\beta$ and $\gamma$ points come together to form the upper boundary toward bursting (dotted line). Whenever the bisectrix line passes through these points, a cusp bifurcation occurs involving two periodic orbits. The cusps formed from $\alpha$ and $\beta$ are labeled with an asterisk (``*''), and the cusps formed from $\beta$ and $\gamma$ are labeled with a plus sign (``+'') in both Fig.~\ref{fig19} and Fig.~\ref{fig20}. The map and corresponding phase space trajectories for this cusp are shown in Figs.~\ref{fig20}D and E. The saddle cusp orbit is red, superimposed over the chaotic attractor in blue.

In the maps, these cusps appear as cubic tangencies that unfold into two saddle-node bifurcations. In the ODE system, each of these cod-2 cusp bifurcations of periodic orbits unfolds into two saddle-node bifurcations of periodic orbits ($\rm SN_{PO}$). Figure~\ref{fig20}A shows a magnified portion of the bifurcation diagram from Fig.~\ref{fig19}, revealing the first few spike-adding transitions, including the $\rm SN_{PO}-\alpha$, $\rm SN_{PO}-\beta$, and $\rm SN_{PO}-\gamma$ curves, which appear as dotted yellow, green, and yellow lines, respectively. Wedges of stability begin at the $\alpha , \beta$ cusps, marked with ``+'' symbols. The map and corresponding phase space trajectories for this cusp are shown in two panels of Fig.~\ref{fig20}D. The stable cusp orbit is red, with a transient shown in blue.

All three $\rm SN_{PO}$ bifurcations appear to be blue-sky catastrophes. In each case, a homoclinic saddle-node bifurcation of periodic orbits leads to the emergence of a stable periodic orbit with infinite period and length. This occurs because the saddle-node periodic orbits possess a homoclinic appendage that becomes stable as the point of intersection is approached (see the background section for more theoretical discussion).

We examined various ``ghost'' periodic orbits of the saddle node near the bifurcation point and found that, in all cases we investigated, these orbits eventually become stable. Computing the stability at the saddle-node bifurcation point can be challenging because it depends on the global features of the map. To assess stability, we tested these orbits by iterating the one-dimensional interval map using linear interpolation between sample points. For example, Fig.~\ref{fig20}C illustrates how iterates of the 1D interval map converge to a stable blue-sky trajectory shown in red. One potential drawback of this method is that false positives for stability are possible due to the finite precision in numerical simulations.

Additionally, there are multiple bistable regions within each wedge. These bistable bursting patterns exist near $\rm SN_{PO}-\gamma$ when the bisectrix (45-degree line) intersects the return map in two places: just to the right of the $\gamma$ point and in the contracting region between $\alpha$ and $\beta$. Since the peaks near $\gamma$ are rather steep, the $\gamma$ point is close to the adjacent point where $\frac{df}{dV_n} = -1$. This proximity means that when the map is perturbed by horizontal compression, the period-doubling cascade closely follows $\rm SN_{PO}-\gamma$. This region of bistability begins near the ``*'' cusp where the rightmost fixed point (FP) of the map emerges, as shown in Fig.~\ref{fig20}E1. The bistable behavior ends when the leftmost FP of the map disappears through $\rm SN_{PO}-\alpha$. This map is visually indistinguishable from the one shown in Fig.~\ref{fig20}F1, with the saddle-node orbit marked by the leftmost ``o'' point (at $\alpha$). The steepness of the peaks is evident at the rightmost ``o''-marked point (at $\gamma$). 

The map shown in Fig.~\ref{fig20}E1 corresponds to the point in the bifurcation diagram (in Fig.~\ref{fig20}A) where $\rm SN_{PO}-\beta$ and $\rm SN_{PO}-\gamma$ coalesce (marked with ``o''). At this point, two saddle-node bifurcations of periodic orbits occur simultaneously. The homoclinic structure here is unique because a ``blue-sky heteroclinic'' cycle occurs, where the unstable and stable manifolds of two saddle-node periodic orbits meet. This heteroclinic cycle appears to coexist with homoclinic trajectories to each of these cycles. Figures~\ref{fig20}E and \ref{fig20}F display the map and the corresponding slow phase projections of the phase space of the SiN model.  The ghosts of two saddle-node orbits (in red) are highlighted in Fig.~\ref{fig20}F2, indicating that bistability is nearby. This heteroclinic point forms a corner of chaos in the parameter space, which separates the strongly chaotic (red) region in Fig.~\ref{fig20}A and the weakly chaotic spike-adding region.

At the end of the spike-adding cascade --   where the 1-spike bursting behavior transitions 
 into 0-spike or subthreshold oscillations -- the boundary of chaotic dynamics is defined by an $\rm SN_{PO}$ curve. This curve originates near the cod-2 Bautin point (the red dot labeled BP in Figs.~\ref{fig20}A and B), which separates the subcritical and supercritical branches of the Andronov-Hopf (AH) bifurcation curve in the parameter plane. We refer to this curve as $\rm SN_{PO}-0$ or $\rm homSN_{PO}$ to emphasize the homoclinic appendage. Upon meticulous examination of the maps, we observe that the Bautin point actually unfolds downward into the tonic spiking region --- adjacent to the bottom-left corner in Fig.~\ref{fig20}A --- where it intersects the $\rm SN_{PO}-0$ curve at a cusp. An exaggerated sketch of this small region of the bifurcation diagram is provided in Fig.~\ref{fig20}B.

\subsection{Transformation of the interval maps through homSF}

On the left side of the saddle-focus homoclinic bifurcation curve, homSF, the saddle-focus lies on the spiking manifold $\rm M_{PO}$. Identifying a suitable return section that yields faithful one-dimensional maps --- maps where the trajectories are directly interpretable as corresponding trajectories in the system of differential equations --- is more challenging in this case. These difficulties make it impractical to numerically compute one-dimensional maps for this region of the parameter plane. However, by utilizing qualitative insights from the topological model, we are able to create sketches of the maps and the corresponding slow phase space trajectories.
The problem that arises on the left side of the parameter plane is that obvious choices for a return section --- such minima of the $\rm [Ca]$-variable --- inevitably become tangent to the flow. This tangency occurs as the return section crosses the junction where $\rm M_{PO}$ squeezes back onto the manifold of stable equilibria in the fast subsystem,  $\rm M_Q$. As a result, the map becomes discontinuous at these points of tangency, making it problematic to interpret the trajectories.
Choosing a return section through phase space that entirely avoids tangency between the section and the discontinuity is not possible. However, there is a natural choice of section that allows the trajectories of the map to be reasonably interpreted.

This section can be constructed by ensuring that the return section intersects the tangency point of the reinsertion loop, such that the discontinuity occurs at local maxima of each arch in the 1D map. This tangency point is referred to as $\rm T$ in section \ref{section:topologicalorganization}, as illustrated in Fig.~\ref{fig10}.

\begin{figure}[t!]
  \begin{center} 
  \includegraphics[width=.9\linewidth]{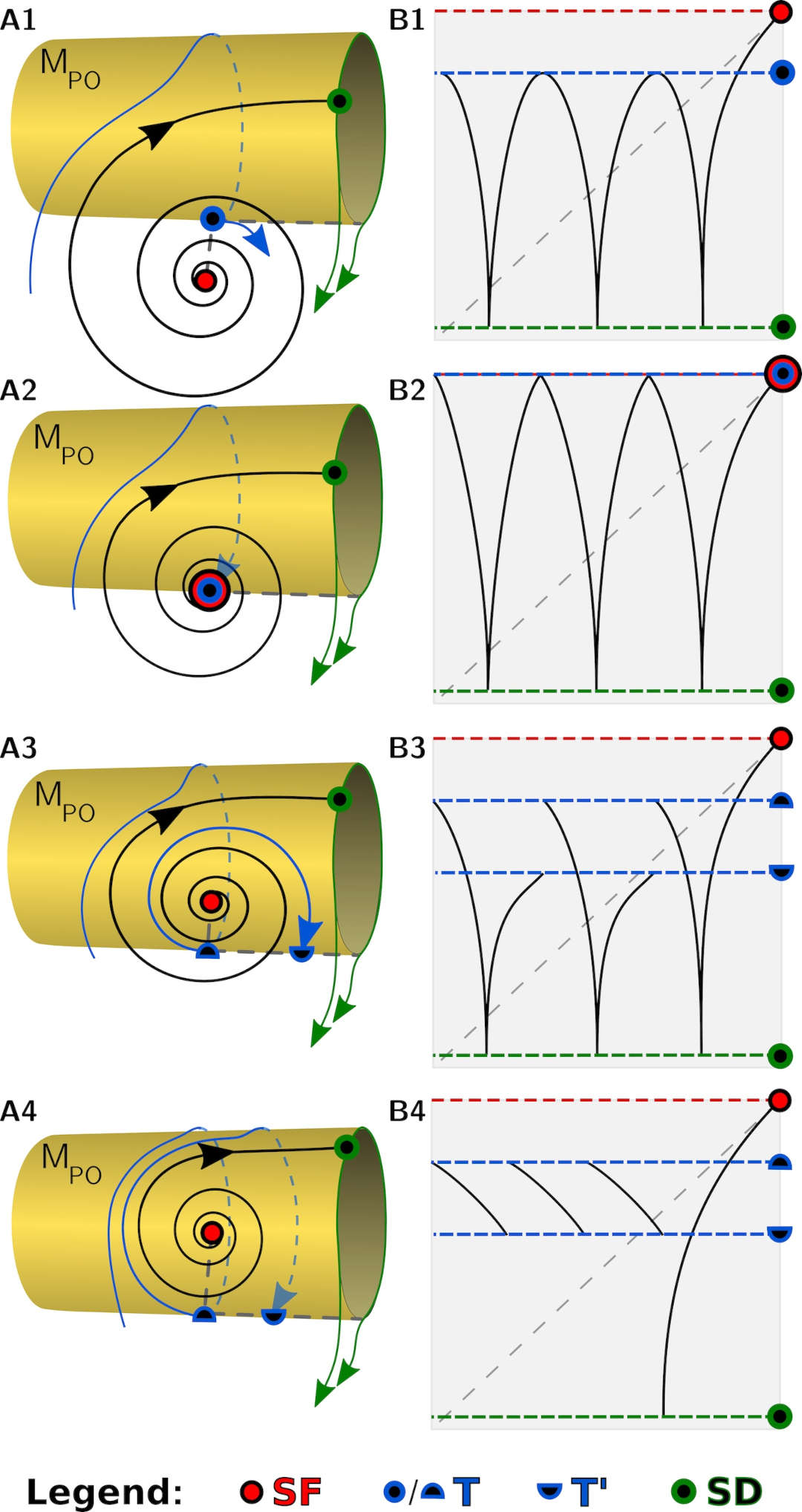} \\
  \caption{Sketch of the transition in return maps across homSF. The 2D phase space is illustrated on the left, with four snapshots of the return maps shown on the right. The scaffolding of the dynamics is represented by a cylinder corresponding to $\rm M_{\rm PO}$, and the persistent connection between the saddle-focus (red dot) and the upper saddle (green dot) is shown as a black line. The unstable manifold of the upper saddle is depicted with green lines. A trajectory that hits the tangency point T is shown in blue. A neighborhood of T projected onto the map is marked by a blue dot. To the right of the homSF bifurcation curve, this neighborhood remains connected, as shown in panels~A1 and A2. To the left of homSF, the neighborhood splits into two, as seen in panels~A3 and A4. In this case, the upper semicircle shows T, and the lower semicircle shows the first forward iterate of T. The bottom of $\rm M_{PO}$ represents the reinsertion loop where $\rm M_{\rm PO}$ meets $\rm M_Q$ (white region below). The return section is represented by a dashed gray line. In the interval return maps, the locations of the saddle-focus SF, tangency point T, and upper saddle SD are marked, along with their pre-images indicated by dashed lines. The bisectrix is shown as a dashed gray line.}\label{fig21} 
  \end{center}
\end{figure}

Figure~\ref{fig21} shows four qualitative sketches at snapshots of the transition from the right side of homSF to the left side of the bifurcation diagram. The movement of the tangency point $\rm T$, and the associated occlusion of the reinsertion loop (corresponding to truncated arches in the map), is shown in Fig.~\ref{fig11}, and there is a close connection between these figures. The solid part of the purple reinsertion loop in Fig.~\ref{fig11} forms the arches of the return map, while the dotted section of the loop corresponds to the discontinuities.
Figs.~\ref{fig21}A1 and A2 illustrate the configuration to the right of homSF where the saddle-focus is on $\rm M_Q$; this situation is detailed throughout the one-dimensional maps section of the paper. Here, trajectories that flow in close proximity on either side of the tangency point (near the blue trajectory) will intersect the return section in near proximity.
Panels A2 and B2 of Fig.~\ref{fig21} depict the homoclinic-to-saddle-focus case (on $\rm homSF$), where the tangency point $\rm T$ coincides with the saddle-focus $\rm SF$ as it reaches the reinsertion loop where $\rm M_Q$ and $\rm M_{PO}$ meet.
Panels C1 and C2 of Fig.~\ref{fig21} show the configuration to the left of the curve homSF. The blue trajectory lands at the point T on the return section. Trajectories that begin just to the left of the blue trajectory will narrowly miss the section near $\rm T$ (upper semicircle) before eventually intersecting the section at $\rm T'$ (lower semicircle). Trajectories that begin just to the right of the blue trajectory will hit the section in a neighborhood of T. This separation is reflected by a discontinuity in the interval return map on the right.
Eventually, at lower $\rm \Delta [Ca]$, the tangency point $\rm T$ will flow directly onto the upper saddle $\rm SD$. At this stage, the positively sloped, orientation-preserving branches of the map are completely eliminated, except for the first branch adjacent to the saddle-focus. This corresponds to an absence of spiking trajectories that flow off the front of the dune; after spiking, all trajectories return to the dune from the back side of $\rm M_{PO}$, following the upper branch of the unstable manifold of the upper saddle.
Panels D1 and D2 of Fig.~\ref{fig21} display the configuration after the tangency point $\rm T$ reaches the upper saddle $\rm SD$. The orientation-reversing, negatively sloped branches begin to shrink and become flatter, eventually undergoing reverse period-doubling as the slope reaches $-1$ and the branches become stable.

\section{Discussion}

The key finding of our study is that the breadth of the triangular chaotic region at the intersection of bursting, tonic-spiking, and quiescent activity in the SiN model is primarily driven by the interaction between the global topological structure of the homoclinic framework surrounding a Shilnikov-Hopf bifurcation. Central to this understanding is the geometry of spiking activity, or the ``Swiss roll'' dynamics about the spiking manifold $\rm M_{PO}$, which we have modeled using multimodal, self-similar, one-dimensional interval return maps. Notably, these maps exhibit increasing steepness at their maxima as both the Shilnikov homoclinic and the subcritical Andronov-Hopf bifurcations are approached in the parameter space. This vertical stretching of the map in these regions can be interpreted as increased expansion across the two-dimensional slow-manifold surface---i.e., the dune $\rm M_Q$, containing the local unstable manifold of the saddle-focus---a common feature critical to understanding spiral chaotic dynamics.

In the classical theory of Shilnikov saddle-focus systems, the principal tool used to demonstrate the presence or absence of chaos is a map capturing distances from the 2D stable or unstable manifold of the saddle-focus upon first returns to a neighborhood of the saddle-focus equilibrium, sometimes referred to as the \emph{Shilnikov map}~\cite{scholarpedia,rossler2020,xing21,sharkov2024}.
In contrast, the one-dimensional interval maps in our study originate directly from first returns to a section of the 2D unstable manifold of the saddle-focus, or more precisely the dune $\rm M_Q$.
These maps do not measure distances from the dune $\rm M_Q$ as would the Shilnikov map, as the SiN system is highly dissipative; accordingly, the Shilnikov chaos is confined in parameter space to a small region about the $\rm homSF$ bifurcation curve, while the broader chaotic region is principally caused by Swiss-roll mixing of trajectories.

On the other hand, near the subcritical Andronov-Hopf bifurcation, the expansion associated with the saddle-focus weakens, resulting in a more coiled spiraling in the outward flow of trajectories on the dune $\rm M_Q$. As we move above the subcritical Andronov-Hopf  bifurcation curve AH$_{\rm sub}$, and enter the bursting region in the parameter space of the SiN model, the saddle-focus ``flings'' trajectories outward more vigorously, as the positive real part of its complex conjugate eigenvalues increases.

However, in this system, trajectories cannot spiral indefinitely due to the bounded nature of the $x$ variable, which is confined within the $[0, 1]$ range and possesses a stable branch in the slow-dune projection. This effectively creates a boundary that prevents trajectories from expanding infinitely downward in phase space. Trajectories approaching this stable branch of the $x$ nullcline encounter regions of contraction on the quiescent slow-motion manifold $\rm M_Q$. A strongly expanding saddle-focus directs more trajectories into this contracting region surrounding the lower stable branch of the $x$ nullcline, whereas a weak saddle-focus allows for more uniform expansion across $\rm M_Q$. 

Our findings suggest that neural systems with bounded gating variables, such as the slow $x$ variable in the SiN model, can support a large stable bursting region coexisting with a Shilnikov saddle-focus possessing a two-dimensional unstable manifold. If the outward spiraling on the unstable manifold of the saddle-focus (and thus on $\rm M_Q$) were to extend infinitely downward into negative $x$ values, and if the spiking manifold $\rm M_{PO}$ extended infinitely to the left into unrealistically negative $\rm [Ca]$ values, we would expect the boundaries of the chaotic region to expand outward, far from the Shilnikov-Hopf bifurcation. Investigating such a scenario may better clarify the role of this cod-2 bifurcation (S--SF), also referred to as the Belyakov type-II bifurcation, describing the transition between homoclinic saddle and saddle-focus, in the onset of chaos which is currently obscured by the dominance of stable bursting occurring in the SiN model.

Variations of this configuration are conceivable in both physically plausible models and mathematically or phenomenologically motivated systems. One possibility is that the beginning of the slow spiking manifold $\rm M_{PO}$ is accessible from the unstable manifold of the Shilnikov saddle-focus. Typically, such a slow two-dimensional stable manifold is initiated and terminated as a result of homoclinic bifurcations of a saddle or saddle-node, or through subcritical or supercritical Andronov-Hopf bifurcations, as well as saddle-node bifurcations of periodic orbits (PO) in the fast subsystem, any of which eventually give rise to a stable tonic-spiking PO. Note that a delayed loss of stability can substantially alter all routes to chaos that a system may undergo.

Another possible variation, following the earlier work on multi-rhythmic bursting~\cite{butera98}, involves weakening the contraction towards the spiking manifold $\rm M_{PO}$, possibly by slowing down the gating $x$ variable or by introducing appropriate bifurcation parameters into the slow subsystem of the SiN model. In one-dimensional maps, this would result in different heights among the arches of the multimodal return map. In the SiN model, however, the saddle-focus homoclinic bifurcations corresponding to different numbers of spikes may occur at parameter values up to or even exceeding numerical precision. In contrast, we anticipate that in the multi-rhythmic bursting model, the homoclinic bifurcations appear as distinct events, revealing intricate connections as the periodic orbits originating from different Shilnikov homoclinics interact. We expect that each Shilnikov homoclinic orbit will terminate through a unique Shilnikov-Hopf bifurcation in this scenario. Studying this configuration may further illuminate the homoclinic structure responsible for the chaos examined in this paper.

Furthermore, interactions between the saddle-focus and the upper saddle may form hetero-dimensional cycles at discrete points along the homoclinic Shilnikov saddle-focus curve (the homSF curve). The existence of such cycles guarantees the presence of blenders in some small neighborhoods, which could give rise to wild chaos, analogous to how Smale horseshoes give rise to classical chaos~\cite{LPbook17,Sciheritage,liblenders24}. However, several lines of evidence make us skeptical about the existence of wild chaos in the SiN model as currently studied. Stability windows appear to be densely packed throughout the parameter space, as shown in Figs.~\ref{fig1} and \ref{fig5}. Zooming in around the homSF curve does not reveal any obvious regions of persistent chaos. This may not be surprising since blenders and hyper-chaos are inherently higher-dimensional phenomena, and the Lyapunov dimension of the attractor in our system is just over two. The fact that the system can be effectively modeled with one-dimensional maps suggests that wild chaos is unlikely to exist in this context. Reducing the timescale separation in at least one variable might allow wild chaos to emerge near the Shilnikov homoclinic bifurcation and near non-transverse homoclinics to periodic orbits.

From biological and neuroscience perspectives, the presence of large chaotic regions could offer advantages for understanding evolution, development, and homeostasis through mathematical modeling. In a chaotic system, small changes in the parameters of a neural model lead to small changes in its statistical behavior, providing a smoother optimization landscape for feedback mechanisms and evolutionary processes. In contrast, a stable bursting system with narrow chaotic zones might respond to small parameter changes with either negligible qualitative change or drastic shifts, which could be less favorable in biological contexts. Therefore, the chaotic SiN model offers a useful framework for studying the role of chaos in evolution and development in future research.

\section{Conclusion}

In this study, we unveiled and examined the structure of a broad region of chaotic dynamics present in the given conductance-based neuronal model. By exploring the interplay between tonic spiking, saddle-focus equilibria, and saddle structures, we demonstrated how intricate homoclinic and heteroclinic connections give rise to widespread chaos. Our bifurcation analysis reveals that the chaotic behavior occupies a much larger parameter region than the narrow bands typically associated with mode transitions, such as spike-adding cascades.

Reducing the complex dynamics to a one-dimensional map allowed for a precise characterization of the observed chaos and its underlying mechanisms. Our findings highlight the robustness and generality of the SiN model within the chaotic region, which persists across a wide range of parameter values, challenging the conventional perception of neuronal models as being rigid. This broad chaotic region may have significant implications for understanding the variability and adaptability of neuronal behaviors, both in healthy and pathological states.

The templates utilized in this study serve to model and explain the underlying topology of the chaotic region, offering a structured approach to understanding the complex dynamics within the neuronal system. Unlike traditional methods that emphasize bifurcation analysis, these templates enable the application of kneading sequences, which are instrumental in characterizing the symbolic dynamics of the system. By providing a detailed representation of the topological structure, the templates allow for an accurate mapping of transitions and invariant sets within the chaotic regime. This approach not only facilitates the identification of chaos but also supports a deeper exploration of the dynamical behavior of the system across different parameter spaces, ensuring that the observed chaotic patterns are consistent with the broader theoretical framework of neuronal dynamics.

Future work will focus on extending these analyses to other neuronal models and exploring the potential functional roles of such widespread chaos in neural computation and information processing. The methods developed here, combining continuation analysis, symbolic dynamics, and Lyapunov exponents, provide a powerful toolkit for investigating complex dynamical systems beyond the context of neuronal modeling.

In conclusion, this work contributes to the growing body of knowledge on chaos in biological systems by providing new theoretical insights and suggesting directions for future exploration. The interplay between global topological structures and local bifurcations, as elucidated in this study, offers a rich avenue for further research, with potential applications that extend across multiple scientific disciplines.

\section{Acknowledgment}
 We thank the Brains \& Behavior initiative of Georgia State University for the graduate fellowships awarded to J. Scully and C. Hinsley. The authors thank the current and former members of the Shilnikov NeurDS (Neuro-Dynamical Systems) lab: Y. Keleta and Drs. H. Ju, D. Alacam, K. Pusuluri and J. Bourahmah for fruitful discussions. We are very grateful to Dr. A. Neiman for his help with Fig.~\ref{fig7}.  This work was partially funded by the National Science Foundation award DMS--2407999. 
 
\section{Data availability}
All data and codes are available from the first two authors upon request and openly available in GitHub at https://github.com/jamesjscully and  https://github.com/hinsley/PlantChaos.

\section{Appendices}
\subsection{Appendix~I: slow nullclines}\label{Ap1}
To locate the positions of the slow nullclines $x'=0$ and $[{\rm Ca}]'=0$ in the $([{\rm Ca}],\,x)$-phase plane
in the figures above, we first parametrize the membrane voltage $V_n$ within an equilibrium range $[-70;\,+20]$mV with some step size. Next, solve the equilibrium state equation $V' = 0$ for the $I_{KCa}$ current:  
$$
 I_{KCa} = -I_{I} - I_K - I_{T}  - I_{leak},  
$$ 
with the corresponding static functions $h_\infty(V_n)$, $n_\infty(V_n)$, and $m_\infty(V_n)$ for equilibria in the currents above, and next using Eq.~(\ref{IKCa}) to define the parameterized array   
$$
[{\rm Ca}](V_n)= 0.5\, \frac{I_{KCa}}{g_{KCa} (V_n-E_k)-I_{KCa} }.   
$$
The ordered pairs $ \left\{  [{\rm Ca}(V_n)], \, x_\infty (V_n) \right \} $ then populate the sought $x$ nullcline $x' = 0$.     
To determine the position of the nullcline $[{\rm Ca}]' = 0$,  first solve the equation $V' = 0$ for the $I_T$ current: 
$$
I_{T} = -I_{I} - I_K -  I_{KCa}  - I_{leak},  
$$ 
next find 
$$
x(V_n)=\frac{I_{T}}{g_T (E_{I}-V_n)},
$$
and then from Eq.~(\ref{eqCa}) find
$$
[{\rm Ca}](V_n) = K_{c}\,x(V_n)\,(E_{Ca}-V_n+\Delta [{\rm Ca}]),  
$$
and use the ordered pairs $\left \{ [{\rm Ca}(V_n)],\, x_\infty (V_n) \right \} $ to locate the  calcium nullcline $[{\rm Ca}]'=0$ in the slow $([{\rm Ca}],\,x)$-phase subspace.

\subsection{Appendix II: MATCONT continuation details}\label{Ap2}
For the numerical bifurcation analysis, we used the Matlab toolbox (version 2022b) MATCONT version 7p4 \cite{matcont,matcont1}, including the homotopy method to initialize connecting orbit continuation. Local bifurcation curves (saddle-node and Andronov-Hopf) were computed with standard settings. This identifies the cod-2 Bautin point (BP), as well as a Bogdanov-Takens (BT) bifurcation at $(\Delta \rm [Ca],\,\Delta V_x)=(-10.81,\,-10.24)$  The saddle-saddle curve SS and the Andronov-Hopf curve subAH meet at the BT-point. We could not connect the homoclinic curves in the region of interest to this BT-point, and hence do not display this. The dark blue curve $\rm homSN_{PO}$ corresponds to a saddle-node bifurcation of a sub-threshold periodic cycle. This branch is computed by branch-switching from the BP-point directly.\\
For the bifurcations of periodic orbits near $\Delta \rm [Ca]\approx -40$ we used a simulation at $\rm (\Delta [Ca],\, \Delta V_x) = (-60, -1.2)$ as initial data to start the continuation of periodic orbits in the $\rm \Delta [Ca]$-parameter. We then detected a first period-doubling bifurcation and several saddle-node bifurcations of cycles oscillating back and forth with increasing period, towards a homoclinic solution. Continuation of the first bifurcation yielded the light blue PD curve. The continuation of the limit cycles for increasing $\rm \Delta V_x$ all approached the Belyakov point, while in the other direction they split, either to the homoclinics with spiking or toward the BT point. The number of discretization points ntst was about 100-200 depending on how many spikes an orbit contained. Maximal step-size was increased to 10 instead of 0.1, while other settings were kept at default values.~\\
Next, we used the homotopy method to initialize the homoclinic orbits. We continued this branch in two parameters, where MATCONT detected another cod-2 Belyakov-II point, where the leading stable directions change (saddle -- saddle-focus transition), and finally the cod-2 Shilnikov-Hopf point or Belyakov-I one where the branch terminates. The orange curves are homoclinics to a saddle corresponding to a state in depolarization block (the red dot in Fig.~\ref{fig5}). The difference is in the number of spikes each orbit exhibits. Using that information, we were able to initialize the homoclinic continuation at $\rm \Delta V_x = 0.7$ where we ensured the final distance to the equilibrium was small (determined by two continuation variables $\varepsilon_1=0.1$, and $\varepsilon_0=0.01$), as that difference mostly corresponds to the potential.~\\
Finally, we noticed that the transition from the region with one spike to two spikes also involved a saddle-node bifurcation curve. Here the branch of periodic orbits originating from the homoclinic orbit becomes stable. For $\rm \Delta V_x = 0.7$, these bifurcations happen for nearly identical parameter values. Continuation for decreasing $\rm \Delta V_x$ shows that near the AH curve they split. The connecting orbits cross the AH curve and get closer and then appear to get closer to the $\rm AH_{sub}$ branch. Instead, the saddle-node curves remained above the AH curve. Near the BP-point they have a swallow-tail structure, and then get closer to the Belyakov II-type S--SF point. We could not reach the S--SF point with continuation. Not all connecting orbits could be continued towards the lower left corner of the parameter region due to convergence errors. We recomputed some curves with a collocation number ntst=300 of points on each and smaller step-sizes, but as each curve takes long several hours in MATCONT, we did not pursue this further.

\subsection{Appendix III: the cellular SIN-model description}\label{Ap3}

The dynamics of the membrane potential, $V$, are governed by the following equation:
\begin{equation}
 	C_{m} {V}^\prime = -I_{I} - I_K - I_{T}  - I_{KCa} - I_{leak}.  \label{a1}
 \end{equation}
 The fast inward sodium and calcium current $I_I$ is given by   
\begin{equation}
I_{I}=g_{I}\,h\, m^{3}_{\infty}(V)(V-E_{I}), 
 \end{equation}
with the reversal potential $E_{I}=30$mV and the maximal conductance value $g_{I}=4$nS, and  
\begin{align}
 m_{\infty}(V) &= \frac{\alpha_{m}(V)}{\alpha_{m}(V) +\beta_{m}(V)}, \qquad  \\
 \alpha_{m}(V) &= 0.1 \frac{50-V_s}{-1+e^{(50-V_s)/10}}\ , \\ 
 \quad \beta_{m}(V) &= 4 e^{(25-V_s)/18}  
\end{align}
 the dynamics of its inactivation  gating variable $h$ is given by      
\begin{align}
 {h}^\prime &=  \frac{h_{\infty}(V)-h}{ {\tau_{h}(V)} }, \quad \mbox{where} \quad \\
h_{\infty}(V) &= \frac{\alpha_{h}(V)}{\alpha_{h}(V) +\beta_{h}(V)}\quad \mbox{and} \\
 \quad \tau_{h}(V) &= \frac{12.5}{\alpha_{h}(V) +\beta_{h}(V)} 
\end{align}
 \begin{align}
 \alpha_{h}(V) &= 0.07 e^{(25-V_s)/20}, \\
 \beta_{h}(V) &= \frac{1}{1+e^{(55-V_s)/10}}, \quad \mbox{and} \quad \\
 V_s &= \frac{127V+8265}{105}{\rm mV}.
 \end{align} 
  The fast potassium current $I_K$ is given by the equation  
\begin{equation}
I_K=g_{K} n^{4}(V-E_{K}),  
\end{equation}
with the reversal potential $E_K=-75$mV and the maximal conductance set as $g_K=0.3$nS. The dynamics of inactivation gating variable are described by 
\begin{align}
{n}^\prime &= \frac{n_{\infty}(V)-n}{\tau_{n}(V)}, \quad \mbox{with} \quad \\
n_{\infty}(V) &= \frac{\alpha_{n}(V)}{\alpha_{n}(V) +\beta_{n}(V)} \quad \mbox{and} \quad \\
\tau_{n}(V) &= \frac{12.5}{\alpha_{n}(V) +\beta_{n}(V)}, \quad \mbox{where} \quad \\
\alpha_{n}(V) &= 0.01 \frac{55-V_s}{e^{(55-V_s)/10}-1} \quad \mbox{} \quad \\
\beta_{n}(V) &= 0.125 e^{(45-V_s)/80}.
\end{align}	
The leak current is given by 
\begin{align}
I_{\rm leak} &=g_L (V-E_L), 
\end{align}
with $E_L =-40{\rm mV}$, and $g_L = 0.003{\rm nS}$. The sub-group of slow currents in the model includes the TTX-resistant sodium and calcium current $I_T$ given by \begin{align}
 I_T &= g_T x (V-E_I),
\end{align}
with $E_I =30{\rm mV}$ and $g_T = 0.01{\rm nS}$, while the dynamics of its slow activation variable are described by
 \begin{align}
{x}^\prime &= \frac{x_{\infty}(V)-x}{\tau_{x}}, \quad \mbox{where} \\ 
\quad x_{\infty} &= \frac{1}{1+e^{-0.15(V+50-\Delta V_x)}},
\end{align}\label{eqttx}
and the time constant $\tau_{x}$ set as 100 or 235ms in this study. The slowest outward $Ca^{2+}$ activated $K^+$ current given by  
\begin{equation}
I_{KCa} =g_{KCa}\frac{[Ca]_i}{0.5+[Ca]_i}(V-E_{K})
 \label{IKCa}
\end{equation}
with  $E_K = -75{\rm mV}$, and $g_{KCa}=0.03$, while the dynamics of the intracellular calcium concentration is governed by 
\begin{equation}
{[{\rm Ca}]}^\prime = \rho \left ( K_{c}\, x\, (E_{\rm Ca} - V + \Delta [{\rm Ca}])-[{\rm Ca}] \right  )  \label{eqCa}
\end{equation}
with the Nernst reversal potential  $E_{Ca}=140$mV, and small constants $\rho= 0.0003{\rm ms}^{-1}$ and  $K_c=0.0085{\rm mV}^{-1}$.

\subsection{Appendix IV: SSCS encoder algorithm}

Much of our analysis of the SiN model involves performing large-scale scans of the parameter space, numerically integrating one or more trajectories at every sampled parameter value.
Because numerical integration of the SiN model equations produces large amounts of data in the form of trajectory time series, we designed an algorithm to encode trajectories as SSCSs during integration.

This encoding scheme permits us to discard trajectory data when we wish to study only the symbolic dynamics of the template associated with the SiN model, so it confers a large improvement in memory efficiency for symbolic scans.
The encoding algorithm involves keeping track of the relative sequential ordering of the occurrences of certain local maxima observed in $V(t)$ and $V'(t)$ time series.
These maximum events are denoted as follows:
\begin{itemize}
	\item $I$: A local maximum in $V'$.
	\item $V^-$: A local maximum in $V$ not exceeding the voltage value $V_{\rm SD}$ of the upper saddle equilibrium point.
	\item $V^+$: A local maximum in $V$ exceeding $V_{\rm SD}$.
\end{itemize}
In general, the value $V^+$ corresponds to the maximum voltage value attained during a spike, while $V^-$ is used to detect the refractory rebound after a spike-train burst.
A DONE event is used to signify that the integration has terminated.
The SSCS encoding is carried out as described in Algorithm~\ref{alg:sscsencoder}.

\begin{algorithm}\label{alg:sscsencoder}
\DontPrintSemicolon
\caption{Signed Spike-Count Sequence Encoder}
\KwIn{A function \FuncSty{NextEvent} which returns the next event detected by the integrator each time it is called}
\KwOut{A list SSCS of integers}
Spikes $\gets 0$\;
Event $\gets$ NULL\;
PrevEvent $\gets$ NULL\;
PrevPrevEvent $\gets$ NULL\;
SSCS $\gets []$\;
\While{Event $\neq$ DONE}{
	Event $\gets$ \FuncSty{NextEvent()}\;
	\eIf{Event $= V^-$}
	{
		\eIf{PrevPrevEvent $= V^+$}
		{
			SSCS\FuncSty{.append(}$-$Spikes\FuncSty{)}\;
		}{
			SSCS\FuncSty{.append(}Spikes\FuncSty{)}\;
		}
		Spikes $\gets 0$\;
	}{
		\If{Event $= V^+$}
		{
			Spikes $\gets$ Spikes$+ 1$\;
		}
	}
	PrevPrevEvent $\gets$ PrevEvent\;
	PrevEvent $\gets$ Event\;
}
\end{algorithm}

Because the events used in this algorithm only depend on the local maxima of $V$ and $V'$ (the latter of which is an explicit function of the current position in state space so that it can be computed directly from the governing equations of the SiN model), it is inexpensive to detect these events with very high accuracy.
In practice, we use the VectorContinuousCallback event callbacks from DifferentialEquations.jl to detect these events and trigger the corresponding subroutines rather than a blocking NextEvent function and while loop; however, the basic encoding algorithm is the same as described here.

One caveat of Algorithm~\ref{alg:sscsencoder} is that it may, for certain choices of system parameters, misclassify bursting trajectories which return to the dune $\rm M_Q$ extremely closely to the tangency point $\rm T$; in this case, one of the nonzero signed spike counts returned in the SSCS list would have incorrect sign.
However, it is remarkably accurate and in most cases---especially near the homSF bifurcation line in the chaotic parameter region---this weakness is exceedingly rare and may be disregarded.
We have not observed any parameter value for which this encoding error occurs with appreciable severity.

\section*{References}
\bibliographystyle{unsrt}

\end{document}